\documentclass[reqno,a4paper]{amsart}
\usepackage{amsmath, amssymb, amscd}

\numberwithin{equation}{section}
\theoremstyle{plain}
\newtheorem{thm}{Theorem}[section]
\newtheorem{cor}[thm]{Corollary}

\newtheorem{lem}[thm]{Lemma}
\newtheorem{prop}[thm]{Proposition}
\newtheorem{Def}[thm]{Definition}

\theoremstyle{remark}
\newtheorem{rem}[thm]{Remark}

\title{Difference Fourier transforms for nonreduced root systems}
\author{Jasper V. Stokman}
\address{Jasper V. Stokman,
KdV Institute for Mathematics, Universiteit van Amsterdam,
Plantage Muidergracht 24, 1018 TV Amsterdam, The Netherlands.}
\email{jstokman@science.uva.nl}
\date{November 19th, 2001}
\begin{document}

\begin{abstract}
In the first part of the paper kernels
are constructed which meromorphically
extend the Macdonald-Koornwinder polynomials in their degrees.
In the second part of the paper the kernels associated with 
rank one root systems are used
to define nonsymmetric variants of the spherical Fourier transform
on the quantum $\hbox{SU}(1,1)$ group. Related Plancherel and 
inversion
formulas are derived using double affine Hecke algebra techniques. 
\end{abstract}
\maketitle
\noindent
{\bf Contents.}
\begin{enumerate}
\item[{\S 1}] Introduction.
\item[{\S 2}] The double affine Hecke algebra.
\item[{\S 3}] General remarks on difference Fourier transforms.
\item[{\S 4}] The Macdonald-Koornwinder transform.
\item[{\S 5}] The construction of the Cherednik kernels.
\item[{\S 6}] Properties of the Cherednik kernels.
\item[{\S 7}] An extension of the Macdonald-Koornwinder transform.
\item[{\S 8}] The (non)symmetric Askey-Wilson function transform.
\item[{\S 9}] Appendix.
\end{enumerate}

\section{Introduction}
The Weyl algebra of differential operators with polynomial
coefficients on Euclidean $n$-space is the image of the Heisenberg
algebra under the Schr{\"o}dinger representation. 
The classical Fourier transform  
induces an automorphism of the Weyl algebra which
interchanges the role of polynomial multiplication operators 
and constant coefficient differential operators. This automorphism
(and its generalizations) is called the {\it duality isomorphism}
in the present paper. 
Plancherel and inversion formulas for the classical Fourier transform
can be easily derived from the above observations by first proving
algebraic versions on the cyclic module of the Weyl algebra
generated by the Gaussian.
See the survey paper \cite{Ho} for more details and references.
We generalize this approach to difference Fourier transforms
associated to nonreduced root systems.

Work of Dunkl \cite{D}, \cite{D1}, Heckman \cite{H}, 
Opdam \cite{O1}, de Jeu \cite{dJ}, Cherednik \cite{C-1}, \cite{C0},
\cite{C4}, Macdonald \cite{M1}, Noumi \cite{N2}, \cite{N} 
and others have led to generalizations 
of the Weyl algebra and the 
underlying Heisenberg algebra, which are naturally
associated to Fourier transforms arising from harmonic analysis on 
Cartan motion groups, Riemannian symmetric spaces and 
compact quantum Riemannian spaces. 
In each case, the Weyl algebra is replaced by an algebra
consisting of differential (or difference) reflection operators and 
multiplication
operators. This generalized Weyl algebra 
may be realized as the faithful image of (an
appropriate degeneration of) the double affine Hecke algebra 
under Cherednik's representation
(the analogue of the Schr{\"o}dinger representation).
In each case, the generalized Weyl algebras form a powerful tool in 
the study of the related Fourier transforms.
The main goal of the present paper is to extend this picture to
include the case of 
Fourier transforms arising from harmonic analysis
on the simplest quantum analogue of a {\it noncompact} simple
Lie-group, namely the quantum $\hbox{SU}(1,1)$ group.

Some remarks are here in order on the theory of locally compact quantum
groups. Despite the fact that {\it compact} quantum groups are
well understood, also from the viewpoint of harmonic analysis
(see e.g. \cite{N}),
this is by far not the case for {\it noncompact} semisimple quantum
groups.
Recent developments though essentially
settled the theory for the quantum $\hbox{SU}(1,1)$ group.
A quantization of the normalizer of 
$\hbox{SU}(1,1)$ in $\hbox{SL}(2,\mathbb{C})$ 
was constructed as a locally compact quantum group by Koelink and
Kustermans \cite{KK} (see Kustermans and Vaes \cite{KV} for a detailed 
account on the general theory of locally compact
quantum groups). The harmonic analytic aspects for the 
quantum $\hbox{SU}(1,1)$ group were analyzed in 
\cite{Kak}, \cite{Kak2} and \cite{KS1}.
The harmonic analysis on the quantum $\hbox{SU}(1,1)$ group
in \cite{KS1} led to an explicit Fourier transform, whose
Plancherel and inversion formula were derived by classical
function-theoretic methods in \cite{KS2}. The transform is called the
Askey-Wilson function transform, since its kernel forms a meromorphic
continuation of the well-known Askey-Wilson polynomials (see
\cite{AW}) in their degrees.

To explore the role of the double affine Hecke algebra
for the Askey-Wilson function transform, we first need to construct
nonsymmetric variants of the transform which
induce the proper analogue of the duality isomorphism on
the double affine
Hecke algebra. The construction of the associated kernel is developed in
this paper for nonreduced root systems of arbitrary rank.
Two basic features of the kernel are its explicit series expansion
in terms of Macdonald-Koornwinder polynomials, and the fact that
it is a meromorphic continuation 
of the Macdonald-Koornwinder polynomials
in their degrees. These results are inspired by 
Cherednik's paper \cite{C1}, in which
such kernels were introduced for reduced root systems.

Restricting attention to rank one, we use these kernels to 
define nonsymmetric variants of the Askey-Wilson function transform
which induces the analogue of the duality isomorphism
on the double affine Hecke algebra. Instead of defining the transforms
on compactly supported
functions (as was done in the classical approach \cite{KS2}), we
define it now on a space of function consisting of a direct
sum of {\it two} cyclic modules of the double affine Hecke algebra. 
{}From harmonic analytic point of view a second cyclic module
is needed to take care of the ``strange part'' of the support of 
the Plancherel measure, i.e. the part of the support
coming from contributions of unitary representations of the 
quantized universal enveloping algebra $\hbox{U}_q(\hbox{su}(1,1))$
which vanish in the classical limit.

Explicit evaluations 
of the images of the cyclic vectors under the nonsymmetric
Askey-Wilson function transform then suffice to prove
algebraic inversion and Plancherel type formulas for 
the nonsymmetric Askey-Wilson function transform.
The image of the cyclic vector of
the first, ``classical'' module follows from
an extension of the polynomial Macdonald-Koornwinder 
theory. This in particular entails explicit formulas 
for the Macdonald-Koornwinder type constant term of the product
of the inverse of a generalized Gaussian and two Macdonald-Koornwinder
polynomials (see \cite{C1} for the analogous results
in the reduced set-up). 
The computation of the image of the cyclic vector of the
second, ``strange'' module (see Proposition \ref{FourierGaussian})
is a key result which combines many of the properties of the 
kernels with some elementary elliptic function theory. 
The particularly large number of parameter freedoms, caused by
the fact that we are dealing with {\it nonreduced} root systems,
is also used here in an essential way.
This may be seen as an extra justification for the special
attention to nonreduced root systems in the present paper.

Finally we show that the results described in the 
previous paragraph easily lead to new 
proofs of the Plancherel and inversion formula for the 
symmetric Askey-Wilson function transform (see
\cite{KS2} and \cite{St3} for the classical function-theoretic approach).

Considering difference analogues of Harish-Chandra transforms
only from the viewpoint of double affine Hecke algebras lead
to several other self-dual difference Fourier transforms,
see e.g. \cite{Cnew} and \cite{CO}.
This flexibility in choices, combined with the lack in 
comprehension of noncompact
semisimple quantum groups, thus poses essential problems in deciding
which difference analogues of the Harish-Chandra transform
arise from harmonic analysis on noncompact quantum symmetric spaces.
The present, detailed study of the Askey-Wilson function transform from the
viewpoint of double affine Hecke algebras hopefully
provides some new insights on this matter.

The paper is organized as follows. In Section 2 we introduce 
the double affine Hecke algebra, the duality isomorphism and the
analogue of the Gaussian.
In Section 3 we introduce the concept of difference Fourier
transforms. We furthermore explain what the main techniques
are going to be in the study of such transforms. 
In Section 4 we show how the polynomial
Macdonald-Koornwinder transform (introduced in \cite{St1}) fits into
this general scheme. In Section 5 we introduce
kernels with which explicit integral transforms can be constructed
that fit into the concept of difference Fourier transforms. We call the
kernels {\it Cherednik kernels}, since they generalize
kernels introduced by Cherednik in \cite{C1}. In Section 6 we 
study the main properties of the Cherednik kernels. In particular, 
we show that the Cherednik kernels 
meromorphically extend the Macdonald-Koornwinder polynomials
in their degrees, and that they satisfy a natural duality property
which extends the duality of the Macdonald-Koornwinder polynomials. 
In Section 7 we study an extension of the Macdonald-Koornwinder
transform, acting on the cyclic module generated
by the inverse of the Gaussian. It forms the second example of
a difference Fourier transform in the sense of Section 3.
In Section 8 we restrict attention to rank one. We define 
the nonsymmetric Askey-Wilson function transform as an integral
transform with kernel given by the (rank one) Cherednik kernel. We 
show that it also qualifies as a difference Fourier transform in the
sense of Section 3. We prove inversion and Plancherel Theorems
on the algebraic and on the $L^2$-level. In the Appendix we finally
discuss certain bounds for Macdonald-Koornwinder polynomials which are
needed for the construction of the Cherednik kernels.

{\it Acknowledgements:} The author is
supported by a fellowship from the Royal
Netherlands Academy of Arts and Sciences (KNAW).


\section{The double affine Hecke algebra}


\subsection{The affine root system of type $C^\vee C_n$}
Denote by $\epsilon_i$ ($i=1,\ldots,n$) the standard orthonormal
basis of the $n$-dimensional Euclidean space $(\mathbb{R}^n,\langle 
\cdot,\cdot\rangle)$. Let
\[\Sigma=\{\pm\epsilon_i\pm\epsilon_j\}_{1\leq i<j\leq n}\cup
\{\pm 2\epsilon_i\}_{i=1}^n\subset \mathbb{R}^n
\]
be the root system of type $C_n$.
Let $\hbox{Aff}(\mathbb{R}^n)$ be the space of
affine linear transformations $f: \mathbb{R}^n\rightarrow \mathbb{R}$.
As a vector space, $\hbox{Aff}(\mathbb{R}^n)\simeq \mathbb{R}^n\oplus
\mathbb{R}\delta$ via the formula
\[ (v+\lambda\delta)(w)=\langle v,w\rangle+\lambda,\qquad v,w\in
\mathbb{R}^n, \quad\lambda\in\mathbb{R}.
\]
We extend the form  $\langle \cdot,\cdot \rangle$ 
to a positive semi-definite form
on $\hbox{Aff}(\mathbb{R}^n)$ by requiring the constant function
$\delta$ to be in the radical of $\langle \cdot,\cdot \rangle$. Then 
\[R=\Sigma+\mathbb{Z}\delta\subset \hbox{Aff}(\mathbb{R}^n)\]
is the {\it affine root system} of type $\widetilde{C}_n$.
Associated with $f\in R$, we have the reflection
$r_f\in \hbox{Gl}_{\mathbb{R}}(\hbox{Aff}(\mathbb{R}^n))$ defined by
\[r_f(g)=g-2\frac{\langle f,g\rangle}{\langle f,f\rangle}f,\qquad
g\in \hbox{Aff}(\mathbb{R}^n).
\]
The {\it affine Weyl group} $\mathcal{W}$ is the
sub-group of $\hbox{Gl}(\hbox{Aff}(\mathbb{R}^n))$
generated by all the reflections $r_f$ ($f\in R$).

There are two important descriptions of $\mathcal{W}$,
namely as a Coxeter group, and as a semi-direct
product of a finite reflection group with a lattice.
For the presentation of $\mathcal{W}$ as a Coxeter group, we
choose the standard basis $\{a_0,a_1,\ldots,a_n\}$
of the affine root system $R$ by
\[a_0=\delta-2\epsilon_1,\quad
a_i=\epsilon_i-\epsilon_{i+1},\quad a_n=2\epsilon_n
\]
for $i=1,\ldots,n-1$, and we set $r_i=r_{a_i}$ for the associated
simple reflections.
The above choice of basis induce a
decomposition of $\Sigma$ and $R$ in positive roots and negative roots,
the positive roots being given by
\[ \Sigma_+=\{\epsilon_i\pm\epsilon_j\}_{i<j}\cup\{2\epsilon_i\}_i,
\qquad R_+=\Sigma_+\cup \{f\in R \, | \, f(0)>0\}.
\]
Furthermore, it is well known that the 
affine Weyl group $\mathcal{W}$ is generated by the simple
reflections $r_i$ ($i=0,\ldots,n$). The fundamental relations
between the simple reflections $r_i$ are given by the Coxeter
relations
\begin{equation*}
\begin{split}
r_ir_{i+1}r_ir_{i+1}&=r_{i+1}r_ir_{i+1}r_i,\qquad i=0,i=n-1,\\
r_ir_{i+1}r_i&=r_{i+1}r_ir_{i+1},\qquad i=1,\ldots,n-2,\\
r_ir_j&=r_jr_i,\qquad |i-j|\geq 2,\\
r_i^2&=1,\qquad i=0,\ldots,n.
\end{split}
\end{equation*}

Let $W_0\subset \mathcal{W}$ be the sub-group of 
$\mathcal{W}$ generated by $r_1,\ldots,r_n$. Then
$W_0$ is the Weyl group of the root system $\Sigma$, hence
isomorphic to $S_n\ltimes (\pm 1)^n$ with $S_n$ the 
symmetric group in $n$ letters. Let $\Lambda$ be the $W_0$-invariant
$\mathbb{Z}$-lattice of $\mathbb{R}^n$ with basis $\{\epsilon_i\}_i$,
\[ \Lambda=\bigoplus_{i=1}^n\mathbb{Z}\epsilon_i.
\]
The lattice $\Lambda$ can be naturally identified with the co-root lattice
as well as with the weight lattice of the root system $\Sigma$.
The second description of the affine Weyl group $\mathcal{W}$
is given by
\[\mathcal{W}\simeq W_0\ltimes \Lambda
\]
with the lattice elements $\lambda\in\Lambda$ acting on 
$\hbox{Aff}(\mathbb{R}^n)$ by the formula
\[ \lambda(f)=f+\langle f,\lambda\rangle\delta,\qquad
\lambda\in\Lambda,\,\,f\in \hbox{Aff}(\mathbb{R}^n).
\]
Finally we note that the set
\[ R_{nr}=R\cup \mathcal{W}\frac{a_0}{2}\cup \mathcal{W}\frac{a_n}{2}
\subset \hbox{Aff}(\mathbb{R}^n)
\]
is a (nonreduced) affine root system having $\mathcal{W}$
as its associated affine Weyl group. 
In Macdonald's \cite{M} terminology, $R_{nr}$ is the {\it affine
root system of type} $C^\vee C_n$. 


\subsection{Difference multiplicity functions}

The nonreduced affine root system $R_{nr}$ of type $C^\vee C_n$
has five $\mathcal{W}$-orbits, namely
\[\mathcal{W}a_0,\qquad\mathcal{W}\frac{a_0}{2},\qquad
\mathcal{W}a_i,\qquad
\mathcal{W}a_n,\qquad\mathcal{W}\frac{a_n}{2},
\]
where $i$ may be arbitrarily chosen from the index set
$\{1,\ldots,n-1\}$.

A $\mathcal{W}$-invariant complex-valued function on $R_{nr}$ is
called a {\it multiplicity function}. We denote a multiplicity
function by $\mathbf{t}=\{t_f\}_{f\in R_{nr}}$, with
$t_f\in\mathbb{C}$ its value at the root $f\in R_{nr}$, and 
we set $t_f=1$ when $f\in\hbox{Aff}(\mathbb{R}^n)\setminus R_{nr}$.

In the theory presented in the next (sub-)sections, 
a sixth generic parameter
$q^{\frac{1}{2}}\in \mathbb{C}^\times$ appears. Its square $q$
is a parameter which arises
from the realization of the group algebra
$\mathbb{C}[\mathcal{W}]$ as the algebra of difference reflection
operators with constant coefficients. 
The pair $\alpha=(\mathbf{t},q^{\frac{1}{2}})$, where 
$\mathbf{t}$ is a
multiplicity function and $q^{\frac{1}{2}}\in \mathbb{C}^\times$, 
is called a {\it difference multiplicity function}. 
Alternatively, a difference multiplicity function
$\alpha$ can be represented by an ordered six-tuple
\[
\alpha=(\mathbf{t},q^{\frac{1}{2}})=(t_0,u_0,t_n,u_n,t,q^{\frac{1}{2}}),
\]
where
\[
t_0=t_{a_0},\qquad u_0=t_{a_0/2},\qquad
t=t_i=t_{a_i},\qquad
t_n=t_{a_n},\qquad
u_n=t_{a_n/2},
\]
where  $i$ may be any integer from the index set $\{1,\ldots,n-1\}$.

We define two ``involutions'' $\sigma$ and $\tau$ on 
difference multiplicity functions, both
preserving the parameter $q^{\frac{1}{2}}$.
The involution
$\sigma$ acts on $\alpha=(\mathbf{t},q^{\frac{1}{2}})$ 
by interchanging the value of $\mathbf{t}$ on the 
$\mathcal{W}a_0$-orbit with its value
on $\mathcal{W}a_n/2$, while $\tau$ acts
by interchanging the value of $\mathbf{t}$ on the
$\mathcal{W}a_0$-orbit with its value on $\mathcal{W}a_0/2$.
The new difference multiplicity functions are denoted by
$\alpha_\sigma$ and $\alpha_\tau$, respectively.

Associated with a generic difference multiplicity function
$\alpha=(\mathbf{t},q^{\frac{1}{2}})$, we write
$\alpha_\ddagger=(\mathbf{t}^{-1},q^{-\frac{1}{2}})$ 
with $\mathbf{t}^{-1}=\{t_f^{-1}\}_{f\in R_{nr}}$ for the difference
multiplicity function with inverted parameters.

If an object $H$ depends on a difference multiplicity function
$\alpha=(\mathbf{t},q^{\frac{1}{2}})$, then 
e.g. $H_{\ddagger\sigma\tau}$ (or $H^{\ddagger\sigma\tau}$)
denotes the object $H$ depending on the
difference multiplicity function 
$\alpha_{\ddagger\sigma\tau}=((\alpha_\ddagger)_\sigma)_\tau$.
For example, if
\[H=H_{(t_0,u_0,t_n,u_n,t,q^{\frac{1}{2}})},
\]
then
\[H_{\ddagger\sigma\tau}=
H_{(u_0^{-1},u_n^{-1},t_n^{-1},t_0^{-1},t^{-1},q^{-\frac{1}{2}})}.
\]
Throughout the remainder of 
this section we fix a generic difference multiplicity
function $\alpha=(\mathbf{t},q^{\frac{1}{2}})$. 


\subsection{Difference reflection operators}

Let $\mathcal{A}=\mathbb{C}[x^{\pm 1}]$
be the algebra of Laurent polynomials in $n$ indeterminates
$x=(x_1,\ldots,x_n)$, $\mathbb{C}(x)$ the
quotient field of $\mathcal{A}$ and $\mathcal{O}=
\mathcal{O}\bigl((\mathbb{C}^\times)^n\bigr)$ the ring of analytic
functions on the complex torus $(\mathbb{C}^\times)^n$. We write 
$\mathcal{M}=\mathcal{M}\bigl((\mathbb{C}^\times)^n\bigr)$ for the field of 
meromorphic functions on $(\mathbb{C}^\times)^n$. 
Since $(\mathbb{C}^\times)^n$ is a connected 
domain of holomorphy, $\mathcal{M}$ is
isomorphic to the quotient field of $\mathcal{O}$ (see e.g. 
\cite[Thm. 7.4.6]{Ho}).
We have natural inclusions
\[\mathcal{A}\subset \mathbb{C}(x)\subset \mathcal{M},\qquad\quad
\mathcal{A}\subset \mathcal{O}\subset \mathcal{M}.
\]
Let $\mathbb{C}[\mathcal{W}]$ be the group algebra of the affine Weyl
group $\mathcal{W}$ over $\mathbb{C}$.
\begin{Def}\label{qdiff0}
The algebra $\mathcal{D}_{q}$ of $q$-difference reflection
operators with meromorphic coefficients is defined by
\[
\mathcal{D}_q=\mathcal{M}\otimes_{\mathbb{C}}{\mathbb{C}}[\mathcal{W}]
\]
with multiplication
\[ (g\otimes v)(h\otimes w)=g(vh)\otimes vw,\qquad g,h\in
\mathcal{M},\,\,\, w,v\in \mathcal{W},
\]
where $\mathcal{W}\simeq W_0\ltimes \Lambda$
acts as field-automorphisms on $\mathcal{M}$ by the
formulas
\begin{equation}\label{Wfund}
\begin{split}
\bigl(\lambda h\bigr)(x)&=
h(q^{\lambda_1}x_1,q^{\lambda_2}x_2,\ldots,q^{\lambda_n}x_n),\\
\bigl(r_ih\bigr)(x)&=h(x_1,\ldots,x_{i-1},x_{i+1},x_i,x_{i+2},\ldots,x_n),\\
\bigl(r_nh\bigr)(x)&=h(x_1,\ldots,x_{n-1},x_n^{-1}),
\end{split}
\end{equation}
for $h(x)=h(x_1,\ldots,x_n)\in\mathcal{M}$, 
$\lambda=\sum_j\lambda_j\epsilon_j\in\Lambda$ and for $i=1,\ldots,n-1$.
\end{Def}
Note that the action \eqref{Wfund} of $\mathcal{W}$ depends on the
sixth parameter $q^{\frac{1}{2}}$ of the underlying difference
multiplicity function $\alpha$. When confusion can arise on
the underlying difference multiplicity function, we say that 
$\mathcal{W}$ acts by constant coefficient $q$-difference
reflection operators when the action is given by \eqref{Wfund} (i.e.
when the sixth parameter of the underlying difference multiplicity
function is $q^{\frac{1}{2}}$).

The algebra $\mathcal{D}_q$ acts on $\mathcal{M}$ by
\begin{equation}\label{natural}
\bigl(\bigl(gw\bigr)h\bigr)(x)=g(x)\bigl(wh\bigr)(x),
\qquad g\in\mathcal{M},\,\, w\in\mathcal{W}
\end{equation}
for $h\in\mathcal{M}$,
with the action of $\mathcal{W}$ on $\mathcal{M}$ as given
by \eqref{Wfund}. Here we use the notation $g(x)w$ or $gw$ for a pure
tensor $g\otimes w\in \mathcal{D}_q$.
We use the terminology that $\mathcal{D}_q$
(or a sub-algebra of $\mathcal{D}_q$) acts as $q$-{\it difference
reflection operators} on some function space on the complex torus 
$(\mathbb{C}^\times)^n$ when the action is given by formula
\eqref{natural}, with the action of $\mathcal{W}$ on functions $h$
defined by \eqref{Wfund}.

We define monomials by
\[ x^f=q^{\frac{m}{2}}x_1^{\lambda_1}x_2^{\lambda_2}\cdots
x_n^{\lambda_n}\in\mathcal{A}
\]
for any $f=\lambda+\frac{1}{2}m\delta\in\Lambda+\frac{1}{2}\mathbb{Z}\delta$,
where $\lambda=\sum_i\lambda_i\epsilon_i$. Again note here that
this definition depends on the sixth parameter $q^{\frac{1}{2}}$
of the underlying difference multiplicity function $\alpha$, 
but it will always be clear from the context what the underlying
difference multiplicity function is. 
Associated to any root $f\in R$ and any difference multiplicity
function $\alpha=(\mathbf{t},q^{\frac{1}{2}})$, we construct 
now explicit difference reflection 
operators $T_f=T_f^\alpha\in\mathcal{D}_q$ as
follows.
\begin{Def}\label{qdiff}
For $f\in R$, the
difference reflection operators $T_f=T_f^{\alpha}\in\mathcal{D}_q$ 
is defined by  
\[
T_f=t_f+t_f^{-1}c_f(x)(r_f-\hbox{id}),
\]
with coefficient $c_f(x)=c_f^\alpha(x)\in \mathbb{C}(x)$ given by
\begin{equation}\label{cf1}
c_f(x)=\frac{(1-t_ft_{f/2}x^{f/2})(1+t_ft_{f/2}^{-1}x^{f/2})}
{(1-x^f)}
\end{equation}
when $f/2\in R_{nr}$ and
\begin{equation}\label{cf2}
c_f(x)=\frac{(1-t_f^2x^f)}{(1-x^f)}
\end{equation}
when $f/2\not\in R_{nr}$.
\end{Def}
Note that the formula \eqref{cf2} for $c_f$ in case that $f/2\not\in R_{nr}$
can formally be written as \eqref{cf1} 
in view of the convention that $t_g=1$ when 
$g\in \hbox{Aff}(\mathbb{R}^n)\setminus R_{nr}$. 

The following crucial theorem, due to Cherednik in the reduced set-up,
was proven by Noumi \cite{N}. 

\begin{thm}\label{basic}
The $q$-difference reflection operators 
$T_j=T_j^\alpha\in \mathcal{D}_q$ defined
by $T_j=T_{a_j}$
\textup{(}$j=0,\ldots,n$\textup{)}, satisfy the fundamental commutation
relations of the affine Hecke algebra of type $\widetilde{C}_n$.
In other words, they satisfy the braid relations
\begin{equation}\label{braid}
\begin{split}
T_iT_{i+1}T_iT_{i+1}&=T_{i+1}T_iT_{i+1}T_i,\qquad i=0,i=n-1,\\
T_iT_{i+1}T_i&=T_{i+1}T_iT_{i+1},\quad\qquad i=1,\ldots,n-2,\\
T_iT_j&=T_jT_i,\qquad\qquad\qquad |i-j|\geq 2,
\end{split}
\end{equation}
and the quadratic relations
\begin{equation}\label{quadratic}
(T_j-t_j)(T_j+t_j^{-1})=0,\qquad\qquad j=0,\ldots,n,
\end{equation}
in the algebra $\mathcal{D}_q$.
\end{thm}

Following Noumi \cite{N}, one can now define the $Y$-operators
$Y_i=Y_i^{\alpha}\in\mathcal{D}_q$ for $i=1,\ldots,n$ by
\begin{equation}
Y_i=T_i\cdots T_{n-1}T_n
T_{n-1}\cdots T_1
T_0T_1^{-1}\cdots T_{i-1}^{-1}.
\end{equation}
The operator $Y_i$ is naturally associated to the lattice element
$\epsilon_i\in\Lambda$ considered as element of the affine Weyl group
$\mathcal{W}\simeq W_0\ltimes \Lambda$, since 
\[ \epsilon_i=r_i\cdots r_{n-1}r_nr_{n-1}\cdots r_1r_0r_1\cdots
r_{i-1}\in\mathcal{W}
\]
is a reduced expression for $\epsilon_i\in\Lambda
\subset \mathcal{W}$.
The following result follows from 
Theorem \ref{basic} and from the algebraic 
structure of affine Hecke algebras (see Lusztig \cite{L} and Noumi \cite{N}).
\begin{cor}
The operators $Y_i\in \mathcal{D}_q$
\textup{(}$i=1,\ldots,n$\textup{)} pair-wise commute.
\end{cor}
The analogue of Cherednik's double affine Hecke algebra can be
defined explicitly in terms of generators and relations, see Sahi
\cite{Sa} (see also \cite{St1} and \cite{St2} for alternative 
presentations). We take a short-cut in the present paper by defining the
double affine Hecke algebra $\mathcal{H}$
directly in its image under Cherednik's 
faithful realization of $\mathcal{H}$ 
as $q$-difference reflection operators and
multiplication operators (the analogue of the Schr{\"o}dinger
realization of the Heisenberg algebra). For this we 
observe that $\mathcal{M}$, 
and hence $\mathcal{A}\subset \mathcal{M}$, is naturally
embedded in $\mathcal{D}_q$ via
\[ g\mapsto g\otimes e,\qquad g\in\mathcal{M},
\]
where $e\in\mathcal{W}$ is the identity element.   
\begin{Def}
The double affine Hecke algebra $\mathcal{H}=\mathcal{H}_\alpha\subset
\mathcal{D}_q$ is the unital
sub-algebra of $\mathcal{D}_q$  generated
by $T_j$ \textup{(}$j=0,\ldots,n$\textup{)} and $\mathcal{A}$.
\end{Def}
It is clear that $\mathcal{H}$ is also generated as algebra by
$Y_i^{\pm 1},x_i^{\pm 1}$ and $T_i$ for $i=1,\ldots,n$.
Observe that $\mathcal{H}$ acts on 
$\mathcal{A}\subset \mathcal{M}$ and on $\mathcal{O}\subset\mathcal{M}$
by restriction of the natural action of $\mathcal{D}_q$ on $\mathcal{M}$.


\subsection{The duality isomorphism}

The Heisenberg algebra has a copy of $\hbox{SL}_2(\mathbb{Z})$
in its automorphism group, generated within the metaplectic
representation by the Fourier transform and by multiplication
by the Gaussian. 
These two isomorphisms, as well as their
realization within the metaplectic representation, generalize
to the set-up of double affine Hecke algebras. 

In this subsection we introduce the isomorphism of the double affine
Hecke algebra associated with (the analogue of) the Fourier transform. 
In the following subsection, we consider the
isomorphism associated with multiplication by the generalized Gaussian.

We first recall the so-called
$\epsilon$-transform, see Sahi \cite{Sa}.
\begin{thm}\label{epsilon}
There exists a unique, unital algebra isomorphism
\[
\epsilon=\epsilon_{\alpha}: \mathcal{H}=\mathcal{H}_\alpha\rightarrow
\mathcal{H}_{\ddagger\sigma}
\]
satisfying
\[
\epsilon(Y_i)=x_i,\qquad
\epsilon(T_i)=T_i^{\ddagger\sigma}{}^{-1},
\qquad
\epsilon(x_i)=Y_i^{\ddagger\sigma}
\]
for $i=1,\ldots,n$. 
Furthermore, $\epsilon^{-1}= \epsilon_{\ddagger\sigma}$.
\end{thm}

By \cite[Lem. 7.3]{St1} the assignment
\[ x_i\mapsto x_i^{-1},\qquad T_j\mapsto T_j^{\ddagger}{}^{-1}
\]
for $i=1,\ldots,n$ and $j=0,\ldots,n$ uniquely extends to a unital
algebra isomorphism $\dagger: \mathcal{H}\rightarrow
\mathcal{H}_{\ddagger}$. Furthermore, if $I: \mathcal{M}\rightarrow
\mathcal{M}$ denotes the involution
\[ \bigl(Ig\bigr)(x)=g(x^{-1}),\qquad g\in\mathcal{M}
\]
where $x^{-1}=(x_1^{-1},x_2^{-1},\ldots,x_n^{-1})$,
then it is easy to check that
\begin{equation}\label{inversionrelation}
I\circ X=\dagger(X)\circ I,\qquad X\in\mathcal{H}.
\end{equation}
Composing now $\epsilon$ with $\dagger$, we
obtain the ``duality isomorphism'' (see \cite[Def. 7.5]{St1}) 
we are looking for.
\begin{cor}\label{sigma}
The map $\sigma:=\dagger_{\ddagger\sigma}\circ\epsilon:
\mathcal{H}\rightarrow \mathcal{H}_{\sigma}$ is an algebra
isomorphism, satisfying
\[ \sigma(x_1^{-1}T_0Y_1^{-1})=T_0^\sigma,\qquad
\sigma(T_i)=T_i^\sigma,\qquad \sigma(Y_i)=x_i^{-1}
\]
for $i=1,\ldots,n$.
\end{cor}
\begin{rem}
Note here that we use the notation $\sigma$ for the isomorphism
$\sigma$ of the double affine Hecke algebra, as well as for the
involution on the (difference) multiplicity function
$\alpha=(\mathbf{t},q^{\frac{1}{2}})$. Which one is used should
always be clear from the context (as a generic rule, $\sigma$ as involution
on difference multiplicity functions will always be written as a
sub-(or super-)index, in contrast with $\sigma$ regarded as
isomorphism of the double affine Hecke algebra). Such conventions
will also hold for other (anti-)isomorphisms 
of the double affine Hecke algebra defined
at later stages.
\end{rem}

In \cite[Prop. 8.8]{St1} it is proven that the nonsymmetric
(polynomial) Macdonald-Koornwinder 
transform induces the automorphism $\sigma$ on
$\mathcal{H}$, see also Section \ref{Pol}.

\subsection{The Gausssian}

We assume in this subsection that the parameter $q^{\frac{1}{2}}$
in the generic difference multiplicity
function $\alpha=(t_0,u_0,t_n,u_n,t,q^{\frac{1}{2}})$ has
modulus unequal to one.

Cherednik introduced an analogue of the Gaussian for double affine
Hecke algebras for reduced root systems, see
\cite{C1} and references therein.
In this subsection we generalize Cherednik's construction to the
nonreduced set-up. 

Recall the standard notations for $q$-shifted factorials (see \cite{GR}),
\begin{equation}\label{qshift}
\bigl(z_1,\ldots,z_m;q\bigr)_{k}=\prod_{i=1}^m\bigl(z_i;q\bigr)_k,
\qquad \bigl(z;q\bigr)_{k}=\prod_{i=0}^{k-1}(1-zq^i),
\end{equation}
where $z_1,\ldots,z_m,z\in\mathbb{C}$ and $k\in \mathbb{N}\cup\{\infty\}$ 
with $\mathbb{N}=\{0,1,2,\ldots\}$. Here the modulus of $q$ should
be taken $<1$ when $k=\infty$. 
\begin{Def}\label{DefGaussian}
{\bf i)} For $|q|<1$ the Gaussian 
$G(\cdot)=G_{\alpha}(\cdot)\in \mathcal{M}$ is defined by the infinite
product
\begin{equation*}
\begin{split}
G(x)&=\prod_{f\in \mathcal{W}a_0\cap
 R_+}(1+t_ft_{f/2}^{-1}x^{f/2})^{-1}\\
&=\prod_{i=1}^n\bigl(-q^{\frac{1}{2}}t_0u_0^{-1}x_i,
-q^{\frac{1}{2}}t_0u_0^{-1}/x_i;q\bigr)_{\infty}^{-1}.
\end{split}
\end{equation*}
{\bf ii)} For $|q|>1$ the Gaussian 
$G(\cdot)=G_{\alpha}(\cdot)\in\mathcal{M}$ is defined by the infinite
product
\begin{equation*}
\begin{split}
G(x)&=\prod_{f\in \mathcal{W}a_0\cap R_-}
(1+t_ft_{f/2}^{-1}x^{f/2})\\
&=\prod_{i=1}^n\bigl(-q^{-\frac{1}{2}}t_0u_0^{-1}x_i,
-q^{-\frac{1}{2}}t_0u_0^{-1}/x_i;q^{-1}\bigr)_{\infty}.
\end{split}
\end{equation*}
\end{Def}

\begin{rem}\label{remabc}
{\bf a)} We used two different ways of representing the Gaussian
in Definition \ref{DefGaussian},
the first one in terms of the root data, the second one in terms
of $q$-shifted factorials. 
The equivalence of the two expressions 
can be easily verified using that
\begin{equation}\label{orbitroots}
\mathcal{W}a_0\cap R_{\pm}= \Sigma^l+(1\pm 2\mathbb{N})\delta
\end{equation}
where $\Sigma^l\subset \Sigma$ is the set of roots of length two:
\[ \Sigma^l=\{ \pm 2\epsilon_i \, | \, i=1,\ldots,n\}.
\]
{\bf b)} The Gaussians for $|q|<1$ and $|q|>1$ are related by
the formula
\begin{equation}\label{Gaussianinversion}
G(x)=G_{\ddagger\tau}(x)^{-1}.
\end{equation}
{\bf c)}
If $t_0=u_0$ and $|q|<1$, then 
\begin{equation}\label{Gless1}
G(x)^{-1}=\prod_{i=1}^n\bigl(-q^{1/2}x_i,-q^{1/2}x_i^{-1};q\bigr)_{\infty}
=\bigl(q;q\bigr)_{\infty}^{-n}
\sum_{\lambda\in\Lambda}q^{\langle\lambda,\lambda\rangle/2}x^{\lambda}
\end{equation}
in view of the Jacobi triple product identity \cite[(1.6.1)]{GR}
for theta-functions.
Up to the (irrelevant) constant $\bigl(q;q\bigr)_{\infty}^{-n}$,
this coincides with the definition of the Gaussian of type
$\widetilde{C}_n$ as given by Cherednik, see \cite{C1} and references
therein.
\end{rem}

\begin{prop}\label{Gaussianconj}
The inner automorphism
\[X\mapsto G\,X\,G^{-1}
\]
of $\mathcal{D}_q$ restricts to an algebra isomorphism
$\tau=\tau_{\alpha}: \mathcal{H}\rightarrow \mathcal{H}_{\tau}$.
Its action on a set of algebraic generators of $\mathcal{H}$ is given
by
\[
\tau(x_i)=x_i, \qquad \tau(T_i)=
T_i^{\tau},\qquad
\tau(T_0)=q^{-\frac{1}{2}}x_1T_0^{\tau}{}^{-1}
\]
for $i=1,\ldots,n$.
\end{prop}
\begin{proof}
It suffices to  show that conjugation by $G$ on the algebraic generators of 
$\mathcal{H}\subset \mathcal{D}_q$
is as stated in the proposition.
Clearly, conjugation by $G$ maps the coordinate
functions $x_i$ to themselves.
Conjugation with $G$ maps $T_i$ to
$T_i^{\tau}$ for $i=1,\ldots,n$ since
$G$ is $W_0$-invariant and since
$T_i=T_i^{\tau}$ for $i=1,\ldots,n$ (indeed,
$\tau$ acts on $\alpha=(\mathbf{t},q^{\frac{1}{2}})$
by interchanging the parameters $t_0=t_{a_0}$
and $u_0=t_{a_0/2}$, which only occur in the generator $T_0$).

It remains to prove that $\tau(T_0)=q^{-\frac{1}{2}}x_1T_0^{\tau}{}^{-1}$.
Observe that the action of the simple reflection
$r_0\in \mathcal{W}\subset \mathcal \mathcal{D}_q$ on the
Gaussian $G\in \mathcal{M}$ is given by
\[
(r_0G)(x)=G(qx_1^{-1},x_2,\ldots,x_n)
=\frac{\bigl(1+t_0u_0^{-1}q^{\frac{1}{2}}x_1^{-1}\bigr)}
{\bigl(1+t_0u_0^{-1}q^{-\frac{1}{2}}x_1\bigr)}G(x),
\]
hence
\begin{equation}\label{eq1}
\begin{split}
t_0^{-1}c_{a_0}(x)\frac{G(x)}{(r_0G)(x)}&=
t_0^{-1}\frac{(1-t_0u_0q^{\frac{1}{2}}x_1^{-1})
(1+t_0u_0^{-1}q^{-\frac{1}{2}}x_1)}
{(1-qx_1^{-2})}\\
&=q^{-\frac{1}{2}}x_1u_0^{-1}c_{a_0}^\tau(x).
\end{split}
\end{equation}
Furthermore, observe that
\begin{equation}\label{eq2}
\begin{split}
t_0-t_0^{-1}c_{a_0}(x)&=\frac{(t_0-t_0^{-1})+
(u_0-u_0^{-1})q^{\frac{1}{2}}x_1^{-1}}
{1-qx_1^{-2}}\\
&=q^{-\frac{1}{2}}x_1\frac{(u_0-u_0^{-1})qx_1^{-2}+
(t_0-t_0^{-1})q^{\frac{1}{2}}x_1^{-1}}
{1-qx_1^{-2}}\\
&=q^{-\frac{1}{2}}x_1\bigl(u_0^{-1}-u_0^{-1}c_{a_0}^\tau(x)\bigr).
\end{split}
\end{equation}
We then compute in $\mathcal{D}_q$, using
\eqref{eq1} and \eqref{eq2},
\begin{equation*}
\begin{split}
G\,T_0\,G^{-1}&=t_0-t_0^{-1}c_{a_0}+t_0^{-1}c_{a_0}\frac{G}{r_0G}r_0\\
&=
q^{-\frac{1}{2}}x_1\bigl(u_0^{-1}-u_0^{-1}c_{a_0}^\tau+u_0^{-1}c_{a_0}^\tau 
r_0\bigr)\\
&=q^{-\frac{1}{2}}x_1T_0^{\tau}{}^{-1},
\end{split}
\end{equation*}
which completes the proof of the proposition.
\end{proof}

\subsection{$\hbox{SL}_2$-type commutation 
relations}
Let $w\in\mathcal{W}$ and choose a reduced expression 
$w=r_{i_1}r_{i_2}\cdots r_{i_l}$. We denote
\[ T_w=T_{i_1}T_{i_2}\cdots T_{i_l}\in \mathcal{H},
\]
which is independent of the choice of reduced expression since
the $T_j$'s satisfy the $\widetilde{C}_n$-braid relations. Similarly,
$t_w=t_{i_1}t_{i_2}\cdots t_{i_l}\in\mathbb{C}$ is independent of
the reduced expression of $w$, where $t_j=t_{a_j}$ is the value of the
multiplicity function $\mathbf{t}$ at the simple root $a_j\in R$.

Let $w_0\in W_0$ be the longest Weyl group element in $W_0$.
Then for all $X\in \mathcal{H}$, we have the relations
\begin{equation}\label{SLrelations}
\begin{split}
(\sigma_\sigma\circ\sigma\circ\sigma_\sigma\circ\sigma)(X)&=
T_{w_0}^{-2}XT_{w_0}^2=
(\sigma_\sigma\circ
\tau_{\sigma\tau}\circ\sigma_{\tau\sigma\tau}\circ
\tau_{\tau\sigma}\circ\sigma_\tau\circ\tau)(X),\\
(\sigma_{\tau\sigma}\circ\sigma_\tau\circ\tau)(X)&=
(\tau\circ\sigma_\sigma\circ\sigma)(X).
\end{split}
\end{equation}
Here the composition of isomorphisms on the
right hand side of the first equality in \eqref{SLrelations}
is well defined since $\sigma\tau\sigma=\tau\sigma\tau$ when acting on 
difference multiplicity functions. 

The relations \eqref{SLrelations} between the isomorphisms $\sigma$ and
$\tau$ of $\mathcal{H}$ should be viewed as the analogues
of the characterizing relations
\begin{equation*}
\widetilde{\sigma}^4=\left(\begin{matrix} 1 & 0\\ 0 &
    1\end{matrix}\right)=(\widetilde{\sigma}\widetilde{\tau})^3,\qquad
\widetilde{\sigma}^2\widetilde{\tau}=\widetilde{\tau}\widetilde{\sigma}^2
\end{equation*}
for the generators
\begin{equation*}
\widetilde{\sigma}=\left(\begin{matrix} 0 & 1\\ -1 &
    0\end{matrix}\right),
\qquad \widetilde{\tau}=\left(\begin{matrix} 1 & 1\\ 0 & 1
\end{matrix}\right)
\end{equation*}
of the modular group $\hbox{SL}_2(\mathbb{Z})$. 

We do not give a detailed proof 
of the relations \eqref{SLrelations} since they are not needed
in the remainder of the paper. We refer to \cite[Thm. 2.3]{C3n} for
the analogous result in case of reduced root systems.
We only note here that the first equality in \eqref{SLrelations}
follows from the simpler formula 
$(\sigma_\sigma\circ\sigma)(X)=T_{w_0}^{-1}XT_{w_0}$ for $X\in\mathcal{H}$.

\section{General remarks on difference Fourier transforms}
\label{General}
In this section the concept of Fourier transforms 
associated with $\sigma$ is 
introduced in an informal manner. Rigorous statements
will follow in subsequent sections. 

Let $V=V_\alpha$ and $W=W_\alpha$ be some function spaces on
the complex torus $(\mathbb{C}^\times)^n$ on which
$\mathcal{H}$ acts as $q$-difference reflection operators.
We assume that $V$ is stable under the inversion operator
$(Ig)(x):=g(x^{-1})$. This in particular implies that $\mathcal{H}_{\ddagger}$
acts as $q^{-1}$-difference reflection operators on $V$. The relation
between the two actions on $V$ is given by the formula 
\eqref{inversionrelation}.

The starting point of our considerations
is the search of explicit linear transformations
\[ F=F_\alpha: V\rightarrow W_\sigma,
\]
which satisfy
\begin{equation}\label{transFF}
 F\circ X=\sigma(X)\circ F,\qquad \forall\, X\in\mathcal{H},
\end{equation}
where $\sigma$ is the duality isomorphism of $\mathcal{H}$, see
Corollary \ref{sigma}. We call such a linear endomorphism 
$F$ a {\it Fourier transform associated with} $\sigma$, or sometimes
just a {\it difference Fourier transform}. 

We are mainly interested in difference Fourier transforms $F$ 
which can be realized as integral transforms. 
The property \eqref{transFF} then formally translates to
explicit transformation properties of the associated 
kernels under the action of the double affine Hecke
algebra. To be more precise, we first need to 
recall certain anti-isomorphisms of the double affine Hecke
algebra $\mathcal{H}$ which play the role of adjoint
maps.

\begin{lem}
{\bf i)}
There exists a unique unital anti-algebra isomorphism $\ddagger=
\ddagger_\alpha:
\mathcal{H}\rightarrow \mathcal{H}_\ddagger$ satisfying
\[\ddagger(T_j)=T_j^{\ddagger}{}^{-1},\qquad
\ddagger(x_i)=x_i^{-1}
\]
for $j=0\ldots,n$ and $i=1,\ldots,n$. Furthermore, 
$\ddagger^{-1}=\ddagger_\ddagger$.\\
{\bf ii)} There exists a unique unital anti-algebra isomorphism 
$\iota=\iota_\alpha:
\mathcal{H}\rightarrow \mathcal{H}$ satisfying
\[\iota(T_j)=T_j,\qquad
\iota(x_i)=x_i
\]
for $j=0\ldots,n$ and $i=1,\ldots,n$. Furthermore, $\iota^{-1}=\iota$.
\end{lem}
\begin{proof}
{\bf i)} See e.g. \cite[Prop. 7.1]{Sa}.\\
{\bf ii)} This follows e.g. from the fact that
$\iota:=\ddagger_\ddagger\circ\dagger=
\dagger_\ddagger\circ\ddagger: \mathcal{H}\rightarrow \mathcal{H}$
is a unital anti-algebra isomorphism which fixes the generators $T_j$
and $x_i$ for $j=0,\ldots,n$ and $i=1,\ldots,n$.
\end{proof}

Suppose now that we are given a linear map $F:V\rightarrow W_\sigma$ 
of the form
\begin{equation}\label{Ftest} 
\bigl(Fg\bigr)(\gamma)=
\bigl(g,\mathfrak{E}_\ddagger(\gamma^{-1},\cdot)\bigr),
\qquad g\in V
\end{equation}
with $\bigl(\cdot,\cdot\bigr)=\bigl(\cdot,\cdot\bigr)_\alpha$ 
some bilinear form satisfying
\[\bigl(Xg,h\bigr)=\bigl(g,\ddagger(X)h\bigr)
\]
for $X\in \mathcal{H}$,
$g\in V$ and $h$ in some completion of $V$.
At a formal level, the condition that the map
$F: V\rightarrow W_\sigma$ defines
a Fourier transform associated with $\sigma$ corresponds to
the transformation behaviour
\begin{equation}\label{Trans}
\bigl(X\mathfrak{E}_\ddagger(\gamma,\cdot)\bigr)(x)=
\bigl(\psi_\ddagger(X)\mathfrak{E}_\ddagger(\cdot,x)\bigr)(\gamma),
\qquad X\in \mathcal{H}_\ddagger
\end{equation}
of the kernel $\mathfrak{E}_\ddagger$,
where $\psi=\psi_\alpha: \mathcal{H}\rightarrow \mathcal{H}_\sigma$
is the anti-isomorphism
\[ \psi=\dagger_{\ddagger\sigma}\circ\sigma_\ddagger\circ\ddagger.
\]
Here the double affine Hecke algebra acts by $q^{-1}$-difference
reflection operators on both sides of \eqref{Trans}. 
The anti-isomorphism $\psi$ is the 
so-called {\it duality anti-isomorphism}, see \cite{Sa} and \cite{St1}. 
In particular, it has the 
following special properties.
\begin{prop}
The map
$\psi=\dagger_{\ddagger\sigma}\circ\sigma_\ddagger\circ\ddagger$
is the unique unital anti-isomorphism $\psi: \mathcal{H}\rightarrow
\mathcal{H}_\sigma$ satisfying
\[ \psi(x_i)=Y_i^{\sigma}{}^{-1},\qquad \psi(T_i)=T_i^{\sigma},
\qquad \psi(Y_i)=x_i^{-1}
\]
for $i=1,\ldots,n$. In particular, $\psi^{-1}=\psi_\sigma$.
\end{prop}
\begin{proof}
This is an easy verification, see e.g. \cite[Sect. 7]{Sa} for details. 
\end{proof}

The transformation behaviour \eqref{Trans} for a kernel
$\mathfrak{E}_\ddagger$ hints at several
important (and desirable) properties for $\mathfrak{E}_\ddagger$. 
For instance the fact that $\psi$ maps the commuting
$Y$-operators to multiplication operators implies that 
$\mathfrak{E}_\ddagger(\gamma,\cdot)$ is a common eigenfunction
of $Y_i^\ddagger\in\mathcal{H}_\ddagger$ with eigenvalue
$\gamma_i^{-1}$ for $i=1,\ldots,n$. Furthermore, the ``involutivity''
$\psi_\ddagger^{-1}=\psi_{\ddagger\sigma}$ of the anti-isomorphism
$\psi_\ddagger$ hints at the symmetric role of the geometric parameter
$x$ and the spectral parameter $\gamma$ in
$\mathfrak{E}_\ddagger(\gamma,x)$.
In fact, the kernels we encounter indeed turn out to 
satisfy the {\it duality property}
\begin{equation}\label{dualitypropp}
\mathfrak{E}_\ddagger(\gamma,x)=\mathfrak{E}_{\ddagger\sigma}(x,\gamma)
\end{equation}
after a proper choice of normalization.

The property that a transform $F$  
satisfies the transformation behaviour \eqref{transFF} 
turns out to be a very strong condition.
In fact, for the transforms we encounter in this paper,  
the corresponding modules $V$ and $W$ are cyclic
$\mathcal{H}$-modules, or they consist of a direct sum of two cyclic
$\mathcal{H}$-modules. In each case, the cyclic vectors are given
explicitly. The transformation behaviour \eqref{transFF} of $F$ 
reduces the study of the transform 
to the explicit evaluation of the image of the cyclic vectors 
under $F$. 
Again by the transformation behaviour \eqref{transFF}
of $F$, the computation
of the image of the cyclic vectors reduce to explicit constant
term evaluations (e.g. Macdonald type constant term identities and 
Macdonald-Mehta type identities).

To analyze inversion formulas for difference Fourier transforms,
we make a similar, formal analysis for Fourier
transforms $J_\sigma=J_{\alpha_\sigma}: W_\sigma\rightarrow V$
associated with the isomorphism $\sigma^{-1}$, i.e. linear maps satisfying
the opposite transformation behaviour
\begin{equation}\label{oppositetransform}
 J_\sigma\circ X=\sigma^{-1}(X)\circ J_\sigma,
\qquad \forall\,X\in\mathcal{H}_\sigma,
\end{equation}
where $V$ and $W$ are as before. When no confusion can arise on 
the underlying isomorphism, we also use the terminology
{\it difference Fourier transforms} for such transforms $J_\sigma$.

We now assume that $J_\sigma$ is of the form
\begin{equation}\label{Jtest}
(J_{\sigma}g)(x)=\lbrack g,\mathfrak{E}(\cdot,x)\rbrack_\sigma,
\end{equation}
with $\lbrack \cdot,\cdot\rbrack=\lbrack \cdot,\cdot\rbrack_\alpha$
a bilinear form satisfying
\[ \lbrack Xg,h\rbrack=\lbrack g,\iota(X)h\rbrack
\]
for $X\in\mathcal{H}$, $g\in W$ and for $h$ in some
completion of $W$. 
The fact that $J_\sigma$ is a Fourier transform associated
with $\sigma^{-1}$ then formally relates to the 
transformation behaviour
\begin{equation}\label{Trans2}
\bigl(X\mathfrak{E}(\gamma,\cdot)\bigr)(x)=\bigl(\psi(X)
\mathfrak{E}(\cdot,x)\bigr)(\gamma),\qquad X\in\mathcal{H}
\end{equation}
of the kernel $\mathfrak{E}$, since $\psi=\iota_\sigma\circ\sigma$.
Here the double affine Hecke algebra acts on both sides of 
\eqref{Trans2} by
$q$-difference reflection operators. In view of the (expected) duality
\eqref{dualitypropp} of the kernels, it is therefore plausible that
a given difference Fourier transform $F: V\rightarrow W_\sigma$ of the
form \eqref{Ftest} can be inverted by an explicit transform $J_\sigma$ 
which has the same kernel as $F$, but depending now on the difference
multiplicity function $\alpha_{\ddagger\sigma}$ instead of $\alpha$.
All these features are shown to be true for the difference Fourier
transforms considered in this paper.

It is clear from the above considerations that the first priority
should be to construct kernels $\mathfrak{E}$ satisfying the
transformation behaviour \eqref{Trans2} under the action of the double
affine Hecke algebra. The first example of 
such a kernel $\mathfrak{E}(\gamma,x)$ can be defined in terms of  
Mac\-do\-nald-Koorn\-win\-der polynomials, but the spectral
parameter $\gamma$ then runs through the discrete, polynomial
spectrum of the operators $Y_i\in\mathcal{H}$. This restrictive
kernel can be used to define the so-called polynomial
Macdonald-Koornwinder transform (see \cite{St1}), which gives a first
example of a difference Fourier transform. In Section 4 we 
explain the concepts introduced in this section for the polynomial 
Macdonald-Koornwinder transform.  

In Section 5 we meromorphically extend this polynomial kernel 
$\mathfrak{E}$ while preserving the desired transformation behaviour
\eqref{Trans2} under the action of the double affine Hecke algebra.   
This kernel was written down explicitly by Cherednik
\cite{C1} for reduced root systems. In Section 7 and Section 8 we
construct and study related difference Fourier transforms
in full detail, following closely the general philosophy as explained
in this section.


\section{The Macdonald-Koornwinder transform}\label{Pol}
In this section we recall a known example of a Fourier transform
$F:V\rightarrow W_\sigma$ 
associated with $\sigma$, the so-called 
Macdonald-Koornwinder transform (see \cite{St1}).
In this case $F$ is defined as an 
integral transform with kernel expressed in terms of 
nonsymmetric Macdonald-Koornwinder polynomials. In subsection 4.1
we introduce the modules $V$ and $W$,
in subsection 4.2 we introduce the kernel 
$\mathfrak{E}$ and in subsection 4.3 we define
the associated bilinear forms $\bigl(\cdot,\cdot\bigr)$ and
$\lbrack \cdot,\cdot\rbrack$. In subsection 4.4
we construct the Macdonald-Koornwinder transform 
and its inverse. Most results
can be found directly or indirectly in \cite{N}, \cite{Sa} or in
\cite{St1} (see also the lecture notes \cite{St2}).
I have decided to be quite detailed in this section, since
the results play a key role throughout this paper. Furthermore, 
the theory is presented in such a way that it directly fits into
the general scheme of difference Fourier transforms as
discussed in Section 3.

Throughout this section we fix a generic multiplicity function
$\alpha=(\mathbf{t},q^{\frac{1}{2}})$. After subsection 4.2 we impose
extra conditions on $\alpha$, see the beginning of subsection 4.3.


\subsection{The modules}
For the Macdonald-Koornwinder transform $F: V\rightarrow W_\sigma$,
we take $V$ to be 
the cyclic $\mathcal{H}$-module $\mathcal{A}$,
with cyclic vector $1\in \mathcal{A}$ the Laurent polynomial
identically equal to one. The target space 
$W=W_\alpha=\mathcal{F}_0(\mathcal{S}_{\ddagger\sigma})$ is
the linear space of functions $g: \mathcal{S}_{\ddagger\sigma}\rightarrow
\mathbb{C}$ with finite support, 
where  $\mathcal{S}=\mathcal{S}_\alpha\subset (\mathbb{C}^\times)^n$
is the spectrum of the commuting elements
$Y_1,\ldots,Y_n\in\mathcal{H}_\alpha$ considered as endomorphism
of $\mathcal{A}$ via their action as $q$-difference reflection
operators.

We now introduce the $\mathcal{H}$-module structure on 
$\mathcal{F}_0(\mathcal{S}_{\ddagger\sigma})$. In \cite[Prop. 8.8]{St1}
the $\mathcal{H}$-module structure on
$\mathcal{F}_0(\mathcal{S}_{\ddagger\sigma})$ was 
constructed in an indirect manner, using
the Macdonald-Koornwinder transform in an essential way. 
In order to emphasize the natural order of definitions and results
in the study of Fourier transforms associated with $\sigma$ (see
Section 3),  we give here a detailed account on a direct
construction of the $\mathcal{H}$-module structure on
$\mathcal{F}_0(\mathcal{S}_{\ddagger\sigma})$.

We first need to
recall the explicit form of the polynomial spectrum 
$\mathcal{S}=\mathcal{S}_\alpha$, see e.g. \cite{Sa} and \cite{St1}.
It is naturally parametrized by the lattice $\Lambda$, 
\[\mathcal{S}=\{s_\lambda\,|\,
\lambda\in\Lambda\},
\]
with $s_\lambda=s_\lambda^\alpha$ given by
\[s_\lambda=(s_{\lambda,1},s_{\lambda,2},\ldots,s_{\lambda,n}),\qquad\quad
s_{\lambda,i}=(t_nt_0)^{(\rho_l(\lambda),\epsilon_i)}
t^{(\rho_m(\lambda),\epsilon_i)}q^{(\lambda,\epsilon_i)},
\]
where $\rho_m(\lambda), \rho_l(\lambda)\in\Lambda$
are given by
\[
\rho_m(\lambda)=\sum_{\alpha\in\Sigma_m^+}
\hbox{sgn}(\langle \lambda,\alpha\rangle)\alpha^\vee,\qquad
\rho_l(\lambda)=\sum_{\alpha\in\Sigma_l^+}
\hbox{sgn}(\langle \lambda,\alpha\rangle)\alpha^\vee.
\]
Here $\hbox{sgn}(m)$ is equal to $1$ if $m\in \mathbb{N}$ and
equal to $-1$ if $m\in \mathbb{Z}_{<0}$, $\Sigma_m^+$ and $\Sigma_l^+$ are
the positive roots of $\Sigma$ of squared length two and four
respectively, and  
$\alpha^\vee=2\alpha/\langle \alpha,\alpha\rangle$
is the co-root of $\alpha$.

Observe that the spectrum $\mathcal{S}$
does not depend on the parameters $u_0$ and $u_n$ of the
difference multiplicity function $\alpha$, and that the spectral
points $s_\lambda=s_\lambda^\alpha$ satisfy
\[ s_\lambda^\ddagger=s_\lambda^{-1},\qquad \forall\,\lambda\in\Lambda.
\]
Note also that $s_0=s_0^\alpha\in\mathcal{S}$
is the spectral point corresponding to the common eigenfunction $1\in
\mathcal{A}$ of the $Y$-operators $Y_1,\ldots,Y_n$.

Let $\mathcal{F}(\mathcal{S}_{\ddagger\sigma})$ be the space
of functions $g: \mathcal{S}_{\ddagger\sigma}\rightarrow \mathbb{C}$
(without finiteness conditions). 
We first introduce an action of $\mathcal{W}$ on 
$\mathcal{F}(\mathcal{S}_{\ddagger\sigma})$, which we call
the dot-action. It is defined by
\[ (w\cdot g)(s_{\lambda}^{\ddagger\sigma})=
g(s_{w^{-1}\cdot\lambda}^{\ddagger\sigma}),\qquad w\in \mathcal{W},\,\,
\lambda\in\Lambda
\]
for $g\in \mathcal{F}(\mathcal{S}_{\ddagger\sigma})$
with the action of $\mathcal{W}$ on $\Lambda$ 
defined by
\begin{equation*}
r_j\cdot\lambda=
\begin{cases}
(-1-\lambda_1,\lambda_2,\lambda_3,\ldots,\lambda_n),\qquad &j=0,\\
(\lambda_1,\ldots,\lambda_{j-1},\lambda_{j+1},\lambda_j,\lambda_{j+2},
\ldots,\lambda_n),\qquad &1\leq j\leq n-1,\\
(\lambda_1,\ldots,\lambda_{n-1},-\lambda_n),\qquad &j=n.
\end{cases}
\end{equation*}
Note that $\lambda\cdot\mu=\lambda+\mu$ for $\lambda,\mu\in\Lambda$,
where $\lambda\in\Lambda$ is viewed as an affine Weyl group element
in $\mathcal{W}\simeq W_0\ltimes \Lambda$. 

\begin{lem}\label{actioncompatibleW}
Fix $j\in\{0,\ldots,n\}$ arbitrary.\\
{\bf a)} For any function $g: (\mathbb{C}^\times)^n\rightarrow \mathbb{C}$
we have
\[\bigl(r_jg\bigr)(s_\lambda^{\ddagger\sigma})=
\bigl(r_j\cdot g|_{\mathcal{S}_{\ddagger\sigma}}\bigr)
(s_{\lambda}^{\ddagger\sigma})
\]
when $\lambda\in\Lambda$ satisfies $r_j\cdot\lambda\not=\lambda$, 
where $r_jg$ is the action of $r_j\in\mathcal{W}$ on $g$ 
as constant coefficient $q$-difference reflection operator
\textup{(}see \eqref{Wfund}\textup{)}.\\
{\bf b)} The rational function $c_{a_j}\in\mathbb{C}(x)$ is
regular at the spectral points $s\in \mathcal{S}_{\ddagger\sigma}$.
Furthermore,
\begin{equation}\label{zerocorr}
c_{a_j}(s_\lambda^{\ddagger\sigma})=0\quad
\Leftrightarrow\quad  r_j\cdot\lambda=\lambda
\end{equation}
for all $\lambda\in\Lambda$.
\end{lem}
\begin{proof}
{\bf a)} See the proof of \cite[Thm 5.3]{Sa}.\\
{\bf b)} It is important to recall here that the difference
multiplicity function $\alpha$ is assumed to be generic.
The regularity follows then from the explicit
expressions for $c_{a_j}$ and 
$s\in\mathcal{S}_{\ddagger\sigma}$. In a similar manner one checks
that $c_{a_0}(s_\lambda^{\ddagger\sigma})\not=0$ for all
$\lambda\in\Lambda$. Since 
$r_0\cdot\lambda\not=\lambda$ for all $\lambda\in\Lambda$,
this proves \eqref{zerocorr} for $j=0$. 

For $j\in\{1,\ldots,n-1\}$ we observe that
\[ c_{a_j}(s_\lambda^{\ddagger\sigma})=0\quad \Leftrightarrow\quad
\bigl(s_\lambda^{\ddagger\sigma})^{a_j}=
(t_nu_n)^{-(\rho_m(\lambda),a_j)}t^{-(\rho_m(\lambda),a_j)}
q^{-(\lambda,a_j)}=t^{-2}
\]
by the explicit expression for $c_{a_j}$ (see \eqref{cf2}).
On the other hand, it is easy to check that
$\bigl(s_\lambda^{\ddagger\sigma}){}^{a_j}=t^{-2}$ if and only if
$r_j\cdot\lambda=\lambda$.

For $j=n$ the condition
$c_{a_j}(s_\lambda^{\ddagger\sigma})=0$ is equivalent to the condition that
$\bigl(s_\lambda^{\ddagger\sigma})^{\epsilon_n}$ is equal to  
$(t_nu_n)^{-1}$ or to $-t_n^{-1}u_n$, in view of \eqref{cf1}. 
By the explicit expression
for $s_\lambda^{\ddagger\sigma}$, only the equality
$\bigl(s_\lambda^{\ddagger\sigma})^{\epsilon_n}=(t_nu_n)^{-1}$ can happen
for some $\lambda\in\Lambda$. Furthermore, one easily checks that
$\bigl(s_\lambda^{\ddagger\sigma})^{\epsilon_n}=(t_nu_n)^{-1}$ if and only
if $r_n\cdot\lambda=\lambda$, which proves \eqref{zerocorr} for $j=n$.
\end{proof}
Using the above lemma we can define an action of $\mathcal{H}$
on $\mathcal{F}(\mathcal{S}_{\ddagger\sigma})$ as follows.

\begin{lem}\label{actioncompatibleH}
For $g\in\mathcal{F}(\mathcal{S}_{\ddagger\sigma})$ and
$X\in \mathcal{H}$, set
\begin{equation}\label{defaction}
\bigl(X\cdot g\bigr)(s):=\bigl(X\widetilde{g}\bigr)(s),\qquad 
s\in\mathcal{S}_{\ddagger\sigma},
\end{equation}
where $\widetilde{g}: (\mathbb{C}^\times)^n\rightarrow \mathbb{C}$ 
is an arbitrary function satisfying 
$\widetilde{g}|_{\mathcal{S}_{\ddagger\sigma}}=g$,
and with the action of $X\in\mathcal{H}$ on $\widetilde{g}$
in the right hand side of \eqref{defaction} by $q$-difference
reflection operators. Then \eqref{defaction} is well defined
\textup{(}i.e. independent of the choice of extension $\widetilde{g}$
of $g$\textup{)}, and it defines a left action of 
$\mathcal{H}$ on $\mathcal{F}(\mathcal{S}_{\ddagger\sigma})$.
\end{lem}
\begin{proof}
It suffices to prove that formula \eqref{defaction} is independent
of the choice of extension $\widetilde{g}$ of 
$g\in\mathcal{F}(\mathcal{S}_{\ddagger\sigma})$.
This is clear when $X\in\mathcal{H}$ is multiplication by a Laurent
polynomial $p\in\mathcal{A}$. 
Lemma \ref{actioncompatibleW} shows that
\begin{equation}\label{Tdot}
\bigl(T_j\widetilde{g}\bigr)(s)=t_jg(s)
+t_j^{-1}c_{a_j}(s)\bigl((r_j\cdot g)(s)-g(s)\bigr),\qquad \forall\,
s\in \mathcal{S}_{\ddagger\sigma}
\end{equation}
for $j\in\{0,\ldots,n\}$, hence \eqref{defaction} is also independent of
the choice of extension $\widetilde{g}$ of $g$ when $X=T_j$ ($j=0,\ldots,n$).
Since $T_j$ ($j=0,\ldots,n$)
and the $p\in\mathcal{A}$ generate 
$\mathcal{H}$ as an algebra, a straightforward
inductive argument shows that \eqref{defaction} is well defined
for all $X\in \mathcal{H}$.  
\end{proof}

The standard basis of $\mathcal{F}_0(\mathcal{S}_{\ddagger\sigma})$
consists of the delta-functions
$\delta_\mu=\delta_\mu^\alpha\in \mathcal{F}_0(\mathcal{S}_{\ddagger\sigma})$ 
for $\mu\in\Lambda$, which are defined by
\begin{equation*}
\delta_\mu(s_\lambda^{\ddagger\sigma})=
\begin{cases}
1,\quad &\hbox{ if }\,\,\, \lambda=\mu,\\
0,\quad &\hbox{ if }\,\,\, \lambda\not=\mu.
\end{cases}
\end{equation*}
\begin{lem}\label{Fcyclic}
The subspace $\mathcal{F}_0(\mathcal{S}_{\ddagger\sigma})\subset 
\mathcal{F}(\mathcal{S}_{\ddagger\sigma})$ of functions with finite
support is a cyclic $\mathcal{H}$-submodule of 
$\mathcal{F}(\mathcal{S}_{\ddagger\sigma})$, with cyclic vector
$\delta_0\in\mathcal{F}_0(\mathcal{S}_{\ddagger\sigma})$.
\end{lem}
\begin{proof}
By \eqref{Tdot} it is clear that
$\mathcal{F}_0(\mathcal{S}_{\ddagger\sigma})\subset 
\mathcal{F}(\mathcal{S}_{\ddagger\sigma})$ is an
$\mathcal{H}$-submodule.

Since $\mathcal{W}$ 
acts transitively on $\Lambda$ under the dot-action, we may define the height 
$h(\lambda)\in\mathbb{N}$ of 
$\lambda\in\Lambda$ to be smallest nonnegative integer $m$ such that
$w\cdot\lambda=0$ for some $w\in\mathcal{W}$ of length 
$m$. To complete the proof of the lemma we show 
that $\delta_\lambda\in \mathcal{H}\cdot \delta_0
\subseteq \mathcal{F}_0(\mathcal{S}_{\ddagger\sigma})$ by 
induction to the height of $\lambda\in\Lambda$. 

It suffices to prove the induction step. Let $m\in\mathbb{Z}_{>0}$
and assume that $\delta_\mu\in\mathcal{H}\cdot\delta_0$ for all
$\mu\in\Lambda$ with $h(\mu)<m$. Let $\lambda\in\Lambda$ 
with $h(\lambda)=m$. Choose a $w\in\mathcal{W}$ satisfying $l(w)=m$
and $w\cdot\lambda=0$. Let $w=r_{i_1}r_{i_2}\cdots r_{i_m}$
be a reduced expression. We define inductively
\[\lambda_j:=r_{i_j}\cdot\lambda_{j+1}\in\Lambda,\qquad j=1,\ldots,m
\]
starting with $\lambda_{m+1}:=\lambda$. Observe that $\lambda_1=0$,
and that the $\lambda_i$ ($i=1,\ldots,m$)
are pair-wise different. Indeed, if the $\lambda_i$ are 
not pair-wise different, then $u\cdot\lambda=0$ for some $u\in \mathcal{W}$
with length strictly smaller than $h(\lambda)$, which is a
contradiction. By \eqref{Tdot} it follows that
\[
T_{w^{-1}}\cdot\delta_0=\bigl(T_{i_m}T_{i_{m-1}}\cdots
T_{i_1}\bigr)\cdot \delta_0
=c_\lambda\delta_\lambda+
\sum_{\stackrel{\scriptstyle{\mu\in\Lambda:}}{h(\mu)<h(\lambda)}}
c_{\mu}\delta_\mu
\]
for some constants $c_{\mu}\in\mathbb{C}$, with leading coefficient
given explicitly by
\[c_{\lambda}=t_w^{-1}\prod_{j=1}^mc_{i_j}(s_{\lambda_{j+1}}^{\ddagger\sigma}),
\]
cf. \cite[Lem. 9.2]{St1}.
Now $c_\lambda\not=0$ by Lemma \ref{actioncompatibleW}{\bf b)},
since $\lambda_i=r_{i_j}\cdot\lambda_{j+1}\not=\lambda_{j+1}$
for $j=1,\ldots,m$. Hence
$\delta_\lambda\in\mathcal{H}\cdot\delta_0$ by the induction hypothesis.
\end{proof}


\subsection{The kernel}

The algebra of Laurent polynomials $\mathcal{A}$ decomposes in common
eigenspaces of $Y_i\in\mathcal{H}$,
\[ \mathcal{A}=\bigoplus_{s\in\mathcal{S}}\mathcal{A}(s)
\]
with $\mathcal{A}(s)$ for $s\in\mathcal{S}=\mathcal{S}_\alpha$
the subspace
\begin{equation*}
\mathcal{A}(s)=\{p\in\mathcal{A} \, | \, r(Y)p=r(s)p,\,\,
\forall\,r\in\mathcal{A}\}.
\end{equation*}
Here $r(Y)\in\mathcal{H}$ for $r\in\mathcal{A}$
is obtained by replacing the variables $x_1,\ldots,x_n$ by
the invertible, commuting operators $Y_1,\ldots,Y_n\in\mathcal{H}_\alpha$.
The common eigenspaces $\mathcal{A}(s)$ are
one-dimensional.  
We fix a unique eigenfunction $E(s;\cdot)=E_\alpha(s;\cdot)
\in \mathcal{A}(s)$ for $s\in\mathcal{S}$ by requiring
the normalization
\begin{equation}\label{normalizationE}
 E(s;s_0^{\ddagger\sigma})=1.
\end{equation}
In particular, $E(s_0;\cdot)=1\in\mathcal{A}$ is the Laurent polynomial
identically equal to one.
The Laurent polynomials $\{E(s;\cdot) \, | \, s\in\mathcal{S}\}$
are called the Macdonald-Koornwinder polynomials,
cf. \cite{Sa}. 
We now define a kernel  
\[
\mathfrak{E}_{\mathcal{A}}(\cdot,\cdot)=
\mathfrak{E}_{\mathcal{A},\alpha}(\cdot,\cdot):
\mathcal{S}_\ddagger\times (\mathbb{C}^\times)^n\rightarrow
\mathbb{C}
\]
by
\[ \mathfrak{E}_{\mathcal{A},\alpha}(s,x)=
E_\alpha(s^{-1};x),
\qquad s\in\mathcal{S}_\ddagger,\,\, x\in (\mathbb{C}^\times)^n.
\]
\begin{prop}\label{allrelations}
For $X\in \mathcal{H}$ and $s\in\mathcal{S}_\ddagger$ we have
\[ \bigl(X\mathfrak{E}_{\mathcal{A}}(s,\cdot)\bigr)(x)=
\bigl(\psi(X)\cdot\mathfrak{E}_{\mathcal{A}}(\cdot,x)\bigr)(s)
\]
where $\mathcal{H}$ acts on the left hand side as $q$-difference
reflection operators.
\end{prop}
\begin{proof}
This follows from \cite[Prop. 7.8]{St1} and from the definition
of $\psi$, using formula \eqref{Tdot}.
\end{proof}
The duality between the spectral and
geometric parameter of the kernel $\mathfrak{E}_{\mathcal{A}}$
now reads as follows.
\begin{thm}
For all $s\in\mathcal{S}_\ddagger$ and
$v\in\mathcal{S}_{\ddagger\sigma}$,
\begin{equation}\label{polduality}
\mathfrak{E}_\mathcal{A}(s,v)=\mathfrak{E}_{\mathcal{A},\sigma}(v,s).
\end{equation}
\end{thm}
\begin{proof}
See Sahi \cite[Thm. 7.4]{Sa}.
\end{proof}


\subsection{The bilinear forms}
In the remainder of this section we choose
$0<q^{\frac{1}{2}}<1$ and $0<t\leq 1$ arbitrarily, and add generic parameters
$t_0,u_0,t_n,u_n\in\mathbb{C}^\times$ to obtain a generic difference
multiplicity function
\[\alpha=(\mathbf{t},q^{\frac{1}{2}})=(t_0,u_0,t_n,u_n,t,q^{\frac{1}{2}}).
\]
We construct suitable bilinear forms 
\[
\bigl( .,. \bigr)_{\mathcal{A}}=
\bigl(\cdot,\cdot\bigr)_{\mathcal{A},\alpha}: \mathcal{A}\times
\mathcal{A}\rightarrow \mathbb{C}
\]
and
\[ \lbrack .,. \rbrack_{\mathcal{A}}=
\lbrack\cdot,\cdot\rbrack_{\mathcal{A},\alpha}:
\mathcal{F}_0(\mathcal{S}_{\ddagger\sigma})\times 
\mathcal{F}_0(\mathcal{S}_{\ddagger\sigma})
\rightarrow \mathbb{C}
\]
satisfying the desired transformation properties 
\[\bigl(Xp,r\bigr)_{\mathcal{A}}=\bigl(p,\ddagger(X)r\bigr)_{\mathcal{A}},
\qquad\quad
\forall\,p,r\in\mathcal{A}, 
\]
respectively
\[\lbrack X\cdot g,h\rbrack_{\mathcal{A}}=\lbrack
g,\iota(X)\cdot h\rbrack_{\mathcal{A}},
\qquad\quad \forall\,g,h\in\mathcal{F}_0(\mathcal{S}_{\ddagger\sigma})
\]
under the action of $X\in\mathcal{H}$.
We start with the definition of $\bigl(\cdot,\cdot\bigr)_{\mathcal{A}}$.
Let $\Delta=\Delta_{\alpha}\in\mathcal{M}$ be the weight function
\begin{equation}\label{weight}
\Delta(x)=\prod_{f\in R_+}\frac{1}{c_f(x)}.
\end{equation}
The fact that $\Delta(x)$ is well defined and meromorphic follows
from the fact that $|q|<1$.

In case that the moduli of the Askey-Wilson parameters
\begin{equation}\label{AWpar}
\{a,b,c,d\}=\{t_nu_n,t_nu_n^{-1},q^{\frac{1}{2}}t_0u_0,
-q^{\frac{1}{2}}t_0u_0^{-1}\}.
\end{equation}
are strictly less than one, we can define the bilinear
form $\bigl(\cdot,\cdot\bigr)_{\mathcal{A}}$ by
\[ \bigl(p,r\bigr)_{\mathcal{A}}=\frac{1}{(2\pi i)^n}
\underset{\mathbb{T}^n}{\iint}p(x)r(x^{-1})\Delta(x)\frac{dx}{x},
\qquad p,r\in\mathcal{A},
\]
where $\mathbb{T}$ is the counterclock-wise oriented
unit circle in the complex plane (centered at zero), and
$\frac{dx}{x}=\frac{dx_1}{x_1}\cdots \frac{dx_n}{x_n}$ is
the standard product measure on $\mathbb{T}^n$.
For the bilinear form $\bigl(\cdot,\cdot\bigr)_{\mathcal{A}}$
in case that some of the Askey-Wilson parameters
have moduli larger than one, one needs to integrate over $\mathcal{T}^n$, 
where $\mathcal{T}=\mathcal{T}_\alpha\subset
\mathbb{C}$ is a suitable deformation of the unit circle
$\mathbb{T}$, in order to avoid certain poles of the weight function 
$\Delta$. For a detailed discussion we refer to the paper \cite{St0},
where $\mathcal{T}$ is called a $\mathbf{t}$-contour. 
Here we only
give the basic properties of such a deformed contour $\mathcal{T}$:
it is assumed to be a rectifiable, closed, counterclock-wise oriented
contour around the origin, which satisfies
$\mathcal{T}^{-1}=\mathcal{T}$ (set-theoretically), and for which the
Askey-Wilson parameters $a,b,c$ and $d$ are contained in the interior of
$\mathcal{T}$.
The following result follows from \cite[Prop. 8.3]{St1}.
\begin{prop}\label{adjoint1}
For all $p,r\in\mathcal{A}$
and $X\in\mathcal{H}=\mathcal{H}_\alpha$,
\[\bigl(Xp,r\bigr)_{\mathcal{A}}=
\bigl(p,\ddagger(X)r\bigr)_{\mathcal{A}}.
\]
\end{prop}
\begin{rem}\label{adjoint1remark}
Proposition \ref{adjoint1} is also valid for
$p$ and $r$ being arbitrary analytic functions on the
complex torus $(\mathbb{C}^\times)^n$.
\end{rem}

It follows from Proposition \ref{adjoint1} that
\begin{equation}\label{ortho}
 \bigl(E(v;\cdot),E_\ddagger(s^{-1};\cdot)\bigr)_{\mathcal{A}}=0,\qquad
v,s\in \mathcal{S},\,\,\, v\not=s
\end{equation}
since $\ddagger(Y_i)=(Y_i^\ddagger)^{-1}$ for $i=1,\ldots,n$.

Next we proceed by introducing $\lbrack
\cdot,\cdot\rbrack_{\mathcal{A}}$. 
The Weyl group $W_0\simeq S_n\ltimes (\pm 1)^n$ acts on $\Lambda$ by
permutations and sign changes of the coordinates.
For $\lambda\in\Lambda$ we write $w_{\lambda}\in S_n\ltimes (\pm 1)^n$
for the element of minimal length such that
$w_{\lambda}^{-1}\lambda\in\Lambda^+$,
where
\[ \Lambda^+=\{\lambda=\sum_i\lambda_i\epsilon_i \, | \,
\lambda_1\geq\lambda_2\geq\cdots\geq\lambda_n\geq 0\}
\]
are the partitions of length $\leq n$. Let $u_{\lambda}$ be the
$S_n$-component of $w_{\lambda}$. Let
$n_{\lambda}$ be the number of parts $\lambda_i$ of $\lambda$
which are strictly smaller than zero.
The discrete weight function
$N=N_{\alpha}: \mathcal{S}_{\ddagger\sigma}\rightarrow \mathbb{C}$
is now defined by
\begin{equation}
N(s_{\lambda}^{\ddagger\sigma})=
\underset{x=s_{\lambda}^{\ddagger\sigma}}{\hbox{\bf{Res}}}
\left(\frac{\Delta(x)}{x_1\cdots x_n}\right),\qquad
\lambda\in\Lambda,
\end{equation}
where the multiple residue {\bf{Res}} is given by
\[
\underset{x=s_{\lambda}^{\ddagger\sigma}}
{\hbox{\bf{Res}}}\bigl(\cdot\bigr)
=(-1)^{n_{\lambda}}
\underset{x_{u_\lambda(1)}=s_{\lambda,u_{\lambda}(1)}^{\ddagger\sigma}}
{\hbox{Res}}\left(
\underset{x_{u_\lambda(2)}=s_{\lambda,u_{\lambda}(2)}^{\ddagger\sigma}}
{\hbox{Res}}\left(\cdots
\underset{x_{u_\lambda(n)}=s_{\lambda,u_{\lambda}(n)}^{\ddagger\sigma}}
{\hbox{Res}}
\Bigl(\cdot\Bigr)\cdots\right)\right).
\]
The bilinear form $\lbrack \cdot,\cdot\rbrack_{\mathcal{A}}=
\lbrack \cdot,\cdot\rbrack_{\mathcal{A},\alpha}$
is defined by
\begin{equation}
\lbrack g,h\rbrack_{\mathcal{A}}=\sum_{s\in\mathcal{S}_{\ddagger\sigma}}
g(s)h(s)N(s),\qquad g,h\in\mathcal{F}_0(\mathcal{S}_{\ddagger\sigma}).
\end{equation}
Note that the definition of the bilinear form 
also makes sense for 
arbitrary functions $g,h\in\mathcal{F}(\mathcal{S}_{\ddagger\sigma})$,
provided that the sum is absolutely convergent. 
The following proposition follows now from the
proof of \cite[Prop. 8.9]{St1}.
\begin{prop}\label{adjoint2}
For all $X\in \mathcal{H}$ and all
$g,h\in\mathcal{F}_0(\mathcal{S}_{\ddagger\sigma})$, 
\begin{equation}\label{adjointdiscreteversion}
\lbrack Xg,h\rbrack_{\mathcal{A}}=\lbrack
g,\iota(X)h\rbrack_{\mathcal{A}} .
\end{equation}
\end{prop}
\begin{rem}\label{adjoint2remark}
Formula \eqref{adjointdiscreteversion} also holds true when
$g,h\in \mathcal{F}(\mathcal{S}_{\ddagger\sigma})$ as long as absolute
convergence of the sums are ensured.
\end{rem}

\subsection{The difference Fourier transforms}
We define the Macdonald-Koorwinder transform
$F_{\mathcal{A}}=F_{\mathcal{A},\alpha}: \mathcal{A}\rightarrow
\mathcal{F}_0(\mathcal{S}_\ddagger)$ by
\[ \bigl(F_{\mathcal{A}}p\bigr)(s)=
\bigl( p,\mathfrak{E}_{\mathcal{A},\ddagger}(s^{-1},\cdot)\bigr)_{\mathcal{A}},
\qquad p\in\mathcal{A}
\]
for $s\in\mathcal{S}_\ddagger$.
Furthermore, we define a linear map 
$J_{\mathcal{A}}=J_{\mathcal{A},\alpha}:
\mathcal{F}_0(\mathcal{S}_{\ddagger\sigma})
\rightarrow \mathcal{A}$ by
\[\bigl(J_{\mathcal{A}}g\bigr)(x)=
\lbrack g, \mathfrak{E}_{\mathcal{A},\sigma}(\cdot,x)\rbrack_{\mathcal{A}},
\qquad g\in\mathcal{F}_0(\mathcal{S}_{\ddagger\sigma}).
\]
The following proposition is now immediate from the previous
subsections and Section \ref{General}, see also \cite[Prop. 8.8 \&
8.9]{St1}.
\begin{prop}\label{Kcase}
{\bf a)} The Macdonald-Koornwinder transform 
$F_{\mathcal{A}}: \mathcal{A}\rightarrow \mathcal{F}_0(\mathcal{S}_\ddagger)$
is a Fourier transform associated with $\sigma$.\\ 
{\bf b)} The map $J_{\mathcal{A}}:
\mathcal{F}_0(\mathcal{S}_{\ddagger\sigma})\rightarrow\mathcal{A}$
is a Fourier transform associated with $\sigma_\sigma^{-1}$.
\end{prop}
We mention as an immediate corollary the following result.
\begin{thm}\label{corpoly}
{\bf a)} $F_{\mathcal{A}}: \mathcal{A}\rightarrow
\mathcal{F}_0(\mathcal{S}_\ddagger)$ is a linear bijection with
inverse $c_{\mathcal{A}}^{-1}J_{\mathcal{A},\sigma}$, 
with the non-zero constant $c_{\mathcal{A}}$ given by
\[ c_{\mathcal{A}}=\bigl(1,1\bigr)_{\mathcal{A}}N_\sigma(s_0^{-1}).
\]
{\bf b)} We have 
\[\frac{\bigl(E(s;\cdot),E_\ddagger(s^{-1};\cdot)\bigr)_{\mathcal{A}}}
{\bigl(1,1\bigr)_{\mathcal{A}}}=\frac{N_\sigma(s_0^{-1})}{N_\sigma(s^{-1})}
\]
for all $s\in\mathcal{S}$.
\end{thm}
\begin{proof}
See \cite[Thm. 8.10]{St1}. Since the proof is illustrative
for later arguments, we shortly recall the proof.\\
{\bf a)} By the orthogonality relations \eqref{ortho}
of the Macdonald-Koornwinder polynomials, we have
\[ F_{\mathcal{A}}(1)=\bigl(1,1\bigr)_{\mathcal{A}}\delta_0^\sigma.
\]
On the other hand, since
$\mathfrak{E}_{\mathcal{A}}(s_0^{-1},x)=E(s_0;x)=1$, 
\[ J_{\mathcal{A},\sigma}(\delta_0^\sigma)=N_\sigma(s_0^{-1})\,1\in\mathcal{A}.
\]
The statement now
now immediately follows from Proposition \ref{Kcase},
since $\mathcal{A}$ (respectively $\mathcal{F}_0(\mathcal{S}_\ddagger)$)
is a cyclic $\mathcal{H}$-module (respectively cyclic
$\mathcal{H}_\sigma$-module) 
with cyclic vector $1$ (respectively $\delta_0^\sigma$).\\
{\bf b)} This follows by computing the right hand side of
\[c_{\mathcal{A}}E(s;\cdot)=
J_{\mathcal{A},\sigma}\bigl(F_{\mathcal{A}}(E(s;\cdot))\bigr),\qquad
s\in\mathcal{S}
\]
directly from the definitions of $F_{\mathcal{A}}$ and
$J_{\mathcal{A},\sigma}$, using the orthogonality relations
\eqref{ortho}.
\end{proof}

One can reformulate the orthogonality relations in terms of an
algebraic Plan\-che\-rel type Theorem
by introducing two additional transforms
$\widetilde{F}_{\mathcal{A}}=\widetilde{F}_{\mathcal{A},\alpha}:
\mathcal{A}\rightarrow \mathcal{F}_0(\mathcal{S}_\ddagger)$
and 
$\widetilde{J}_{\mathcal{A}}=\widetilde{J}_{\mathcal{A},\alpha}: 
\mathcal{F}_0(\mathcal{S}_{\ddagger\sigma})
\rightarrow \mathcal{A}$,
\begin{equation*}
\begin{split}
\bigl(\widetilde{F}_{\mathcal{A}}p\bigr)(s)&=
\bigl(\mathfrak{E}_{\mathcal{A}}(s,\cdot),p\bigr)_{\mathcal{A}},\\
\bigl(\widetilde{J}_{\mathcal{A}}g\bigr)(x)&=
\lbrack I\mathfrak{E}_{\mathcal{A},\ddagger\sigma}(\cdot,x),
g\rbrack_{\mathcal{A}}
\end{split}
\end{equation*}
for $p\in\mathcal{A}$, $g\in\mathcal{F}_0(\mathcal{S}_{\ddagger\sigma})$ and
$s\in\mathcal{S}_\ddagger$, where $I$ is the inversion operator
$(Ig)(s)=g(s^{-1})$ mapping $\mathcal{F}_0(\mathcal{S}_{\ddagger\sigma})$
onto $\mathcal{F}_0(\mathcal{S}_\sigma)$. 
By the previous theorem, $\widetilde{F}_{\mathcal{A}}:
\mathcal{A}\rightarrow \mathcal{F}_0(\mathcal{S}_\ddagger)$
is a linear bijection with inverse 
$c_{\mathcal{A}}^{-1}\widetilde{J}_{\mathcal{A},\sigma}: 
\mathcal{F}_0(\mathcal{S}_\ddagger) \rightarrow
\mathcal{A}$. Furthermore, $\widetilde{J}_{\mathcal{A},\sigma}$ (respectively
$\widetilde{F}_{\mathcal{A}}$) is the adjoint of $F_{\mathcal{A}}$
(respectively $J_{\mathcal{A},\sigma}$) in the sense that
\begin{equation}\label{nonPlancherel}
\begin{split}
\lbrack F_{\mathcal{A}}p,g\rbrack_{\mathcal{A},\sigma}&=
\bigl(p,\widetilde{J}_{\mathcal{A},\sigma}g\bigr)_{\mathcal{A}},\\
\bigl(J_{\mathcal{A},\sigma}g,p\bigr)_{\mathcal{A}}&=
\lbrack g,\widetilde{F}_{\mathcal{A}}p\rbrack_{\mathcal{A},\sigma}
\end{split}
\end{equation}
for all $p\in\mathcal{A}$ and $g\in
\mathcal{F}_0(\mathcal{S}_\ddagger)$. This leads to the following
Plancherel type Theorem.
\begin{cor}\label{corPlancherel}
{\bf a)} For all $p,r\in \mathcal{A}$,
\[
\lbrack F_{\mathcal{A}}p,
\widetilde{F}_{\mathcal{A}}r\rbrack_{\mathcal{A},\sigma}
=c_{\mathcal{A}}\bigl(p,r\bigr)_{\mathcal{A}}.
\]
{\bf b)} For all $g,h\in\mathcal{F}_0(\mathcal{S}_\ddagger)$,
\[
\bigl(J_{\mathcal{A},\sigma}g,
\widetilde{J}_{\mathcal{A},\sigma}h\bigr)_{\mathcal{A}}
=c_{\mathcal{A}}\lbrack g,h\rbrack_{\mathcal{A},\sigma}.
\]
\end{cor}
\begin{proof}
{\bf a)} For $p,r\in\mathcal{A}$ we compute
\[
c_{\mathcal{A}}\bigl(p,r\bigr)_{\mathcal{A}}=
\bigl(J_{\mathcal{A},\sigma}(F_{\mathcal{A}}p),r\bigr)_{\mathcal{A}}=\lbrack
F_{\mathcal{A}}p,\widetilde{F}_{\mathcal{A}}r\rbrack_{\mathcal{A}}.
\]
The proof of {\bf b)} is similar.
\end{proof}

We end this section by recalling the explicit form of the
weights $\Delta$ and $N(s)$ ($s\in\mathcal{S}_{\ddagger\sigma}$) in terms 
of $q$-shifted factorials.
The expressions take their most natural form by splitting off the 
nonsymmetric part of the weight functions $\Delta\in\mathcal{M}$ 
and $N(s)$ ($s\in\mathcal{S}_{\ddagger\sigma})$. Explicitly, we can
write for the weight function $\Delta$,
\begin{equation}\label{decomposition}
\Delta(x)=\mathcal{C}(x)\Delta^+(x),
\end{equation}
with $\mathcal{C}=\mathcal{C}_\alpha\in\mathbb{C}(x)$ and
$\Delta^+=\Delta^+_\alpha\in\mathcal{M}$ given by
\begin{equation}\label{CDform}
 \mathcal{C}(x)=\prod_{\alpha\in\Sigma_-}c_\alpha(x),\qquad\quad
\Delta^+(x)=\prod_{\stackrel{\scriptstyle{f\in R:}}{f(0)\geq 0}}
\frac{1}{c_f(x)},
\end{equation}
and $\Delta^+$ is easily seen to be $W_0$-invariant. In fact, $\Delta^+$ is
given explicitly in terms of $q$-shifted factorials by
\begin{equation}\label{Deltaform}
\begin{split}
\Delta^+(x)&=\prod_{i=1}^n\frac{\bigl(x_i^2,x_i^{-2};q\bigr)_{\infty}}
{\bigl(ax_i,ax_i^{-1},bx_i,bx_i^{-1},cx_i,cx_i^{-1},dx_i,dx_i^{-1};
q\bigr)_{\infty}}\\
&\qquad\times\prod_{1\leq i<j\leq n}
\frac{\bigl(x_ix_j,x_ix_j^{-1},x_i^{-1}x_j,x_i^{-1}x_j^{-1};q\bigr)_{\infty}}
{\bigl(t^2x_ix_j,t^2x_ix_j^{-1},t^2x_i^{-1}x_j,t^2x_i^{-1}x_j^{-1};
q\bigr)_{\infty}},
\end{split}
\end{equation}
see e.g. \cite[Lem. 3.12]{St2}.

For the discrete weight $N(s)$ ($s\in\mathcal{S}_{\ddagger\sigma})$
splitting of the nonsymmetric part gives the formula 
\begin{equation}\label{Ndecomposition}
N(s_\lambda^{\ddagger\sigma})=\mathcal{C}(s_\lambda^{\ddagger\sigma})
N^+(s_{\lambda^+}^{\ddagger\sigma})
\end{equation}
for $\lambda\in\Lambda$, where $\lambda^+\in\Lambda^+$
is the unique element in $(W_0\cdot\lambda)\cap\Lambda^+$, and 
with $N^+=N_\alpha^+$ given by
\[N^+(s_\mu^{\ddagger\sigma})=
\underset{x=s_{\mu}^{\ddagger\sigma}}{\hbox{\bf{Res}}}\left(
\frac{\Delta^+(x)}{x_1\cdots x_n}\right),\qquad 
\forall\,\mu\in\Lambda^+,
\]
see \cite[(8.17) \& (8.20)]{St1}. 
Then results in \cite{St0} (see also
\cite[Rem. 8.12]{St1}) lead to the expression
\begin{equation}\label{Nplus}
\begin{split}
\frac{N^+(s_\mu^{\ddagger\sigma})}{N^+(s_0^{\ddagger\sigma})}
&=\prod_{i=1}^n\left\{\frac{\bigl(qa^2t^{4(n-i)};q\bigr)_{2\mu_i}
\bigl(q^{-1}abcdt^{4(n-i)}\bigr)^{-\mu_i}}
{\bigl(a^2t^{4(n-i)};q\bigr)_{2\mu_i}}\right.\\
&\left.\qquad\times
\frac{\bigl(a^2t^{2(n-i)},abt^{2(n-i)},act^{2(n-i)},adt^{2(n-i)}
;q\bigr)_{\mu_i}}
{\bigl(qt^{2(n-i)},qat^{2(n-i)}/b, qat^{2(n-i)}/c, 
qat^{2(n-i)}/d;q\bigr)_{\mu_i}}\right\}\\
&\times\prod_{1\leq i<j\leq n}\left\{
\frac{\bigl(qa^2t^{2(2n-i-j)},
a^2t^{2(2n-i-j+1)};q\bigr)_{\mu_i+\mu_j}}
{\bigl(qa^2t^{2(2n-i-j-1)},a^2t^{2(2n-i-j)};q\bigr)_{\mu_i+\mu_j}}
\right.\\
&\left.\qquad\qquad\qquad\qquad\qquad\qquad\qquad\times
\frac{\bigl(qt^{2(j-i)},t^{2(j-i+1)};q\bigr)_{\mu_i-\mu_j}}
{\bigl(qt^{2(j-i-1)},t^{2(j-i)};q\bigr)_{\mu_i-\mu_j}}\right\}
\end{split}
\end{equation}
for all $\mu=\sum_i\mu_i\epsilon_i\in\Lambda^+$.


\subsection{The symmetric theory}
The {\it symmetrizer}
$C_+=C_+^\alpha\in\mathcal{H}$ is defined by
\[ C_+=\frac{1}{\sum_{w\in W_0}t_w^2}\sum_{w\in W_0}t_wT_w.
\]
Observe that $C_+\in H_0=H_0^\alpha$, where $H_0$ is the sub-algebra of
$\mathcal{H}$ generated by $T_1,\ldots,T_n$ (which is 
isomorphic to the finite Hecke algebra of type $C_n$). In fact, 
the symmetrizer $C_+$ is the idempotent of $H_0$ corresponding to the
trivial character $T_i\mapsto t_i$ ($i=1,\ldots,n$) of $H_0$.

Observe that $H_0$
consists of reflection-operators only, since the $q$-difference
operators only arise from the affine root $a_0$.
In particular, $C_+$ only depends on the values $t,t_n$ and $u_n$ of
the multiplicity function $\mathbf{t}$.

Let $\mathcal{A}_+\subset \mathcal{A}$ be the algebra of
$W_0$-invariant Laurent polynomials in $\mathcal{A}$.
Then $C_+$, acting as reflection operator,
defines a projection
\[ C_+: \mathcal{A}\rightarrow \mathcal{A}_+
\]
since $T_iC_+(g)=t_iC_+(g)$ for $i=1,\ldots,n$ and $g\in\mathcal{A}$.
Applying the symmetrizer $C_+$ leads directly
to symmetric variants of the results in the previous subsections,
mainly due to the stability of 
$C_+$ under the (anti-)isomorphisms we have encountered
so far:
\begin{equation}\label{stableC}
\begin{split}
\iota(C_+)=C_+,\qquad &\dagger(C_+)=\ddagger(C_+)=C_+^{\ddagger},\\
\tau(C_+)=C_+^\tau,\qquad &\sigma(C_+)=\psi(C_+)=C_+^\sigma.
\end{split}
\end{equation}
In particular, applying the symmetrizer $C_+\in\mathcal{H}$ 
to the Macdonald-Koornwinder
polynomials leads to the following results.
\begin{prop}\label{symmpol}
{\bf a)} For all $\lambda\in\Lambda$, we have
\begin{equation}\label{invertiblestable}
C_+E(s_\lambda;\cdot)=C_+^{\ddagger}E_{\ddagger}(s_\lambda^\ddagger;\cdot)
\end{equation}
in $\mathcal{A}_+$. Furthermore, 
the expression \eqref{invertiblestable} only depends
on the $W_0$-orbit $W_0\cdot\lambda$ of $\lambda\in\Lambda$.\\
{\bf b)} Denote $\mathcal{S}^+=\mathcal{S}_\alpha^+=
\{s_\lambda\,| \, \lambda\in\Lambda^+\}$, and define $E^+(s;\cdot)=
E^+_\alpha(s;\cdot)$ for $s\in\mathcal{S}^+$ by
\[ E^+(s;\cdot)=C_+E(s;\cdot),\qquad \forall\,s\in
\mathcal{S}^+.
\]
Then $\{ E^+(s;\cdot)\, | \, s\in\mathcal{S}^+\}$ is the unique basis of
$\mathcal{A}_+$ whose basis elements satisfy the conditions
\begin{equation*}
\begin{split}
p(Y)E^+(s;\cdot)&=p(s)E^+(s;\cdot),\qquad
\forall\,p\in\mathcal{A}_+,\\
E^+(s;s_0^\sigma)&=1.
\end{split}
\end{equation*}
{\bf c)} $E^+(s;v)=E^+_\sigma(v;s)$ for $s\in\mathcal{S}^+$
and $v\in \mathcal{S}_\sigma^+$ \textup{(}duality\textup{)}.
\end{prop}
\begin{proof}
{\bf a)} This follows from \cite[(8.12)]{St1} and from the proof
of \cite[Cor. 8.11]{St1}.\\
{\bf b)} This is well known, see e.g. \cite{N}, \cite{Sa} or \cite{St1}.\\
{\bf c)} See Sahi \cite[Thm. 7.4]{Sa}, and references therein.
\end{proof}
 
We note that the coefficient of $x^\lambda$ in the monomial expansion of
$E^+(s_\lambda;x)$ for $\lambda\in\Lambda^+$ is known explicitly
in product form, see \cite[Cor. 9.4]{St1} and references therein. 
It is known as the
evaluation Theorem. 
In our present notations it reads as follows.
\begin{thm}\label{evaluation}
Let $\lambda=\sum_i\lambda_i\epsilon_i\in\Lambda^+$. 
The coefficient $c_\lambda=c_\lambda^\alpha\in\mathbb{C}$
of $x^\lambda$ in the monomial expansion of $E^+(s_\lambda;x)$ is
given explicitly by
\begin{equation*}
\begin{split}
c_\lambda=&\prod_{i=1}^n\frac{\bigl(q^{-1}abcdt^{4(n-i)};q\bigr)_{2\lambda_i}
(at^{2(n-i)})^{\lambda_i}}
{\bigl(abt^{2(n-i)},act^{2(n-i)}, adt^{2(n-i)},
q^{-1}abcdt^{2(n-i)};q\bigr)_{\lambda_i}}\\
&\times\prod_{1\leq i<j\leq n}\frac
{\bigl(q^{-1}abcdt^{2(2n-i-j)};q\bigr)_{\lambda_i+\lambda_j}
\bigl(t^{2(j-i)};q\bigr)_{\lambda_i-\lambda_j}}
{\bigl(q^{-1}abcdt^{2(2n-i-j+1)};
q\bigr)_{\lambda_i+\lambda_j}\bigl(t^{2(j-i+1)};q\bigr)_{\lambda_i-\lambda_j}},
\end{split}
\end{equation*}
where we used the Askey-Wilson parametrization \eqref{AWpar} for
part of the multiplicity function $\alpha$.
\end{thm}

The basis elements $E^+(s;\cdot)$ ($s\in\mathcal{S}^+$) are known
as the (normalized) symmetric Macdonald-Koornwinder polynomials.
Up to an explicit multiplicative constant, they coincide with
Koornwinder's \cite{K} multivariable analogues of the 
Askey-Wilson polynomials. This can be proven by relating the
$q$-difference reflection operator
\begin{equation}\label{KAWeq}
m_{\epsilon_1}(Y):=Y_1+\cdots +Y_n+Y_1^{-1}+\cdots +Y_n^{-1}\in\mathcal{H}
\end{equation}
acting on $\mathcal{A}_+$ to Koornwinder's \cite{K} multivariable
second-order $q$-difference operator of Askey-Wilson type, see
Noumi \cite{N} (see \cite{St2} for a detailed treatment in english).
In particular, in the rank one set-up ($n=1$), $E^+(s_m;\cdot)$ for
$m\in\Lambda^+\simeq \mathbb{N}$ is the well known Askey-Wilson
polynomial \cite{AW} of degree $m$. 
In terms of the standard notation \cite{GR}
for basic hypergeometric series,
\[
{}_r\phi_s\left(\begin{matrix} a_1,a_2,\ldots,a_r\\
b_1,b_2,\ldots,b_s\end{matrix}\,;\,q,z\right)=
\sum_{k=0}^{\infty}
\frac{\bigl(a_1,a_2,\ldots,a_r;q\bigr)_k}
{\bigl(q,b_1,b_2,\ldots,b_s;q\bigr)_k}
\lbrack (-1)^kq^{\frac{1}{2}k(k-1)}\rbrack^{1+s-r}\,z^k,
\]
this leads to the explicit series expansion
\begin{equation}\label{AWpolynomial}
E^+(s_m;x)={}_4\phi_3\left(\begin{matrix} q^{-m}, q^{m-1}abcd,ax,a/x\\
ab,ac,ad\end{matrix}\,;\, q,q\right)
\end{equation}
for the symmetric, rank one Macdonald-Koornwinder polynomial, 
see \cite{NS} for a detailed account. The Askey-Wilson 
parameters $a,b,c,d$ play
a symmetrical role in both the Askey-Wilson second-order
$q$-difference operator as well as in the polynomial spectrum
$\mathcal{S}^+$. Consequently, $E^+(s_m;x)$ is, up to a
multiplicative constant, invariant under permuting the Askey-Wilson parameters
$\{a,b,c,d\}$. The corresponding identity for the balanced 
${}_4\phi_3$ is known as Sear's transformation formula, see
\cite[(2.10.4)]{GR}. This can be generalized to the higher rank
setup. We formulate
here two cases, the first corresponds to interchanging the role of
$a$ and $b$, the second corresponds to interchanging the role of $c$
and $d$.

\begin{prop}\label{symmetrypar}
Let $\{a,b,c,d\}$ be the Askey-Wilson parametrization \eqref{AWpar}
of part of the
difference multiplicity function 
$\alpha=(t_0,u_0,t_n,u_n,t,q^{\frac{1}{2}})$.\\
{\bf a)} Set $\beta=(t_0,u_0,t_n,-u_n^{-1},t,q^{\frac{1}{2}})$. Then
$s_\lambda^\beta=s_\lambda^\alpha$ for all $\lambda\in\Lambda^+$ and
\[E_\beta^+(s_\lambda;x)=
\left(\prod_{i=1}^n\frac{\bigl(act^{2(n-i)},adt^{2(n-i)};q\bigr)_{\lambda_i}}
{\bigl(bct^{2(n-i)}, bdt^{2(n-i)};q\bigr)_{\lambda_i}}
\left(\frac{b}{a}\right)^{\lambda_i}\right)E_\alpha^+(s_\lambda;x),\qquad
\forall\,\lambda\in\Lambda^+.
\]
{\bf b)} Set $\gamma=(t_0,-u_0^{-1},t_n,u_n,t,q^{\frac{1}{2}})$. Then
$s_\lambda^\gamma=s_\lambda^\alpha$ for all $\lambda\in\Lambda^+$ and
\[E_\gamma^+(s_\lambda;x)=
E_\alpha^+(s_\lambda;x),\qquad\forall\,\lambda\in\Lambda^+.
\]
\end{prop}
\begin{proof}
For $\lambda\in\Lambda^+$ we 
write $E_\alpha^+(s_\lambda;x)=c_\lambda^\alpha 
P_\lambda^\alpha(x)$ with $c_\lambda=c_\lambda^\alpha$ 
given as in Theorem \ref{evaluation}.
The Laurent polynomial
$P_\lambda$ is the unique eigenfunction of $p(Y)$ 
with eigenvalue $p(s)$ for all $p\in\mathcal{A}_+$ whose coefficient of 
$x^\lambda$ in its monomial expansion is equal to one. 

Now the operators $Y_i$ ($i=1,\ldots,n$) 
as well as the spectral points $s_\lambda\in\mathcal{S}^+$ 
($\lambda\in\Lambda^+$)
are invariant under replacement of $u_0$ by $-u_0^{-1}$  
(respectively $u_n$ by $-u_n^{-1}$). Indeed, for the $Y_i$ this
follows from the fact that the $q$-difference reflection operators
$T_j\in\mathcal{H}$ ($j=0,\ldots,n$) are invariant under replacement of 
$u_0$ by $-u_0^{-1}$  (respectively $u_n$ by $-u_n^{-1}$).
For the spectral points $s_\lambda\in\mathcal{S}^+$ this simply
follows from the fact that the spectrum is independent
of $u_0$ and of $u_n$. We conclude that  
$P_\lambda^\beta(x)=P_\lambda^\gamma(x)=P_\lambda^\alpha(x)$
for all $\lambda\in\Lambda^+$, hence
\[E_{\beta}^+(s_\lambda;x)=\frac{c_{\lambda}^\beta}{c_{\lambda}^\alpha}\,
E_\alpha^+(s_\lambda;x),\qquad
E_{\gamma}^+(s_\lambda;x)=\frac{c_{\lambda}^\gamma}{c_{\lambda}^\alpha}\,
E_\alpha^+(s_\lambda;x)
\]
for all $\lambda\in\Lambda^+$. The proposition now follows from 
the explicit expressions for the
coefficients $c_\lambda$, see Theorem \ref{evaluation}. 
\end{proof}

Let $\mathcal{F}^+_0(\mathcal{S})$ be the space of functions 
$g\in \mathcal{F}_0(\mathcal{S})$ which are $W_0$-invariant
under the dot-action. 
Observe that $\mathcal{F}_0^+(\mathcal{S})$ may be identified with
the space $\mathcal{F}_0(\mathcal{S}^+)$ of functions
$g: \mathcal{S}^+\rightarrow\mathbb{C}$ of finite support
by the natural restriction map.

The restriction of the Macdonald-Koornwinder transform
$F_{\mathcal{A}}$ to $\mathcal{A}_+$ can be written as integral
transform with $W_0$-invariant integrand as follows.
By Proposition \ref{adjoint1}, \eqref{stableC} and Proposition
\ref{symmpol} we have
\[\bigl(F_{\mathcal{A}}p\bigr)(s_\lambda^{-1})=
\bigl(p,E^+(s_{\lambda^+};\cdot)\bigr)_{\mathcal{A}},\qquad
p\in\mathcal{A}_+,\,\, \lambda\in\Lambda,
\]
where $\lambda^+$ is the unique element in $(W_0\cdot\lambda)\cap
\Lambda^+$. The restriction of the bilinear form
$\bigl(\cdot,\cdot\bigr)_{\mathcal{A}}$ to $\mathcal{A}_+$
can be rewritten using the formula
\begin{equation}\label{symmistrivial}
\sum_{w\in W_0}(w\mathcal{C})(x)=
\mathcal{C}(s_0^{\ddagger\sigma}),
\end{equation}
see \cite[Lem. 8.1 \& Lem. 8.2]{St1}, leading to the identity
\begin{equation}\label{symmcontinuous}
 \bigl(p,r\bigr)_{\mathcal{A}}=
\frac{\mathcal{C}(s_0^{\ddagger\sigma})}{2^nn!}\bigl(p,r\bigr)_{\mathcal{A},+},
\qquad p,r\in\mathcal{A}_+,
\end{equation}
with the bilinear form 
$\bigl(\cdot,\cdot\bigr)_{\mathcal{A},+}=\bigl(\cdot,\cdot
\bigr)_{\mathcal{A},\alpha,+}$ on $\mathcal{A}_+$ given by
\begin{equation}\label{pairingsymm}
\bigl(p,r\bigr)_{\mathcal{A},+}=
\frac{1}{(2\pi i)^n}\underset{\mathcal{T}^n}{\iint}
p(x)r(x)\Delta^+(x)\frac{dx}{x},\qquad p,r\in\mathcal{A}_+.
\end{equation}
In particular, the restriction of the Macdonald-Koornwinder transform
$F_{\mathcal{A}}$ to $\mathcal{A}_+$ becomes
\begin{equation}\label{FAsym}
\bigl(F_{\mathcal{A}}p\bigr)(s_\lambda^{-1})=
\frac{\mathcal{C}(s_0^{\ddagger\sigma})}{2^nn!}
\bigl(p,E^+(s_{\lambda^+};\cdot)\bigr)_{\mathcal{A},+},\qquad
p\in\mathcal{A}_+,\,\, \lambda\in\Lambda^+.
\end{equation}
Similar arguments show that
$\widetilde{F}_{\mathcal{A}}|_{\mathcal{A}_+}=
\widetilde{F}_{\mathcal{A}}|_{\mathcal{A}_+}$. In particular,
the Plancherel type formulas for the Macdonald-Koornwinder transform
$F_{\mathcal{A}}$ (see Corollary \ref{corPlancherel}) reduce to 
genuine Plancherel formulas when restricting to $W_0$-invariant
functions.

The same arguments can now be applied to handle the restriction of
the transform $J_{\mathcal{A}}:
\mathcal{F}_0(\mathcal{S}_{\ddagger\sigma})
\rightarrow \mathcal{A}$ to the subspace
$\mathcal{F}_0^+(\mathcal{S}_{\ddagger\sigma})$ of $W_0$-invariant
functions. In this case,
one needs Proposition \ref{adjoint2} and the 
discretized version 
\begin{equation}\label{symmistrivialdiscrete}
\sum_{\mu\in W_0\cdot\lambda}\mathcal{C}(s_{\mu}^{\ddagger\sigma})=
\mathcal{C}(s_0^{\ddagger\sigma})
\end{equation}
of \eqref{symmistrivial}, see e.g. \cite[Lem. 8.2]{St1} for a proof.
It leads to the formula
\begin{equation}\label{JAsym}
J_{\mathcal{A}}g=
\widetilde{J}_{\mathcal{A}}g=\mathcal{C}(s_0^{\ddagger\sigma})
\sum_{\lambda\in\Lambda^+}g(s_\lambda^{\ddagger\sigma})
E_\sigma^+(s_\lambda^\sigma;\cdot)N^+(s_\lambda^{\ddagger\sigma})
\end{equation}
for $g\in\mathcal{F}_0^+(\mathcal{S}_{\ddagger\sigma})$. 
By Theorem \ref{corpoly}{\bf a)} and formulas \eqref{FAsym},
\eqref{JAsym} and \eqref{symmistrivial}, the orthogonality
relations of the nonsymmetric Macdonald-Koornwinder polynomials
(see \eqref{ortho} and Theorem \ref{corpoly}{\bf b)}) 
reduce to the orthogonality relations 
\begin{equation}
\frac{\bigl(E^+(s;\cdot),E(v;\cdot)\bigr)_{\mathcal{A},+}}
{\bigl(1,1\bigr)_{\mathcal{A},+}}=\frac{N_\sigma^+(s_0^{-1})}
{N_\sigma^+(s^{-1})}\,
\delta_{s,v},\qquad \forall\,s,v\in\mathcal{S}^+
\end{equation}
for the symmetric Macdonald-Koornwinder polynomials,
where $\delta_{s,v}$ is the Kro\-ne\-cker delta function on $\mathcal{S}^+$.

Under suitable values of the difference multiplicity function
$\alpha$, the integration contour can be shifted to $\mathbb{T}^n$
on the cost of some extra (partly discrete) contributions to
the measure. This in particular leads to a positive orthogonality measure 
for the symmetric Macdonald-Koornwinder polynomials, see \cite{St0}
for more details. We also remark that the constant term 
$\bigl(1,1\bigr)_{\mathcal{A},+}$ is exactly the multivariable
Askey-Wilson integral, evaluated by Gustafson in \cite{G}.
Explicitly, it reads
\begin{equation}\label{Gustafson}
\begin{split}
\bigl(1,&1\bigr)_{\mathcal{A},+}=2^nn!\prod_{i=1}^n\left\{
\frac{\bigl(t^2,abcdt^{2(2n-i-1)};q\bigr)_{\infty}}
{\bigl(q,t^{2(n-i+1)};q\bigr)_{\infty}}\right. \\
&\,\left.\times\frac{1}{\bigl(
abt^{2(n-i)},act^{2(n-i)},adt^{2(n-i)},bct^{2(n-i)},
bdt^{2(n-i)},cdt^{2(n-i)};q\bigr)_{\infty}}\right\}.
\end{split}
\end{equation}


\section{The construction of the Cherednik kernels}

We assume throughout this section that 
\[\alpha=(\mathbf{t},q^{\frac{1}{2}})=(t_0,u_0,t_n,u_n,t,q^{\frac{1}{2}})
\]
is a difference multiplicity function with
$0<q^{\frac{1}{2}}<1$ and $0<t\leq 1$, and with generic parameters  
$t_0,u_0,t_n,u_n\in\mathbb{C}^\times$. 
As was explained in Section 3, the construction of non-polynomial 
Fourier transforms associated with $\sigma$ hinges on the existence 
of non-zero kernels
\begin{equation*}
\begin{split}
\mathfrak{E}&=\mathfrak{E}_\alpha: (\mathbb{C}^\times)^n\times
(\mathbb{C}^\times)^n\rightarrow \mathbb{C}\\
\mathfrak{E}_\ddagger&=\mathfrak{E}_{\alpha_\ddagger}:
(\mathbb{C}^\times)^n\times
(\mathbb{C}^\times)^n\rightarrow \mathbb{C}
\end{split}
\end{equation*}
satisfying the transformation behaviour
\begin{equation}\label{transformation1}
\begin{split}
\bigl(X\mathfrak{E}(\gamma,\cdot)\bigr)(x)&=
\bigl(\psi(X)\mathfrak{E}(\cdot,x)\bigr)(\gamma),\qquad
X\in\mathcal{H},\\
\bigl(X\mathfrak{E}_\ddagger(\gamma,\cdot)\bigr)(x)&=
\bigl(\psi_\ddagger(X)\mathfrak{E}(\cdot,x)\bigr)(\gamma),\qquad
X\in\mathcal{H}_\ddagger
\end{split}
\end{equation}
with respect to the duality anti-isomorphism $\psi$.
Here e.g. $\bigl(X\mathfrak{E}(\gamma,\cdot)\bigr)(x)$ means
the action of $X\in\mathcal{H}$ as $q$-difference reflection operator
on the function $x\mapsto \mathfrak{E}(\gamma,x)$, which also makes
sense when $\mathfrak{E}$ is assumed to be a meromorphic function on
$(\mathbb{C}^\times)^n\times (\mathbb{C}^\times)^n$, cf. Subsection 2.3.
 
If we require $\mathfrak{E}$ and $\mathfrak{E}_\ddagger$
to be meromorphic and ``as regular
as possible'', then kernels $\mathfrak{E}$ and $\mathfrak{E}_\ddagger$
satisfying \eqref{transformation1} 
turn out to be unique up to a multiplicative constant. 
We call them the {\it Cherednik kernels},
since they generalize kernels introduced by 
Cherednik in \cite{C1} to the set-up of nonreduced root systems.

We note that 
the classical approach to study existence of such a 
kernel $\mathfrak{E}(\gamma,x)=
\mathfrak{E}_\alpha(\gamma,x)$ 
is by analyzing the solution space of the spectral problem
\begin{equation}\label{spectralproblem}
 Y_ig=\gamma_i^{-1}g,\qquad i=1,\ldots,n,
\end{equation}
for the commuting $Y$-operators $Y_i\in\mathcal{H}$. 
We follow here a drastically different approach, which is
inspired by Cherednik's paper \cite{C1}. The philosophy
is to {\it construct} the special solutions $\mathfrak{E}(\gamma,x)$ 
of the spectral problem \eqref{spectralproblem} as an explicit series
expansion in terms of the Macdonald-Koornwinder polynomials.
The important step in the construction is to associate  
a so-called {\it auxiliary kernel}
$\mathfrak{F}=\mathfrak{F}_\alpha$ to 
$\mathfrak{E}$, for which the spectral problem
\eqref{spectralproblem} translates into 
a transformation property of the type
\[ \bigl(Y_i\mathfrak{F}(\gamma,\cdot)\bigr)(x)=
\bigl(Y_i\mathfrak{F}(\cdot,x)\bigr)(\gamma),\qquad i=1,\ldots,n
\]
(with the $Y$-operators on each side depending on different
multiplicity functions).
Invoking the polynomial theory, it becomes now plausible that
such an auxiliary kernel $\mathfrak{F}$ can be constructed 
as explicit series expansions involving
the Macdonald-Koornwinder polynomials in the variables $x$ as well as in 
the variables $\gamma$. This is indeed the case, as is explained in
full detail in this section.


\subsection{Auxiliary kernels}
Instead of studying kernels
$\mathfrak{E},\mathfrak{E}_\ddagger\in
\mathcal{M}\bigl((\mathbb{C}^\times)^n\times (\mathbb{C}^\times)^n
\bigr)$ satisfying the transformation behaviour \eqref{transformation1}
directly, it turns out to be convenient to study certain
auxiliary kernels
$\mathfrak{F}=\mathfrak{F}_\alpha, \mathfrak{F}_\ddagger=
\mathfrak{F}_{\alpha_\ddagger}\in\mathcal{M}
\bigl((\mathbb{C}^\times)^n\times(\mathbb{C}^\times)^n\bigr)$ first,
which are related to $\mathfrak{E}$ and $\mathfrak{E}_\ddagger$
by the formulas
\begin{equation}\label{link1}
\begin{split}
\mathfrak{F}(\gamma,x)&=G_{\sigma\tau}(\gamma)^{-1}
G_\tau(x)^{-1}\mathfrak{E}(\gamma,x),\\
\mathfrak{F}_\ddagger(\gamma,x)&=
G_{\ddagger\sigma}(\gamma)G_{\ddagger}(x)
\mathfrak{E}_\ddagger(\gamma^{-1},x),
\end{split}
\end{equation}
where $G=G_\alpha\in\mathcal{M}$ is the Gaussian.
The kernel $\mathfrak{F}_\ddagger$ may be alternatively written as
\[
\mathfrak{F}_\ddagger(\gamma,x)=G_{\sigma\tau}(\gamma)^{-1}G_\tau(x)^{-1}
\mathfrak{E}_\ddagger(\gamma^{-1},x)
\]
in view of \eqref{Gaussianinversion}.

\begin{lem}\label{EFconnection}
Let $\kappa=\kappa_\alpha: \mathcal{H}\rightarrow
\mathcal{H}_{\tau\sigma\tau}$ and 
$\kappa^I=\kappa^I_\alpha: \mathcal{H}\rightarrow
\mathcal{H}_{\ddagger\tau\sigma\tau}$ be unital anti-isomorphisms
defined by
\begin{equation}\label{kappnew}
\kappa=\tau_{\tau\sigma\tau}^{-1}\circ\psi_\tau\circ\tau,\qquad
\kappa^I=\dagger_{\tau\sigma\tau}\circ\tau_{\tau\sigma}\circ
\psi_\tau\circ\tau_\tau^{-1}.
\end{equation}
Under the correspondence \eqref{link1}, the transformation behaviour
\eqref{transformation1} for the meromorphic kernels $\mathfrak{E}$
and $\mathfrak{E}_\ddagger$ are equivalent to the transformation
behaviour
\begin{equation}\label{transformationauxnew1}
\bigl(X\mathfrak{F}(\gamma,\cdot)\bigr)(x)=
\bigl(\kappa_\tau(X)\mathfrak{F}(\cdot,x)\bigr)(\gamma),
\qquad X\in\mathcal{H}_\tau,
\end{equation}
\begin{equation}\label{transformationauxnew2}
\bigl(X\mathfrak{F}_{\ddagger}(\gamma,\cdot)\bigr)(x)=
\bigl(\kappa^I_{\ddagger\tau}(X)\mathfrak{F}_{\ddagger}(\cdot,x)
\bigr)(\gamma),\qquad X\in\mathcal{H}_{\ddagger\tau}
\end{equation}
for the auxiliary 
meromorphic kernels $\mathfrak{F}$ and $\mathfrak{F}_\ddagger$.
\end{lem}
\begin{proof}
This follows from the fact that conjugation by the Gaussian
induces the isomorphism $\tau$ on $\mathcal{H}$ (see Proposition
\ref{Gaussianconj}), 
and that the inversion operator
$(Ig)(\gamma)=g(\gamma^{-1})$ induces the isomorphism $\dagger$ 
on $\mathcal{H}$. As an example, we prove the transformation behaviour
for $\mathfrak{F}_\ddagger$ under the assumption that
$\mathfrak{E}_\ddagger$ satisfies \eqref{transformation1}.
We use the formula
\[\mathfrak{F}_\ddagger(\gamma,x)=G_\ddagger(x)
\bigl(I\bigl(G_{\ddagger\sigma}(\cdot)
\mathfrak{E}_\ddagger(\cdot,x)\bigr)\bigr)(\gamma),
\]
which follows from \eqref{link1}
since the Gaussian is $W_0$-invariant.
Then we compute for $X\in\mathcal{H}_{\ddagger\tau}$,
\begin{equation*}
\begin{split}
\bigl(X\mathfrak{F}_\ddagger(\gamma,\cdot)\bigr)(x)&=
G_\ddagger(x)
\bigl(I\bigl(G_{\ddagger\sigma}(\cdot)
\bigl(\psi_\ddagger\circ\tau_\ddagger^{-1}\bigr)(X)
\mathfrak{E}_\ddagger(\cdot,x)\bigr)\bigr)(\gamma)\\
&=G_\ddagger(x)
\bigl(I\bigl(
\bigl(\tau_{\ddagger\sigma}\circ\psi_\ddagger\circ\tau_\ddagger^{-1}\bigr)(X)
G_{\ddagger\sigma}(\cdot)
\mathfrak{E}_\ddagger(\cdot,x)\bigr)\bigr)(\gamma)\\
&=\bigl(\bigl(\dagger_{\ddagger\sigma\tau}\circ
\tau_{\ddagger\sigma}\circ\psi_\ddagger\circ\tau_\ddagger^{-1}\bigr)(X)
\mathfrak{F}_\ddagger(\cdot,x)\bigr)(\gamma)\\
&=\bigl(\kappa_{\ddagger\tau}^I(X)\mathfrak{F}_\ddagger(\cdot,x)\bigr)(\gamma),
\end{split}
\end{equation*}
which is the desired result.
\end{proof}

In the next lemma we 
compute the anti-isomorphisms $\kappa$ and $\kappa^I$
explicitly on suitable algebraic generators of $\mathcal{H}$. In particular,
we show that these anti-isomorphisms have the desired 
property that $Y$-operators are mapped to $Y$-operators.

We denote
\begin{equation}\label{U}
U=U_{\alpha}=T_1\cdots T_{n-1}T_nT_{n-1}\cdots T_1\in\mathcal{H},
\end{equation}
so that $Y_1=UT_0$ in the double affine Hecke algebra
$\mathcal{H}$.
\begin{lem}\label{kappalem}
{\bf a)} For $i=1,\ldots,n$ we have
\[\kappa(T_i)=T_i^{\tau\sigma\tau},\qquad
\kappa(Y_i)=Y_i^{\tau\sigma\tau}.
\]
Furthermore, $\kappa(T_0)=U_{\tau\sigma\tau}T_0^{\tau\sigma\tau}
U_{\tau\sigma\tau}^{-1}$ and
$\kappa(x_1)=q^{\frac{1}{2}}x_1^{-1}T_0^{\tau\sigma\tau}
U_{\tau\sigma\tau}^{-1}$.\\
{\bf b)} For $i=1,\ldots,n$ we have
\[\kappa^I(T_i)=T_i^{\ddagger\tau\sigma\tau}{}^{-1},\qquad
\kappa^I(Y_i)=Y_i^{\ddagger\tau\sigma\tau}{}^{-1}.
\]
Furthermore, $\kappa^I(T_0)=T_0^{\ddagger\tau\sigma\tau}{}^{-1}$
and $\kappa^I(x_1)=q^{\frac{1}{2}}T_0^{\ddagger\tau\sigma\tau}{}^{-1}
x_1U_{\ddagger\tau\sigma\tau}$.
\end{lem}
\begin{proof}
{\bf a)} The identities $\kappa(T_i)=T_i^{\tau\sigma\tau}$
for $i=1,\ldots,n$ are immediately clear from the definitions 
of $\tau$ and $\psi$. Since $Y_{i+1}=T_i^{-1}Y_iT_i^{-1}$
for $i=1,\ldots,n-1$, the identities $\kappa(Y_i)=Y_i^{\tau\sigma\tau}$
for $i=1,\ldots,n$ will follow from $\kappa(Y_1)=Y_1^{\tau\sigma\tau}$.
Omitting the parameter dependencies we compute,
\begin{equation*}
\begin{split}
\kappa(Y_1)&=\tau^{-1}\psi\tau(UT_0)\\
&=\tau^{-1}\psi(q^{-\frac{1}{2}}Ux_1T_0^{-1})\\
&=\tau^{-1}\psi(q^{-\frac{1}{2}}Ux_1Y_1^{-1}U)\\
&=\tau^{-1}(q^{-\frac{1}{2}}Ux_1Y_1^{-1}U)\\
&=\tau^{-1}(q^{-\frac{1}{2}}Ux_1T_0^{-1})=UT_0=Y_1,
\end{split}
\end{equation*}
which is the desired result. Furthermore, omitting the 
parameter dependencies,
\[\kappa(T_0)U=\kappa(UT_0)=\kappa(Y_1)=UT_0,
\]
hence $\kappa(T_0)=
U_{\tau\sigma\tau}T_0^{\tau\sigma\tau}U_{\tau\sigma\tau}^{-1}$.
The identity 
$\kappa(x_1)=q^{\frac{1}{2}}x_1^{-1}T_0^{\tau\sigma\tau}
U_{\tau\sigma\tau}^{-1}$ follows from similar (but easier)
computations.\\
{\bf b)} The identities $\kappa^I(T_i)=
T_i^{\ddagger\tau\sigma\tau}{}^{-1}$ for $i=1,\ldots,n$
are clear. For the identities 
$\kappa^I(Y_i)=Y_i^{\ddagger\tau\sigma\tau}{}^{-1}$, it then
suffices to prove that
$\kappa^I(T_0)=T_0^{\ddagger\tau\sigma\tau}{}^{-1}$. 
Omitting the parameter dependencies
we compute,
\begin{equation*}
\begin{split}
\kappa^I(T_0)&=\dagger\tau\psi\tau^{-1}(T_0)\\
&=\dagger\tau\psi(q^{-\frac{1}{2}}T_0^{-1}x_1)\\
&=\dagger\tau\psi(q^{-\frac{1}{2}}Y_1^{-1}Ux_1)\\
&=\dagger\tau(q^{-\frac{1}{2}}Y_1^{-1}Ux_1)\\
&=\dagger\tau(q^{-\frac{1}{2}}T_0^{-1}x_1)\\
&=\dagger(T_0)=T_0^{-1},
\end{split}
\end{equation*}
which is the desired result. The proof of the identity
$\kappa^I(x_1)=
q^{\frac{1}{2}}
T_0^{\ddagger\tau\sigma\tau}{}^{-1}x_1U_{\ddagger\tau\sigma\tau}$ 
is left to the reader.
\end{proof}

\subsection{Auxiliary transforms}
Recall that $\mathcal{O}=
\mathcal{O}\bigl((\mathbb{C}^\times)^n\bigr)$ is the ring of
analytic functions on the complex torus $(\mathbb{C}^\times)^n$.
\begin{Def}
{\bf a)} A non-zero linear mapping $L=L_{\alpha}:
\mathcal{A}\rightarrow \mathcal{O}$ is called an auxiliary
transform associated to $\alpha$ when
\begin{equation}\label{aux1}
 L\circ X=\bigl(\kappa_{\tau}\circ\ddagger_{\tau}^{-1}\bigr)(X)\circ
L,\qquad \forall\, X\in\mathcal{H}_{\ddagger\tau},
\end{equation}
with the usual action of the double affine Hecke algebras on the
function spaces $\mathcal{A}\subset \mathcal{O}\subset \mathcal{M}$.\\
{\bf b)} A non-zero linear mapping $L_\ddagger=L_{\alpha_\ddagger}:
\mathcal{A}\rightarrow \mathcal{O}$ is called an auxiliary
transform associated to $\alpha_\ddagger$ when
\begin{equation}\label{aux2} 
L_\ddagger\circ X=
\bigl(\kappa_{\ddagger\tau}^I\circ \ddagger_\tau\bigr)(X)\circ 
L_\ddagger,\qquad \forall\,X\in\mathcal{H}_\tau,
\end{equation}
with the usual action of the double affine Hecke algebras on the
function spaces $\mathcal{A}\subset \mathcal{O}\subset \mathcal{M}$.
\end{Def} 

In the following lemma we link auxiliary transforms to {\it analytic} 
kernels $\mathcal{F}$ and $\mathcal{F}_\ddagger$
satisfying the transformation behaviour \eqref{transformationauxnew1}
and \eqref{transformationauxnew2}, respectively.
\begin{lem}\label{Auxiliarylink}
Suppose that $L,L_\ddagger:\mathcal{A}\rightarrow\mathcal{O}$
are linear mappings of the form
\[\bigl(Lp\bigr)(\gamma)=\bigl(\mathfrak{F}(\gamma,\cdot),
p\bigr)_{\mathcal{A},\tau},\qquad
\bigl(L_{\ddagger}p\bigr)(\gamma)=
\bigl(p,\mathfrak{F}_\ddagger(\gamma,\cdot)\bigr)_{\mathcal{A},\tau}
\]
for $p\in\mathcal{A}$, with 
kernels $\mathfrak{F},\mathfrak{F}_\ddagger\in
\mathcal{O}\bigl((\mathbb{C}^\times)^n\times
(\mathbb{C}^\times)^n\bigr)$ \textup{(}which are then necessarily
unique\textup{)}.\\
{\bf a)} $L$ 
is an auxiliary transform associated to $\alpha$ if and only if 
$\mathfrak{F}$ satisfies the transformation behaviour
\eqref{transformationauxnew1}.\\
{\bf b)} $L_\ddagger$ is an auxiliary transform associated to
$\alpha_\ddagger$ if and only if $\mathfrak{F}_\ddagger$
satisfies the transformation behaviour \eqref{transformationauxnew2}.
\end{lem}
\begin{proof}
By Proposition \ref{adjoint1} and Remark \ref{adjoint1remark} we have
that $\bigl(\cdot,\cdot\bigr)_{\mathcal{A},\tau}$ induces
the anti-isomorphism $\ddagger_\tau$ on $\mathcal{H}_\tau$, i.e.
\[\bigl(Xp,r\bigr)_{\mathcal{A},\tau}=\bigl(p,\ddagger_\tau(X)r
\bigr)_{\mathcal{A},\tau},\qquad X\in\mathcal{H}_\tau
\]
when $p,r\in\mathcal{O}$. The lemma is now immediate.
\end{proof}

In this and the next subsection we study auxiliary transforms
in detail. We prove that they exist and that
they are unique up to a 
multiplicative constant. Furthermore, we explicitly compute
the auxiliary transforms on
suitable bases of $\mathcal{A}$. This leads, via the previous
lemma, to the explicit construction of {\it analytic} kernels
$\mathfrak{F}$ and $\mathfrak{F}_\ddagger$ satisfing
the transformation behaviour \eqref{transformationauxnew1}
and \eqref{transformationauxnew2}, respectively.

We have the following key lemma.
\begin{lem}\label{L}
{\bf a)} If an auxiliary transform $L: \mathcal{A}\rightarrow
\mathcal{O}$ associated to $\alpha$ exists, 
then it is unique up to a multiplicative constant. Furthermore,
$L(1)=c 1$ for some constant
$c\in\mathbb{C}^\times$, 
where $1\in\mathcal{A}$ is the Laurent polynomial
identically equal to one, and $L(\mathcal{A})\subseteq\mathcal{A}$.\\
{\bf b)} If an auxiliary transform $L_{\ddagger}: \mathcal{A}
\rightarrow \mathcal{O}$ associated to $\alpha_\ddagger$ exists,
then it is unique up to a multiplicative constant.
Furthermore,
$L_\ddagger(1)=c_\ddagger 1$ for some constant
$c_\ddagger\in\mathbb{C}^\times$ 
and $L_\ddagger(\mathcal{A})\subseteq\mathcal{A}$.
\end{lem}
\begin{proof}
Let $L,L_\ddagger: \mathcal{A}\rightarrow \mathcal{O}$ be non-zero
linear maps satisfying the transformation property \eqref{aux1} and
\eqref{aux2}, respectively.
Since $\mathcal{A}$ is a cyclic module for the action of the double
affine Hecke algebra with cyclic vector $1\in\mathcal{A}$, the lemma
follows from the formulas
$L(1)=c 1$ and $L_\ddagger(1)=c_\ddagger 1$
for some constants $c,c_\ddagger\in\mathbb{C}$.

We first show that the functions
$L(1), L_\ddagger(1)\in \mathcal{O}$ 
are $W_0$-invariant.
We focus on $L(1)$ 
(the proof for $L_\ddagger(1)\in\mathcal{O}$
is the same). Let $i\in\{1,\ldots,n\}$, then it suffices to prove that
$r_i(L(1))=L(1)$.
By the explicit expression of $T_i\in
\mathcal{H}\subset \mathcal{D}_q$ as a $q$-difference reflection
operator, we have $T_i^{\ddagger\tau}(1)=t_i^{-1}1$, 
where (recall) $t_i=t$ for $i=1,\ldots,n-1$.
It follows that
\[L(T_i^{\ddagger\tau}1)=t_i^{-1}L(1).
\]
On the other hand, by the transformation behaviour 
\eqref{transformationauxnew1}
of $L$ under the action of the double affine Hecke algeba, we have
\begin{equation*}
\begin{split}
L(T_i^{\ddagger\tau}1)&=
\bigl(\kappa_\tau\circ\ddagger_\tau^{-1}\bigr)(T_i^{\ddagger\tau})
L(1)\\
&=T_i^{\sigma\tau}{}^{-1}L(1)\\
&=t_i^{-1}L(1)+t_i^{-1}c_{a_i}^{\sigma\tau}
\bigl(r_i L(1)-L(1)\bigr).
\end{split}
\end{equation*}
Combining the two outcomes we obtain
$r_i(L(1))=L(1)$ in $\mathcal{O}$.

The next step is to show that $L(1), 
L_\ddagger(1)\in \mathcal{O}$ are in fact
$\mathcal{W}$-invariant, where $\mathcal{W}$
acts as constant coefficient
$q$-difference reflection operators
on $\mathcal{O}$ by \eqref{Wfund}. 
It suffices to show that $r_0(L(1))=L(1)$
and $r_0(L_\ddagger(1))=
L_{\ddagger}(1)$, where 
$r_0\in\mathcal{W}$
is the simple reflection associated to the simple root $a_0$.
For $L_{\ddagger}(1)$ it follows
from the identity $\bigl(\kappa^I_{\ddagger\tau}\circ\ddagger_\tau
\bigr)(T_0^\tau)=T_0^{\sigma\tau}$ 
that $r_0(L_\ddagger(1))=L_{\ddagger}(1)$ 
by repeating the
argument of the previous paragraph. For $L(1)$, 
we observe that
\[L\bigl(Y_1^{\ddagger\tau}1\bigr)=
L\bigl(Y_1^{\ddagger\tau}
E_{\ddagger\tau}(s_0^{\ddagger\tau};\cdot)\bigr)=
s_{0,1}^{\ddagger\tau}L(1)\]
on the one hand, where $s_{0,1}^{\ddagger\tau}$ is the first
coordinate of the spectral point $s_0^{\ddagger\tau}\in
\mathcal{S}_{\ddagger\tau}$,
\[ s_{0,1}^{\ddagger\tau}=u_0^{-1}t_n^{-1}t^{2(1-n)}.
\]
On the other hand,
\begin{equation*}
\begin{split}
L(Y_1^{\ddagger\tau}1)&=
\bigl(\kappa_\tau\circ\ddagger_\tau^{-1}\bigr)(Y_1^{\ddagger\tau})
L(1)\\
&=Y_1^{\sigma\tau}{}^{-1}L(1)\\
&=T_0^{\sigma\tau}{}^{-1}U_{\sigma\tau}^{-1}L(1)\\
&=t_n^{-1}t^{2(1-n)}T_0^{\sigma\tau}{}^{-1}L(1),
\end{split}
\end{equation*}
since we have seen that 
\[U_{\sigma\tau}=T_1^{\sigma\tau}T_2^{\sigma\tau}\cdots
T_{n-1}^{\sigma\tau}T_n^{\sigma\tau}T_{n-1}^{\sigma\tau}\cdots
T_2^{\sigma\tau}T_1^{\sigma\tau}
\]
acts as multiplication by $t_nt^{2(n-1)}$ on $L(1)$.
Hence $T_0^{\sigma\tau}{}^{-1}L(1)=u_0^{-1}L(1)$,
and a similar argument as before yields $r_0(L(1))=L(1)$.

The completion of the proof of the lemma is now standard:
by the change of variables
\[x=e^{2\pi iw}:=(e^{2\pi iw_1}, e^{2\pi iw_2},\cdots,e^{2\pi iw_n}),
\]
we obtain analytic functions
\[ w\mapsto L(1)(e^{2\pi iw}),\qquad w\mapsto
L_{\ddagger}(1)(e^{2\pi iw})
\]
on $\mathbb{C}^n/\Gamma_q$,
where $\Gamma_q=\mathbb{Z}^n+\mathbb{Z}^nv$ and $v$ is an element
in the upper half plane satisfying $q=e^{2\pi iv}$.
By Liouville's Theorem, these functions must be constants.
\end{proof}
We fix a convenient normalization for auxiliary transforms
as follows.
\begin{Def} 
We call an auxiliary transform $L$ \textup{(}respectively
$L_\ddagger$\textup{)} 
associated to $\alpha$ \textup{(}respectively 
$\alpha_\ddagger$\textup{)}
normalized when $L(1)=G_{\tau\sigma\tau}(s_0^\tau)1$ 
\textup{(}respectively 
$L_{\ddagger}(1)=G_{\tau\sigma\tau}(s_0^{\tau})1$\textup{)}.
\end{Def} 
\begin{rem}\label{Gaussianevaluation}
The normalization constant $G_{\tau\sigma\tau}(s_0^\tau)$
can be written as
\[
G_{\tau\sigma\tau}(s_0^\tau)=
\prod_{i=1}^n\bigl(bct^{2(n-i)}, dt^{2(i-n)}/a;q\bigr)_{\infty}^{-1}
\]
in terms of the Askey-Wilson parametrization \eqref{AWpar}
for part of the multiplicity function $\alpha$. 
\end{rem}
The fact that the above choice of normalization is the natural one
becomes clear in the next subsection. 
{}From the previous lemma it follows
that a normalized auxiliary transform associated to
$\alpha$ (respectively $\alpha_\ddagger$) is unique, provided
that it exists.

In the following statement we use the stablility condition
\begin{equation}\label{stabilitys}
s_{\lambda}^{\tau}=
s_{\lambda}^{\sigma\tau},\qquad\forall\, \lambda\in\Lambda
\end{equation}
of the polynomial spectrum, which follows from the fact that
the spectral points $s_\lambda\in\mathcal{S}=\mathcal{S}_\alpha$
only depend on the values of the multiplicity function 
$\mathbf{t}$ at the inmultipliable roots $R$ of $R_{nr}$.

\begin{cor}\label{auxcor}
Suppose that $L,L_\ddagger: \mathcal{A}\rightarrow \mathcal{O}$
are normalized auxiliary transforms associated to
$\alpha$ and $\alpha_\ddagger$, respectively. Then
for $s\in\mathcal{S}_{\tau}$ we have
\begin{equation*}
\begin{split}
L\bigl(E_{\ddagger\tau}(s^{-1};\cdot)\bigr)&=
d^\tau(s)E_{\sigma\tau}(s;\cdot),\\
L_\ddagger\bigl(E_\tau(s;\cdot)\bigr)&=
d^{\ddagger\tau}(s^{-1})E_{\sigma\tau}(s;\cdot)
\end{split}
\end{equation*}
for unique functions $d^\tau\in \mathcal{F}(\mathcal{S}_\tau)$
and $d^{\ddagger\tau}\in \mathcal{F}(\mathcal{S}_{\ddagger\tau})$
satisfying $d^\tau(s_0^\tau)=d^{\ddagger\tau}(s_0^{\ddagger\tau})=
G_{\tau\sigma\tau}(s_0^{\tau})$.
\end{cor}
\begin{proof}
For $i=1,\ldots,n$ we have
\[ \bigl(\kappa_\tau\circ\ddagger_\tau^{-1}\bigr)(Y_i^{\ddagger\tau})=
Y_i^{\sigma\tau}{}^{-1},
\]
hence for $s\in\mathcal{S}_\tau$,
\[ Y_i^{\sigma\tau}{}^{-1}
L\bigl(E_{\ddagger\tau}(s^{-1};\cdot)\bigr)=
L\bigl(Y_i^{\ddagger\tau}E_{\ddagger\tau}(s^{-1};\cdot)\bigr)=
s_i^{-1}L\bigl(E_{\ddagger\tau}(s^{-1};\cdot)\bigr).
\]
Furthermore, $L\bigl(E_{\ddagger\tau}(s^{-1};\cdot)\bigr)
\in\mathcal{A}$, hence by the polynomial Mac\-do\-nald-Koorn\-win\-der
theory (see Subsection 4.2) and by \eqref{stabilitys},
\[ L\bigl(E_{\ddagger\tau}(s^{-1};\cdot)\bigr)=
d^\tau(s)E_{\sigma\tau}(s;\cdot)
\]
for some constant $d^\tau(s)\in\mathbb{C}$. The normalization
$d^\tau(s_0^\tau)=G_{\tau\sigma\tau}(s_0^\tau)$ follows from the
normalization of $L$ and the fact that 
$E_{\ddagger\tau}(s_0^{\ddagger\tau};\cdot)=
E_{\sigma\tau}(s_0^{\sigma\tau};\cdot)=1$ is the Laurent
polynomial identically equal to one. The formulas for
$L_\ddagger$ are proven similarly, now using the fact that
\[\bigl(\kappa_{\ddagger\tau}^I\circ\ddagger_\tau\bigr)(Y_i^\tau)=
Y_i^{\sigma\tau}
\]
for $i=1,\ldots,n$.
\end{proof}

\begin{Def}
We call the functions $d^\tau\in \mathcal{F}(\mathcal{S}_\tau)$
and $d^{\ddagger\tau}\in \mathcal{F}(\mathcal{S}_{\ddagger\tau})$ 
in Corollary \ref{auxcor} the generalized eigenvalues associated 
to the normalized auxiliary transform $L$ and $L_\ddagger$,
respectively.
\end{Def}

We prove in the next subsection the {\it existence} of the normalized 
auxiliary transforms associated to $\alpha$ and $\alpha_\ddagger$
by determining the related generalized eigenvalues
$d^\tau$ and $d^{\ddagger\tau}$ explicitly.


\subsection{The generalized eigenvalues}
In this subsection we give full details on determining
the generalized eigenvalue
$d^\tau\in \mathcal{F}(\mathcal{S}_\tau)$. It leads to
the existence of the normalized auxiliary transform
associated to $\alpha$.
At the end of the subsection we indicate how a similar procedure leads
to the explicit expression for the generalized eigenvalue
$d^{\ddagger\tau}\in \mathcal{F}(\mathcal{S}_{\ddagger\tau})$.

Recall that the space $\mathcal{F}(\mathcal{S}_\tau)$ consisting
of functions 
$\phi: \mathcal{S}_\tau\rightarrow
\mathbb{C}$ is an $\mathcal{H}_{\ddagger\tau\sigma}$-module
in view of Lemma \ref{actioncompatibleH}. 
Since $\mathcal{S}_\tau=\mathcal{S}_{\sigma\tau}$, the space
$\mathcal{F}(\mathcal{S}_\tau)$ also has the structure of an 
$\mathcal{H}_{\ddagger\sigma\tau\sigma}$-module. We consider
$\mathcal{F}(\mathcal{S}_\tau)$ as a commutative algebra
by point-wise multiplication, and we view 
$\mathcal{F}(\mathcal{S}_\tau)$ as sub-algebra of
$\hbox{End}_{\mathbb{C}}(\mathcal{F}(\mathcal{S}_\tau))$ by
identifying $g\in \mathcal{F}(\mathcal{S}_\tau)$ with the
endomorphism of $\mathcal{F}(\mathcal{S}_\tau)$
defined as multiplication by
$g$.
\begin{lem}\label{eigenconj}
Let $d^\tau\in \mathcal{F}(\mathcal{S}_\tau)$ and let
$L: \mathcal{A}\rightarrow \mathcal{A}$ be the linear map
defined by
\[L(E_{\ddagger\tau}(s^{-1};\cdot)\bigr)=
d^\tau(s)E_{\sigma\tau}(s;\cdot),
\qquad \forall\,s\in\mathcal{S}_{\tau}.
\]
If
\begin{equation}\label{conj}
X\circ d^\tau=d^\tau\circ \nu_{\ddagger\tau\sigma}(X),
\qquad \forall\, X\in \mathcal{H}_{\ddagger\tau\sigma}
\end{equation}
within the endomorphism ring
$\hbox{End}_{\mathbb{C}}(\mathcal{F}(\mathcal{S}_\tau))$,
where $\nu_{\ddagger\tau\sigma}:
\mathcal{H}_{\ddagger\tau\sigma}\rightarrow
\mathcal{H}_{\ddagger\sigma\tau\sigma}$ is the unital algebra isomorphism
defined by
\[
\nu_{\ddagger\tau\sigma}=
\dagger_{\sigma\tau\sigma}\circ\psi_{\sigma\tau}\circ
\kappa_\tau\circ\ddagger_\tau^{-1}\circ\psi_{\ddagger\tau}^{-1},
\]
then $L$ satisfies the transformation behaviour \eqref{aux1}
under the action of the double affine Hecke algebra 
$\mathcal{H}_{\ddagger\tau}$.
\end{lem}
\begin{proof}
The space $\mathcal{F}(\mathcal{S}_\tau,\mathcal{A})$
consisting of functions $\phi: \mathcal{S}_\tau\rightarrow
\mathcal{A}$ is an
$\mathcal{H}_{\ddagger\tau\sigma}$-module by
\[ \bigr((X\cdot\phi)(\gamma)\bigr)(x)=\bigl(X\cdot
\phi_x\bigr)(\gamma),\qquad
x\in (\mathbb{C}^\times)^n,\,\,\,\gamma\in\mathcal{S}_\tau
\]
with $\phi_x\in \mathcal{F}(\mathcal{S}_\tau)$ defined
by $\phi_x(\gamma)=\bigl(\phi(\gamma)\bigr)(x)$. Since 
$\mathcal{S}_\tau=\mathcal{S}_{\sigma\tau}$, the same formula
defines an $\mathcal{H}_{\ddagger\sigma\tau\sigma}$-module structure
on $\mathcal{F}(\mathcal{S}_\tau,\mathcal{A})$. Furthermore, point-wise
multiplication defines a commutative algebra structure on
$\mathcal{F}(\mathcal{S}_\tau,\mathcal{A})$.

The canonical embedding
$\mathcal{F}(\mathcal{S}_\tau) \hookrightarrow
\mathcal{F}(\mathcal{S}_\tau,\mathcal{A})$ of algebras via the 
map $\mathbb{C}\rightarrow \mathcal{A}$, 
$\lambda\mapsto \lambda 1$ 
is compatible with the above actions of the double affine
Hecke algebras $\mathcal{H}_{\ddagger\tau\sigma}$ and 
$\mathcal{H}_{\ddagger\sigma\tau\sigma}$.
For $\phi\in\mathcal{F}(\mathcal{S}_\tau,\mathcal{A})$
and $s\in\mathcal{S}_\tau$ we write $(X\cdot\phi)(s)$ for the action of
$X$ on $\phi$ as described above, while we write $X(\phi(s))$ for the 
action of $X$ on the element $\phi(s)\in\mathcal{A}$.

Fix $d^\tau\in \mathcal{F}(\mathcal{S}_\tau)$ and define a linear map
$L:\mathcal{A}\rightarrow\mathcal{A}$ by
\[L\bigl(E_{\ddagger\tau}(s^{-1};\cdot)\bigr)=
d^\tau(s)E_{\sigma\tau}(s;\cdot),\qquad \forall\, s\in
\mathcal{S}_\tau.
\]
Assume furthermore that \eqref{conj} holds true for some
(yet to be determined) algebra isomorphism $\nu_{\ddagger\tau\sigma}:
\mathcal{H}_{\ddagger\tau\sigma}\rightarrow 
\mathcal{H}_{\ddagger\sigma\tau\sigma}$. 
We define $\phi_1,\phi_2\in \mathcal{F}(\mathcal{S}_\tau,\mathcal{A})$ by
\begin{equation*}
\begin{split}
\phi_1(s)&=L\bigl(E_{\ddagger\tau}(s^{-1};\cdot)\bigr)=
L\bigl(\mathfrak{E}_{\mathcal{A},\ddagger\tau}(s,\cdot)\bigr),\\
\phi_2(s)&=E_{\sigma\tau}(s;\cdot)=
\mathfrak{E}_{\mathcal{A},\sigma\tau}(s^{-1},\cdot)
\end{split}
\end{equation*}
for $s\in\mathcal{S}_\tau$, so that $\phi_1=d^\tau\phi_2$
in the algebra $\mathcal{F}(\mathcal{S}_\tau,\mathcal{A})$. 
Then we compute for $s\in\mathcal{S}_\tau$ and
$X\in\mathcal{H}_{\ddagger\tau\sigma}$ by Proposition
\ref{allrelations},
\[
(X\cdot\phi_1)(s)
=L\bigl(\psi_{\ddagger\tau}^{-1}(X)
\mathfrak{E}_{\mathcal{A},\ddagger\tau}(s,\cdot)\bigr),
\]
with the action of $\psi_{\ddagger\tau}^{-1}(X)$ on $\mathcal{A}$ 
in the right hand side of the equality.
On the other hand, under the above assumptions,
\begin{equation*}
\begin{split}
(X\cdot\phi_1)(s)&=\bigl(X\cdot(d^\tau\phi_2)\bigr)(s)\\
&=\Bigl(d^\tau\bigl(\nu_{\ddagger\tau\sigma}(X)\cdot
\phi_2\bigr)\Bigr)(s)\\
&=d^\tau(s)
\bigl(\psi_{\sigma\tau}^{-1}\circ\dagger_{\ddagger\sigma\tau\sigma}
\circ\nu_{\ddagger\tau\sigma}\bigr)(X)\bigl(\phi_2(s)\bigr)\\
&=\bigl(\psi_{\sigma\tau}^{-1}\circ\dagger_{\ddagger\sigma\tau\sigma}
\circ\nu_{\ddagger\tau\sigma}\bigr)(X)\bigl(\phi_1(s)\bigr).
\end{split}
\end{equation*}
Equating both outcomes yields
\[L\bigl(X\mathfrak{E}_{\mathcal{A},\ddagger\tau}(s,\cdot)\bigr)
=\bigl(\psi_{\sigma\tau}^{-1}\circ\dagger_{\ddagger\sigma\tau\sigma}
\circ\nu_{\ddagger\tau\sigma}\circ\psi_{\dagger\tau}\bigr)(X)
\bigl(\phi_1(s)\bigr)\]
for $s\in\mathcal{S}_\tau$ and
$X\in\mathcal{H}_{\ddagger\tau}$.
On the other hand, $L$ satisfies the intertwining property
\eqref{aux1} when  
\[L\bigl(X\mathfrak{E}_{\mathcal{A},\ddagger\tau}(s,\cdot)\bigr)=
\bigl(\kappa_\tau\circ\ddagger_\tau^{-1}\bigr)(X)\bigl(\phi_1(s)\bigr)
\]
for all $s\in\mathcal{S}_{\tau}$ and 
$X\in\mathcal{H}_{\ddagger\tau}$.
This will thus be the case when
\[\psi_{\sigma\tau}^{-1}\circ\dagger_{\ddagger\sigma\tau\sigma}
\circ\nu_{\ddagger\tau\sigma}\circ\psi_{\dagger\tau}=
\kappa_\tau\circ\ddagger_\tau^{-1},\]
i.e. when the isomorphism $\nu_{\ddagger\tau\sigma}$ is given by
 \[
\nu_{\ddagger\tau\sigma}=
\dagger_{\sigma\tau\sigma}\circ\psi_{\sigma\tau}\circ
\kappa_\tau\circ\ddagger_\tau^{-1}\circ\psi_{\ddagger\tau}^{-1}.
\]
\end{proof}

We thus search for a function 
$d^\tau\in \mathcal{F}(\mathcal{S}_\tau)$, normalized by
$d^\tau(s_0^\tau)=G_{\tau\sigma\tau}(s_0^\tau)$, which induces the
isomorphism $\nu_{\ddagger\tau\sigma}:
\mathcal{H}_{\ddagger\tau\sigma}\rightarrow
\mathcal{H}_{\ddagger\sigma\tau\sigma}$ via the formula \eqref{conj}. 
\begin{prop}\label{Lexists}
{\bf i)} We have $\nu_{\ddagger\tau\sigma}=\tau_{\ddagger\tau\sigma}$.\\
{\bf ii)} The normalized auxiliary transform
$L:\mathcal{A}\rightarrow \mathcal{A}$ associated
to $\alpha$ exists. Its generalized eigenvalue $d^\tau$ is given by
\[ d^\tau(s)=G_{\tau\sigma\tau}(s),\qquad \forall\,s\in\mathcal{S}_\tau.
\]
In other words,
\[ L\bigl(E_{\ddagger\tau}(s^{-1};\cdot)\bigr)=
G_{\tau\sigma\tau}(s)E_{\sigma\tau}(s;\cdot),\qquad
\forall\,s\in\mathcal{S}_\tau.
\]
\end{prop}
\begin{proof}
{\bf i)} First observe that $\tau_{\ddagger\tau\sigma}$ defines an 
algebra isomorphism
$\tau_{\ddagger\tau\sigma}:\mathcal{H}_{\ddagger\tau\sigma}\rightarrow
\mathcal{H}_{\ddagger\sigma\tau\sigma}$ since
$\sigma\tau\sigma=\tau\sigma\tau$ when acting upon 
difference multiplicity functions. Hence we only have to show that
the action of $\nu_{\ddagger\tau\sigma}$ and $\tau_{\ddagger\tau\sigma}$
coincides on the algebraic generators $T_i^{\ddagger\tau\sigma}$
($i=0,\ldots,n$) and $x_1$ of $\mathcal{H}_{\ddagger\tau\sigma}$.
We omit parameter dependencies in the following computations.
It is easy to check that
\[\nu(T_i)=T_i=\tau(T_i),\qquad i=1,\ldots,n.
\]
Furthermore, we have
\begin{equation*}
\begin{split}
\nu(x_1)&=\dagger\psi\kappa\ddagger^{-1}\psi^{-1}(x_1)=
\dagger\psi\kappa\ddagger^{-1}(Y_1^{-1})\\
&=\dagger\psi\kappa(Y_1)
=\dagger\psi(Y_1)
=\dagger(x_1^{-1})=x_1=\tau(x_1).
\end{split}
\end{equation*}
Since $Y_1=UT_0$ with $U\in\mathcal{H}$ given by
\eqref{U}, we furthermore obtain
\begin{equation*}
\begin{split}
\nu(T_0)&=\dagger\psi\kappa\ddagger^{-1}\psi^{-1}(T_0)\\
&=\dagger\psi\kappa\ddagger^{-1}\psi^{-1}(U^{-1}Y_1)\\
&=\dagger\psi\kappa\ddagger^{-1}(x_1^{-1}U^{-1})\\
&=\dagger\psi\kappa(Ux_1)\\
&=\dagger\psi(q^{\frac{1}{2}}x_1^{-1}T_0)\\
&=\dagger\psi(q^{\frac{1}{2}}x_1^{-1}U^{-1}Y_1)\\
&=\dagger(q^{\frac{1}{2}}x_1^{-1}U^{-1}Y_1)\\
&=q^{\frac{1}{2}}x_1T_0^{-1}=\tau(T_0)
\end{split}
\end{equation*}
where we use in the last equality that the sixth coordinate of the
underlying difference multiplicity function is $q^{-\frac{1}{2}}$, 
so that
$\tau_{\ddagger\tau\sigma}(T_0^{\ddagger\tau\sigma})=
q^{\frac{1}{2}}x_1T_0^{\ddagger\sigma\tau\sigma}{}^{-1}$.
We conclude that $\nu_{\ddagger\tau\sigma}=
\tau_{\ddagger\tau\sigma}$.\\ 
{\bf ii)} Using the realization of the isomorphism
$\tau$ as conjugation by the Gaussian and using Lemma
\ref{actioncompatibleH}, we observe that
the function
\[
d^\tau(s)=G_{\ddagger\tau\sigma}(s)^{-1}=G_{\tau\sigma\tau}(s),\qquad
\forall\,s\in\mathcal{S}_\tau
\]
satisfies \eqref{conj}. Furthermore, for this function $d^\tau$ we
have the desired normalization 
$d^\tau(s_0^\tau)=G_{\tau\sigma\tau}(s_0^\tau)$.
The result now follows from Lemma \ref{eigenconj}.
\end{proof}

We leave it as an exercise to the reader to repeat the arguments of
this subsection to determine the normalized auxiliary transform 
$L_{\ddagger}$ associated to $\alpha_\ddagger$. It leads to
the transformation behaviour 
\begin{equation*}
\begin{split}
X\circ d^{\ddagger\tau}&=d^{\ddagger\tau}\circ
\bigl(\psi_{\sigma\tau}\circ\kappa_{\ddagger\tau}^I\circ\ddagger_\tau\circ
\psi_\tau^{-1}\bigr)(X)\\
&=d^{\ddagger\tau}\circ\tau_{\tau\sigma\tau}^{-1}(X),
\qquad\qquad\qquad\qquad\qquad\forall\, X\in\mathcal{H}_{\tau\sigma}
\end{split}
\end{equation*}
in $\hbox{End}_{\mathbb{C}}
\bigl(\mathcal{F}(\mathcal{S}_{\ddagger\tau})\bigr)$
for the associated generalized eigenvalue
$d^{\ddagger\tau}\in \mathcal{F}(\mathcal{S}_{\ddagger\tau})$,
and hence to the following result.
\begin{prop}\label{Ldaggerexists} 
The normalized auxiliary transform 
$L_{\ddagger}: \mathcal{A}\rightarrow
\mathcal{A}$ exists. Its generalized eigenvalue $d^{\ddagger\tau}\in 
\mathcal{F}(\mathcal{S}_{\ddagger\tau})$ is given by
\[d^{\ddagger\tau}(s)=G_{\tau\sigma\tau}(s)
=G_{\tau\sigma\tau}(s^{-1}),\qquad \forall\,
s\in\mathcal{S}_{\ddagger\tau}.
\]
In other words,
\[L_\ddagger(E_{\tau}(s;\cdot))=G_{\tau\sigma\tau}(s)
E_{\sigma\tau}(s;\cdot),\qquad \forall\, s\in\mathcal{S}_\tau.
\]
\end{prop}


\subsection{The Cherednik kernels}
In this subsection we construct meromorphic kernels $\mathfrak{E},
\mathfrak{E}_\ddagger$ satisfying the transformation behaviour
\eqref{transformation1} under the action of the double affine Hecke
algebra. For this
it suffices to construct 
kernels $\mathfrak{F},\mathfrak{F}_\ddagger\in
\mathcal{O}\bigl((\mathbb{C}^\times)^n\times 
(\mathbb{C}^\times)^n\bigr)$ such that the normalized
auxiliary transforms $L$ and $L_\ddagger$ associated to
$\alpha$ and $\alpha_\ddagger$ are given by
\begin{equation}\label{integral}
\bigl(Lp\bigr)(\gamma)=\bigl(\mathfrak{F}(\gamma,\cdot),
p\bigr)_{\mathcal{A},\tau},\qquad
\bigl(L_{\ddagger}p\bigr)(\gamma)=
\bigl(p,\mathfrak{F}_\ddagger(\gamma,\cdot)\bigr)_{\mathcal{A},\tau},
\qquad \forall\,p\in\mathcal{A}
\end{equation}
respectively, cf. Lemma \ref{EFconnection} and 
Lemma \ref{Auxiliarylink}. Formal series expansions in
Mac\-do\-nald-Koorn\-win\-der polynomials 
can now immediately be written down for such
kernels $\mathfrak{F}$ and $\mathfrak{F}_\ddagger$
in view of Proposition \ref{Lexists}, Proposition \ref{Ldaggerexists} 
and the orthogonality relations \eqref{ortho} 
for the Macdonald-Koornwinder polynomials.
To ensure that these formal series expansions define 
analytic kernels, we need to determine bounds 
for the Mac\-do\-nald-Koorn\-win\-der
polynomials $E(s_\lambda;\cdot)$ in their degree $\lambda\in\Lambda$.
The following proposition suffices for our purposes.

\begin{prop}\label{MKbounds}
For any compacta $K\subset (\mathbb{C}^\times)^n$, there exists a
constant $C=C_K>0$ such that
\[ |E(s_\lambda;x)|\leq C^{N(\lambda)},\qquad \forall\,x\in K,\quad
\forall\,\lambda\in\Lambda,
\]
where $N(\lambda)=\sum_{i=1}^n|\lambda_i|$
for $\lambda=\sum_{i=1}^n\lambda_i\epsilon_i\in\Lambda$.
\end{prop}

The proof of the proposition, which is a bit technical, is based on
explicit recurrence relations for the Macdonald-Koornwinder
polynomials. These recurrence relations are a direct consequence
of Proposition \ref{allrelations}. For details of the proof, we
refer the reader to the appendix.

\begin{rem}\label{symmbound}
For the symmetric Macdonald-Koornwinder polynomials,
we have for any
compact set $K\subset (\mathbb{C}^\times)^n$,
\[|E^+(s_\lambda;x)|\leq C^{N(\lambda)},\qquad
\forall\, x\in K,\quad \forall\, \lambda\in\Lambda^+,
\]
for some constant $C=C_K>0$. This follows for instance 
from Proposition \ref{MKbounds}, 
Lemma \ref{boundedcoeff}, the fact that $N(w\cdot\mu)=N(\mu)$
for $w\in W_0$ and $\mu\in\Lambda$ and from 
the explicit expansion 
\[E^+(s_\lambda;x)=\sum_{\mu\in W_0\cdot\lambda}
\frac{\mathcal{C}_\sigma(s_\mu^\ddagger)}
{\mathcal{C}_\sigma(s_0^\ddagger)}\,E(s_\mu;x),\qquad \lambda\in\Lambda^+,
\]
of the symmetric Macdonald-Koornwinder
polynomial as linear combination of the nonsymmetric ones, 
see \cite[Thm. 6.6]{St1} or \cite[Thm. 3.27]{St2}.
\end{rem}
\begin{cor}\label{convcor}
{\bf a)}
The series expansion
\begin{equation}\label{Fseries}
\mathfrak{F}(\gamma,x):=\sum_{s\in\mathcal{S}_\tau}
\frac{G_{\tau\sigma\tau}(s)}{\bigl(E_\tau(s;\cdot),
E_{\ddagger\tau}(s^{-1};\cdot)\bigr)_{\mathcal{A},\tau}}
E_\tau(s;x)E_{\sigma\tau}(s;\gamma)
\end{equation}
converges absolutely and uniformly on compacta
of $(\mathbb{C}^\times)^n\times (\mathbb{C}^\times)^n$. 
The series expansion
$\mathfrak{F}$ defined by \eqref{Fseries} is the unique
analytic kernel on $(\mathbb{C}^\times)^n\times (\mathbb{C}^\times)^n$
such that $\bigl(Lp\bigr)(\gamma)=\bigl(\mathfrak{F}(\gamma,\cdot),p
\bigr)_{\mathcal{A},\tau}$ for all $p\in\mathcal{A}$,
where $L$ is the normalized auxiliary transform associated to $\alpha$.\\
{\bf b)} The series expansion
\begin{equation}\label{Fdaggerseries}
\mathfrak{F}_\ddagger(\gamma,x):=\sum_{s\in\mathcal{S}_\tau}
\frac{G_{\tau\sigma\tau}(s)}{\bigl(E_\tau(s;\cdot),
E_{\ddagger\tau}(s^{-1};\cdot)\bigr)_{\mathcal{A},\tau}}
E_{\ddagger\tau}(s^{-1};x)E_{\sigma\tau}(s;\gamma)
\end{equation} 
converges absolutely and uniformly on compacta
of $(\mathbb{C}^\times)^n\times (\mathbb{C}^\times)^n$.
The series expansion
$\mathfrak{F}_\ddagger$ defined by \eqref{Fdaggerseries} is the
unique analytic kernel on $(\mathbb{C}^\times)^n\times (\mathbb{C}^\times)^n$
such that $\bigl(L_{\ddagger}p\bigr)(\gamma)=\bigl(p,
\mathfrak{F}_\ddagger(\gamma,\cdot)\bigr)_{\mathcal{A},\tau}$ 
for all $p\in\mathcal{A}$, where $L_\ddagger$ 
is the normalized auxiliary transform associated to $\alpha_\ddagger$.
\end{cor}
\begin{proof}
We first prove that the series expansions \eqref{Fseries} 
and \eqref{Fdaggerseries}
converge absolutely and uniformly on compacta of
$(\mathbb{C}^\times)^n\times (\mathbb{C}^\times)^n$.
For this we need to consider the behaviour of the coefficients
\begin{equation}\label{coefficient}
\frac{G_{\tau\sigma\tau}(s)}
{\bigl(E_{\tau}(s;\cdot),E_{\ddagger\tau}(s^{-1};\cdot)
\bigr)_{\mathcal{A},\tau}}=\frac{1}{N_{\tau\sigma}(s_0^{\ddagger\tau})
\bigl(1,1\bigr)_{\mathcal{A},\tau}}G_{\tau\sigma\tau}(s)
N_{\tau\sigma}(s^{-1})
\end{equation}
as function of $s\in\mathcal{S}_\tau$ 
in the expansion sums \eqref{Fseries} and \eqref{Fdaggerseries}, 
see Theorem \ref{corpoly} for the second equality in 
\eqref{coefficient}.
By the $W_0$-invariance of the Gaussian and by Lemma
\ref{actioncompatibleW}, we have for $\lambda\in\Lambda$,
\begin{equation}\label{WnulGaussian}
 G_{\tau\sigma\tau}(s_{\lambda}^\tau)=
G_{\tau\sigma\tau}(s_{\lambda^+}^\tau),
\end{equation}
where $\lambda^+$ is the unique element in $(W_0\cdot\lambda)\cap
\Lambda^+$. Using the explicit expression of the Gaussian $G$ and using
that
\[s_\mu^\tau=(t_nu_0t^{2(n-1)}q^{\mu_1},
t_nu_0t^{2(n-2)}q^{\mu_2},\ldots,
t_nu_0q^{\mu_n}),
\qquad \mu=\sum_{i=1}^n\mu_i\epsilon_i\in\Lambda^+
\]
it is now easy to show that
\begin{equation}\label{Gaussform}
\frac{G_{\tau\sigma\tau}(s_\mu^\tau)}
{G_{\tau\sigma\tau}(s_0^\tau)}=
\prod_{i=1}^n\frac{\bigl(bct^{2(n-i)};q\bigr)_{\mu_i}}
{\bigl(qat^{2(n-i)}/d;q\bigr)_{\mu_i}}
\left(\frac{-q^{\frac{1}{2}}at^{2(n-i)}}{d}\right)^{\mu_i}
q^{\mu_i^2/2}
\end{equation}
for $\mu\in\Lambda^+$ (cf. Remark \ref{Gaussianevaluation}), 
where we used the Askey-Wilson parametrization
\eqref{AWpar} for part of the multiplicity function $\alpha$.
We thus obtain the bounds 
\begin{equation}\label{Gaussianbound}
 |G_{\tau\sigma\tau}(s_{\lambda}^\tau)|\leq 
c_1c_2^{N(\lambda)}q^{\langle \lambda,\lambda\rangle/2},\qquad
\forall\,\lambda\in\Lambda
\end{equation}
for certain $\lambda$-independent constants $c_1,c_2>0$.

We have seen that the factor 
$N_{\tau\sigma}(s_\lambda^{\ddagger\tau})$ can be evaluated explicitly,
see \eqref{Ndecomposition}
and \eqref{Nplus}. {}From these explicit expressions and
from Lemma \ref{boundedcoeff} it follows that
\[
|N_\tau(s_\lambda^{\ddagger\tau})|\leq d_1d_2^{N(\lambda)},\qquad
\forall\,\lambda\in\Lambda\]
for certain $\lambda$-independent coefficients $d_1,d_2>0$.
Combined with Proposition \ref{MKbounds} and the fact that $0<q<1$,
we see that the ``Gaussian term'' $q^{\langle \lambda,\lambda\rangle/2}$
in the bounds for $G_{\tau\sigma\tau}$ ensure the absolute and
uniform convergence
of the sums \eqref{Fseries} and \eqref{Fdaggerseries} on compact
subsets of $(\gamma,x)$ in $(\mathbb{C}^\times)^n\times 
(\mathbb{C}^\times)^n$. In particular, the series expansions
\eqref{Fseries} and \eqref{Fdaggerseries} define analytic kernels
$\mathfrak{F},\mathfrak{F}_\ddagger\in\mathcal{O}\bigl((\mathbb{C}^\times)^n
\times (\mathbb{C}^\times)^n\bigr)$.

By the orthogonality relations \eqref{ortho} for the
Macdonald-Koornwinder polynomials and the explicit series expansions
\eqref{Fseries} and \eqref{Fdaggerseries} for $\mathfrak{F}$ and 
$\mathfrak{F}_\ddagger$, it is immediate that
\[\bigl(\mathfrak{F}(\gamma,\cdot),E_{\ddagger\tau}(s^{-1};\cdot)
\bigr)_{\mathcal{A},\tau}=G_{\tau\sigma\tau}(s)E_{\sigma\tau}(s;\gamma)=
\bigl(E_\tau(s;\cdot), \mathfrak{F}_\ddagger(\gamma,\cdot)
\bigr)_{\mathcal{A},\tau},\qquad \forall\,s\in\mathcal{S}_\tau.
\]
In view of Proposition \ref{Lexists} and Proposition \ref{Ldaggerexists},
we conclude that the normalized auxiliary transforms $L$ and $L_\ddagger$
associated to $\alpha$ and $\alpha_\ddagger$ are given by
the integral transforms \eqref{integral}, with 
$\mathfrak{F}$ and $\mathfrak{F}_\ddagger$ the analytic kernels
defined by the
series expansions \eqref{Fseries} and \eqref{Fdaggerseries}
respectively. Clearly, the linear mappings \eqref{integral}
determine the analytic kernels $\mathfrak{F}$ and
$\mathfrak{F}_\ddagger$ uniquely (cf. Lemma \ref{Auxiliarylink}).
\end{proof}

\begin{Def}
{\bf a)}
The kernel $\mathfrak{E}=\mathfrak{E}_\alpha\in\mathcal{M}
\bigl((\mathbb{C}^\times)^n\times (\mathbb{C}^\times)^n\bigr)$
defined by
\[\mathfrak{E}(\gamma,x)=
G_{\sigma\tau}(\gamma)G_\tau(x)\mathfrak{F}(\gamma,x),
\]
with $\mathfrak{F}$ the explicit series
expansion \eqref{Fseries}, 
is called the Cherednik kernel associated to $\alpha$.\\
{\bf b)} The kernel
$\mathfrak{E}_\ddagger=\mathfrak{E}_{\alpha_\ddagger}\in\mathcal{M}
\bigl((\mathbb{C}^\times)^n\times (\mathbb{C}^\times)^n\bigr)$
defined by
\[
\mathfrak{E}_\ddagger(\gamma,x)=
G_{\sigma\tau}(\gamma)G_\tau(x)\mathfrak{F}_\ddagger(\gamma^{-1},x),
\]
with $\mathfrak{F}_\ddagger$ the explicit series
expansion \eqref{Fdaggerseries},
is called the Cherednik kernel associated to
$\alpha_\ddagger$.
\end{Def}
We can now state the following main result of this section. 
\begin{thm}\label{Cherednikkernels}
{\bf a)}
Up to a multiplicative constant, the 
Cherednik kernel $\mathfrak{E}$ associated to $\alpha$ is the unique 
non-zero meromorphic kernel such that
\begin{enumerate}
\item[{\bf i)}] The function
$G_{\sigma\tau}(\gamma)^{-1}G_\tau(x)^{-1}\mathfrak{E}(\gamma,x)$  
depends analytically on $(\gamma,x)\in (\mathbb{C}^\times)^n\times
(\mathbb{C}^\times)^n$.
\item[{\bf ii)}] For all $X\in \mathcal{H}=\mathcal{H}_\alpha$
we have $\bigl(X\mathfrak{E}(\gamma,\cdot)\bigr)(x)=
\bigl(\psi(X)\mathfrak{E}(\cdot,x)\bigr)(\gamma)$.
\end{enumerate}
{\bf b)} Up to a multiplicative constant, 
the Cherednik kernel $\mathfrak{E}_\ddagger$ associated to
$\alpha_\ddagger$ is the unique 
non-zero meromorphic kernel such that
\begin{enumerate}
\item[{\bf i)}] The function 
$G_{\sigma\tau}(\gamma)^{-1}G_\tau(x)^{-1}\mathfrak{E}_\ddagger(\gamma,x)$ 
depends analytically on 
$(\gamma,x)\in (\mathbb{C}^\times)^n\times (\mathbb{C}^\times)^n$.
\item[{\bf ii)}] For all $X\in \mathcal{H}_\ddagger$
we have $\bigl(X\mathfrak{E}_\ddagger(\gamma,\cdot)\bigr)(x)=
\bigl(\psi_\ddagger(X)\mathfrak{E}_\ddagger(\cdot,x)\bigr)(\gamma)$.
\end{enumerate}
\end{thm}
\begin{proof}
{\bf a)} It follows immediately from Corollary \ref{convcor},
Lemma \ref{Auxiliarylink} and Lemma \ref{EFconnection} 
that $\mathfrak{E}$ satisfies 
properties {\bf i)} and {\bf ii)}. Suppose that $\widetilde{\mathfrak{E}}$
is another non-zero meromorphic kernel satisfying {\bf i)} and {\bf
ii)}. Then 
\[\widetilde{\mathfrak{F}}(\gamma,x):=G_{\sigma\tau}(\gamma)^{-1}
G_{\tau}(x)^{-1}\widetilde{\mathfrak{E}}(\gamma,x)
\]
depends analytically on $(\gamma,x)\in (\mathbb{C}^\times)^n\times
(\mathbb{C}^\times)^n$. Furthermore, the linear map
$\widetilde{L}: \mathcal{A}\rightarrow \mathcal{O}$ 
defined by $\bigl(\widetilde{L}p\bigr)(\gamma)=
\bigl(\widetilde{\mathfrak{F}}(\gamma,\cdot),p\bigr)_{\mathcal{A},\tau}$
for $p\in\mathcal{A}$ is an auxiliary transform associated to $\alpha$
in view of Lemma \ref{EFconnection} and
Lemma \ref{Auxiliarylink}. 
Lemma \ref{L} then shows that  $\widetilde{L}=c\,L$
for some constant $c\in\mathbb{C}^\times$, where $L$ is the normalized
auxiliary transform associated to $\alpha$. By Corollary
\ref{convcor} we conclude that 
$\widetilde{\mathfrak{F}}=c\,\mathfrak{F}$,
where $\mathfrak{F}$ is given by the series expansion
\eqref{Fseries}, and hence
$\widetilde{\mathfrak{E}}=c\,\mathfrak{E}$.\\
{\bf b)} The proof is similar to the proof of {\bf a)}.
\end{proof}

\section{Properties of the Cherednik kernels}

In this section we analyze the Cherednik kernels further.
In Subsection 6.1 we prove an evaluation formula, which allows to normalize
the Cherednik kernels in such a manner that they 
meromorphic extend
the Macdonald-Koornwinder polynomials in their degrees (as will be
shown in Subsection 6.2).
In Subsection 6.1 we furthermore prove the duality of the 
normalized Cherednik kernels between
their geometric and spectral parameters.
In Subsection 6.3 we consider symmetric Cherednik kernels.
In case of reduced root systems, many of the results in this section
reduce to statements in Cherednik's paper \cite{C1}.

We keep the same generic conditions on the multiplicity function
$\alpha=(\mathbf{t},q^{\frac{1}{2}})$ as in Section 5.


\subsection{Evaluation formula and duality}

Let $\mathfrak{F}$ and $\mathfrak{F}_\ddagger$ be the auxiliary
kernels as defined by the series expansions 
\eqref{Fseries} and \eqref{Fdaggerseries}, respectively.
By the normalization \eqref{normalizationE} of the Macdonald-Koornwinder
polynomials and \eqref{coefficient}, we have
\begin{equation}\label{evaluation1}
\mathfrak{F}(s_0^{\ddagger\sigma\tau\sigma},s_0^{\ddagger\tau\sigma})
=\frac{G_{\tau\sigma\tau}(s_0^\tau)}{\bigl(1,1\bigr)_{\mathcal{A},\tau}}
\sum_{s\in \mathcal{S}_\tau}
\frac{G_{\tau\sigma\tau}(s)N_{\tau\sigma}(s^{-1})}
{G_{\tau\sigma\tau}(s_0^{\tau})N_{\tau\sigma}(s_0^{\ddagger\tau})}
=\mathfrak{F}_\ddagger(s_0^{\ddagger\sigma\tau\sigma},s_0^{\tau\sigma}).
\end{equation}
By \eqref{Ndecomposition}, \eqref{symmistrivialdiscrete} and by
\eqref{WnulGaussian}, the sum in 
\eqref{evaluation1} can be rewritten as
\begin{equation}\label{sumplus}
\sum_{s\in\mathcal{S}_\tau}
\frac{G_{\tau\sigma\tau}(s)N_{\tau\sigma}(s^{-1})}
{G_{\tau\sigma\tau}(s_0^\tau)N_{\tau\sigma}(s_0^{\ddagger\tau})}=
\sum_{s\in\mathcal{S}_\tau^+}\frac{G_{\tau\sigma\tau}(s)
N_{\tau\sigma}^+(s^{-1})}{G_{\tau\sigma\tau}(s_0^\tau)
N_{\tau\sigma}^+(s_0^{\ddagger\tau})}.
\end{equation}
This sum can be evaluated as follows.
\begin{prop}\label{evaluationprop}
In terms of the Askey-Wilson parametrization \eqref{AWpar} for part
of the difference multiplicity function $\alpha$,
\begin{equation}\label{GNsum}
\sum_{s\in\mathcal{S}_\tau^+}\frac{G_{\tau\sigma\tau}(s)
N_{\tau\sigma}^+(s^{-1})}
{G_{\tau\sigma\tau}(s_0^\tau)N_{\tau\sigma}^+(s_0^{\ddagger\tau})}=
\prod_{i=1}^n\frac{\bigl(qabct^{2(2n-i-1)}/d,
qt^{2(i-n)}/ad;q\bigr)_{\infty}}
{\bigl(qbt^{2(n-i)}/d, qct^{2(n-i)}/d;q\bigr)_{\infty}}.
\end{equation}
\end{prop}
\begin{proof}
We combine \eqref{Nplus} and \eqref{Gaussform} to obtain the 
explicit expression
\begin{equation}\label{lhs}
\begin{split}
\sum_{\lambda\in\Lambda^+}
&\prod_{i=1}^n\left\{\frac{\bigl(qabct^{4(n-i)}/d;q\bigr)_{2\lambda_i}}
{\bigl(abct^{4(n-i)}/d;q\bigr)_{2\lambda_i}}
\left(\frac{-q^{\frac{1}{2}}}{adt^{2(n-i)}}\right)^{\lambda_i}
q^{\lambda_i^2/2}\right.\\
&\left.\qquad\times
\frac{\bigl(abt^{2(n-i)},act^{2(n-i)},
abct^{2(n-i)}/d;q\bigr)_{\lambda_i}}
{\bigl(qt^{2(n-i)}, qbt^{2(n-i)}/d, 
qct^{2(n-i)}/d;q\bigr)_{\lambda_i}}\right\}\\
&\quad\times\prod_{1\leq i<j\leq n}\left\{
\frac{\bigl(qabct^{2(2n-i-j)}/d,
abct^{2(2n-i-j+1)}/d;q\bigr)_{\lambda_i+\lambda_j}}
{\bigl(qabct^{2(2n-i-j-1)}/d,abct^{2(2n-i-j)}/d;
q\bigr)_{\lambda_i+\lambda_j}}
\right.\\
&\left.\qquad\qquad\qquad\qquad\qquad\qquad\qquad\times
\frac{\bigl(qt^{2(j-i)},t^{2(j-i+1)};q\bigr)_{\lambda_i-\lambda_j}}
{\bigl(qt^{2(j-i-1)},t^{2(j-i)};q\bigr)_{\lambda_i-\lambda_j}}\right\}
\end{split}
\end{equation}
for the left hand side of \eqref{GNsum}, where 
$\lambda_i=\langle \lambda,\epsilon_i\rangle$ for $i=1,\ldots,n$.
Now this can be evaluated using the limit case $g_d\rightarrow -\infty$
of van Diejen's \cite[Thm. 2.2]{vD2}
multiple Roger's ${}_6\phi_5$-sum. After straightforward
simplications (see also \cite[(3.7a)]{vD2}), 
we obtain the desired evaluation formula.
\end{proof}
It follows from the product formula \eqref{GNsum}
that the special values \eqref{evaluation1} of the auxiliary kernels
are non-zero (in view of the
generic conditions on the difference multiplicity function
$\alpha$). 

Observe that the spectrum $\mathcal{S}$ satisfies the stability
conditions 
\[ s_\lambda^{\tau\sigma}=s_\lambda^\sigma,\qquad
s_\lambda^{\sigma\tau\sigma}=s_\lambda^{\tau\sigma\tau}=s_\lambda
\]
for all $\lambda\in\Lambda$, since $s\in \mathcal{S}$ only depends
on the values of $\mathbf{t}$ at the inmultipliable roots  $R$ of $R_{nr}$.
In particular, the product formula \eqref{GNsum}
leads to an explicit evaluation
of the non-zero values $\mathfrak{F}(s_0^\ddagger,s_0^{\ddagger\sigma})$
and $\mathfrak{F}_\ddagger(s_0^\ddagger,s_0^\sigma)$. 
Since $G_{\tau}$  (respectively $G_{\sigma\tau}$) is regular
at $s_0^\sigma$ (respectively $s_0$) and $G_{\tau}(s_0^\sigma)\not=0$
(respectively $G_{\sigma\tau}(s_0)\not=0$), we may define {\it normalized
Cherednik kernels} $\mathfrak{E}$ and $\mathfrak{E}_\ddagger$
associated to $\alpha$ and $\alpha_\ddagger$ respectively by
requiring
\begin{equation}\label{fraknormalization}
\mathfrak{E}(s_0^{\ddagger},s_0^{\ddagger\sigma})=1=
\mathfrak{E}_\ddagger(s_0,s_0^{\sigma}).
\end{equation}
This choice of normalization determines the Cherednik kernels $\mathfrak{E}$
and $\mathfrak{E}_\ddagger$ uniquely in view of 
Theorem \ref{Cherednikkernels}. Furthermore, 
by Corollary \ref{convcor}, Lemma \ref{EFconnection} and
\eqref{coefficient}, the normalized Cherednik kernels
$\mathfrak{E}$ and $\mathfrak{E}_\ddagger$ may be explicitly written as
\begin{equation}\label{Eexpansions}
\begin{split}
\mathfrak{E}(\gamma,x)&=G_{\sigma\tau}(\gamma)G_\tau(x)
\sum_{s\in\mathcal{S}_\tau}\mu_\tau(s)E_{\sigma\tau}(s;\gamma)E_\tau(s;x),\\
\mathfrak{E}_\ddagger(\gamma,x)&=
G_{\sigma\tau}(\gamma)G_\tau(x)\sum_{s\in\mathcal{S}_\tau}
\mu_\tau(s)E_{\sigma\tau}(s;\gamma^{-1})E_{\ddagger\tau}(s^{-1};x),
\end{split}
\end{equation}
with $\mu=\mu_{\alpha}\in\mathcal{F}(\mathcal{S})$
defined by
\begin{equation}\label{mutau}
\mu(s)=C_0\frac{G_{\sigma\tau}(s)N_{\sigma}(s^{-1})}
{G_{\sigma\tau}(s_0)N_{\sigma}(s_0^{\ddagger})},
\qquad s\in\mathcal{S}
\end{equation}
and with the (non-zero) normalization constant
$C_0=C_0^{\alpha}\in\mathbb{C}^\times$ chosen in such a way
that \eqref{fraknormalization} holds true. In other words, $C_0$
is chosen in such a way that
\begin{equation}\label{summutau}
\sum_{s\in\mathcal{S}}\mu(s)=
\frac{1}{G_{\tau\sigma\tau}(s_0^\tau)G(s_0^{\tau\sigma})}.
\end{equation}
By Proposition \ref{evaluationprop}
and Remark \ref{Gaussianevaluation},
the constant $C_0$ can be evaluated explicity in terms of
the Askey-Wilson parametrization \eqref{AWpar} of part of the difference
multiplicity function $\alpha$,
\begin{equation}\label{C0plus}
C_0=\prod_{i=1}^n\frac{\bigl(adt^{2(n-i)}, bdt^{2(n-i)},
cdt^{2(n-i)},bct^{2(n-i)},dt^{2(i-n)}/a;q\bigr)_{\infty}}
{\bigl(abcdt^{2(2n-i-1)};q\bigr)_{\infty}}.
\end{equation}
Combined with \eqref{Ndecomposition},
\eqref{Nplus} and \eqref{Gaussform}, this entails an
explicit expression for the weight
$\mu\in\mathcal{F}(\mathcal{S})$
in terms of $q$-shifted factorials, cf. the expression \eqref{lhs}.
\begin{thm}[Duality]\label{dualitytheorem}
The normalized Cherednik kernels $\mathfrak{E}$ and 
$\mathfrak{E}_\ddagger$ associated to $\alpha$
and $\alpha_\ddagger$ respectively, satisfy the duality property
\[\mathfrak{E}(\gamma,x)=\mathfrak{E}_\sigma(x,\gamma),\qquad
\mathfrak{E}_\ddagger(\gamma,x)=\mathfrak{E}_{\ddagger\sigma}(x,\gamma).
\]
\end{thm}
\begin{proof}
This follows immediately from the characterization Theorem 
\ref{Cherednikkernels}
for Che\-red\-nik kernels, the fact that $\psi^{-1}=\psi_\sigma$ for the
duality anti-isomorphism $\psi:\mathcal{H}\rightarrow
\mathcal{H}_\sigma$ and the fact that the chosen normalization
of the Che\-red\-nik kernels is self-dual.
\end{proof}


\subsection{The polynomial reduction}

Let $v\in \mathcal{S}$ be a spectral point of the $Y$-operators
$Y_i\in\mathcal{H}$ ($i=1,\ldots,n$), 
considered as endomorphisms of $\mathcal{A}$.
Then we have observed that the Macdonald-Koornwinder
polynomial $E(v;\cdot)\in\mathcal{A}$ is the unique Laurent polynomial
satisfying
\begin{equation}\label{polkara}
 p(Y)E(v;\cdot)=p(v)E(v;\cdot),\qquad \forall\, p\in\mathcal{A}
\end{equation}
and satisfying the normalization $E(v;s_0^{\ddagger\sigma})=1$.

On the other hand, the normalized Cherednik function 
$\mathfrak{E}(\gamma,x)$
associated to $\alpha$  is regular at $\gamma=v^{-1}$, and the
meromorphic function $\mathfrak{E}(v^{-1},\cdot)\in\mathcal{M}$ satisfies
\[ p(Y)\mathfrak{E}(v^{-1};\cdot)=p(v)\mathfrak{E}(v^{-1},\cdot),\qquad
\forall\, p\in\mathcal{A}.
\]
Whence it is tempting to believe that $\mathfrak{E}(v^{-1},\,\cdot)$
is a constant multiple of the normalized Macdonald-Koornwinder
polynomial $E(v;\,\cdot)$. In this section, 
we show that this is indeed the case.
The argument essentially amounts to the fact that 
the meromorphic common eigenfunction $\mathfrak{E}(v^{-1},\,\cdot)$
of the $Y$-operators is ``regular enough'' to ensure that 
$\mathfrak{E}(v^{-1},\,\cdot)\in\mathcal{A}$ by standard elliptic
function theory, and hence it can only be a 
constant multiple of the Macdonald-Koornwinder polynomial $E(v;\,\cdot)$.

To make the arguments rigorous and transparent, we 
study properties of certain explicit function transforms in detail,
following similar lines of reasoning as for the auxiliary transforms
(see Subsections 5.2 \& 5.3).

The starting point for the definitions of these transforms
forms the explicit series expansions \eqref{Eexpansions} for the
normalized Cherednik kernels $\mathfrak{E}$ and $\mathfrak{E}_\ddagger$.
The bound on the associated weight 
$\mu_\tau\in\mathcal{F}(\mathcal{S}_\tau)$ as derived in 
the proof of Corollary \ref{convcor} allows us to define
the following two transforms in terms of series expansions which 
absolutely and uniformly converge
on compacta of $(\mathbb{C}^\times)^n$.
\begin{Def}\label{Hdef}
We define linear mappings $H,H_{\ddagger}: \mathcal{A}\rightarrow
\mathcal{O}\,G_{\tau}\subset\mathcal{M}$ by
\begin{equation}\label{H}
\begin{split}
\bigl(Hp\bigr)(x)&=G_{\tau}(x)
\sum_{s\in\mathcal{S}_{\tau}}\mu_\tau(s)p(s^{-1})
E_\tau(s;x),\\
\bigl(H_{\ddagger}p\bigr)(x)&=G_{\tau}(x)
\sum_{s\in\mathcal{S}_{\tau}}\mu_\tau(s)p(s^{-1})
E_{\ddagger\tau}(s^{-1};x)
\end{split}
\end{equation}
for all $p\in\mathcal{A}$.
\end{Def}
The transforms $H$ and $H_\ddagger$ are linked to
$\mathfrak{E}(v^{-1},\,\cdot)$ and $\mathfrak{E}_\ddagger(v,\,\cdot)$
for $v\in\mathcal{S}$ as follows. Observe that 
the factor $E_{\sigma\tau}(s;v^{-1})$
for $s\in\mathcal{S}_\tau$ occurring in the explicit expansion 
sum \eqref{Eexpansions} for $\mathfrak{E}(v^{-1},\,\cdot)$ can be
rewritten by the duality \eqref{polduality} as
\[E_{\sigma\tau}(s;v^{-1})=E_{\sigma\tau\sigma}(v;s^{-1})=
E_{\tau\sigma\tau}(v;s^{-1})
\]
since $\sigma\tau\sigma=\tau\sigma\tau$ when acting upon
difference multiplicity functions, and 
$s_\lambda^\tau=s_\lambda^{\sigma\tau}$,
$s_\lambda=s_\lambda^{\sigma\tau\sigma}$ for all $\lambda\in\Lambda$.
This implies that  
\begin{equation}\label{formulatransform1}
\mathfrak{E}(v^{-1},x)
=G_{\sigma\tau}(v)
H\bigl(E_{\tau\sigma\tau}(v;\,\cdot)\bigr)(x),\qquad \forall\,
v\in \mathcal{S}
\end{equation}
with $H:\mathcal{A}\rightarrow \mathcal{M}$ the linear
map defined in \eqref{H}.
In exactly the same fashion, we can write
\begin{equation}\label{formulatransform2}
\mathfrak{E}_\ddagger(v^{-1},x)=G_{\sigma\tau}(v)
H_\ddagger\bigl(E_{\tau\sigma\tau}(v^{-1};\,\cdot)\bigr)(x),
\qquad \forall\, v\in\mathcal{S}_\ddagger
\end{equation}
with $H_\ddagger:\mathcal{A}\rightarrow \mathcal{M}$
the linear map defined in \eqref{H}.
The main step now is to prove that the images of the transforms
$H$ and $H_\ddagger$ are contained in $\mathcal{A}$. For this, we
first compute the intertwining properties of $H$ and $H_{\ddagger}$
under the action of the double affine Hecke algebra.
Recall the anti-isomorphisms $\kappa$ and $\kappa^I$ defined
in Lemma \ref{EFconnection}.
\begin{lem}\label{transfH}
{\bf a)}
For $X\in\mathcal{H}_{\sigma\tau\sigma}=\mathcal{H}_{\tau\sigma\tau}$,
\begin{equation}\label{auxH1} 
H\circ
X=\bigl(\dagger_{\ddagger}\circ\kappa_{\tau\sigma\tau}^I\circ
\iota_{\tau\sigma\tau}\bigr)(X)\circ H,
\end{equation}
with the double affine Hecke algebra acting as $q$-difference
reflection operators on both sides.\\
{\bf b)} For
$X\in\mathcal{H}_{\sigma\tau\sigma}=\mathcal{H}_{\tau\sigma\tau}$,
\begin{equation}\label{auxH2}
H_{\ddagger}\circ X=
\bigl(\kappa_{\ddagger\tau\sigma\tau}
\circ\ddagger_{\tau\sigma\tau}\bigr)(X)\circ H_{\ddagger},
\end{equation}
with the double affine Hecke algebra acting as $q$-difference
operators on the left hand side and as $q^{-1}$-difference operators
on the right hand side.
\end{lem}
\begin{proof}
{\bf a)} Up to an (irrelevant) multiplicative constant, $Hp$
can be written as
\[ \bigl(Hp\bigr)(x)=G_\tau(x)\lbrack G_{\tau\sigma\tau}p, 
\mathfrak{E}_{\mathcal{A},\tau}(\cdot\,,x)\rbrack_{\mathcal{A},\tau\sigma}
\]
by \eqref{mutau}, see Subsections 4.2 \& 4.3 for the notations.
We thus obtain for $X\in \mathcal{H}_{\tau\sigma\tau}$
\[ H\circ X=
\bigl(\tau_\tau\circ\psi_{\tau\sigma}\circ\iota_{\tau\sigma}
\circ\tau_{\tau\sigma\tau}\bigr)(X)\circ H
\]
by Proposition \ref{Gaussianconj}, Proposition \ref{adjoint2},
Remark \ref{adjoint2remark}, Proposition \ref{allrelations} 
and Lemma \ref{actioncompatibleH}.
Now computing the isomorphism
$\tau_\tau\circ\psi_{\tau\sigma}\circ\iota_{\tau\sigma}
\circ\tau_{\tau\sigma\tau}: \mathcal{H}_{\tau\sigma\tau}\rightarrow
\mathcal{H}$ on the algebraic generators $T_j^{\tau\sigma\tau}$
($j=0,\ldots,n$) and $x_1$ of 
$\mathcal{H}_{\tau\sigma\tau}$ and comparing the outcome with Lemma
\ref{kappalem} shows that
\[ \tau_\tau\circ\psi_{\tau\sigma}\circ\iota_{\tau\sigma}
\circ\tau_{\tau\sigma\tau}=\dagger_{\ddagger}\circ
\kappa_{\tau\sigma\tau}^I\circ\iota_{\tau\sigma\tau}
\]
(both sides map $T_j^{\tau\sigma\tau}$ to $T_j$ for $j=0,\ldots,n$ and
$x_1$ to $q^{\frac{1}{2}}T_0x_1^{-1}U^{-1}$).\\
{\bf b)} By \eqref{Gaussianinversion}, we can write $H_{\ddagger}p$
up to an (irrelevant) multiplicative constant as
\[ \bigl(H_{\ddagger}p\bigr)(x)=
G_\ddagger(x)^{-1}\lbrack
G_{\tau\sigma\tau}p,I\mathfrak{E}_{\mathcal{A},\ddagger\tau}(\cdot\,,x)
\rbrack_{\mathcal{A},\tau\sigma},
\]
where $(Ig)(s)=g(s^{-1})$ for 
$g\in \mathcal{F}(\mathcal{S}_\tau,\mathcal{A})$ and
$s\in\mathcal{S}_{\ddagger\tau}$. A similar computation as for
{\bf a)} now shows that
\[ H_\ddagger\circ X=
\bigl(\tau_\dagger^{-1}\circ
\psi_{\ddagger\tau\sigma}\circ\dagger_{\tau\sigma}\circ\iota_{\tau\sigma}
\circ\tau_{\tau\sigma\tau}\bigr)(X)\circ H_\ddagger
\]
for $X\in\mathcal{H}_{\tau\sigma\tau}$.
The lemma follows from the fact that
\[ \tau_\dagger^{-1}\circ
\psi_{\ddagger\tau\sigma}\circ\dagger_{\tau\sigma}\circ\iota_{\tau\sigma}
\circ\tau_{\tau\sigma\tau}=\kappa_{\ddagger\tau\sigma\tau}\circ
\ddagger_{\tau\sigma\tau}
\]
as unital algebra isomorphisms from $\mathcal{H}_{\tau\sigma\tau}$ onto 
$\mathcal{H}_{\ddagger}$, which again can be checked by computing
the images of the algebraic generators $T_i^{\tau\sigma\tau}$
($i=1,\ldots,n$), $Y_1^{\tau\sigma\tau}$ and $x_1$ of
$\mathcal{H}_{\tau\sigma\tau}$ explicitly for both isomorphisms.
\end{proof}

Note the similarities between the maps $H,H_\ddagger:
\mathcal{A}\rightarrow \mathcal{M}$ and the auxiliary
transforms $L,L_{\ddagger}:\mathcal{A}\rightarrow \mathcal{O}$
defined before: they nearly satisfy the same intertwining properties
with respect to the action of the double affine Hecke algebra (up to
some elementary (anti-)isomorphisms $\iota,\ddagger,\dagger$).
For the auxiliary transforms $L,L_\ddagger:\mathcal{A}\rightarrow
\mathcal{O}$ we used the intertwining properties to prove 
that $L1$ and $L_{\ddagger}1$ are 
constant multiples of $1\in\mathcal{A}$, 
see Lemma \ref{L}. A similar argument can now
be applied to the transforms $H,H_{\ddagger}: \mathcal{A}\rightarrow
\mathcal{M}$. 
\begin{lem}\label{H1}
The transforms $H,H_\ddagger:\mathcal{A}\rightarrow \mathcal{M}$
satisfy
\[ H1=G_{\sigma\tau}(s_0)^{-1}1=H_{\ddagger}1,
\]
where $1\in\mathcal{A}$ is the Laurent polynomial identically equal to
one.
\end{lem}
\begin{proof}
In view of Lemma \ref{transfH} we can repeat the
arguments of the proof of Lemma \ref{L} to prove that
$H1\in\mathcal{M}$ is $\mathcal{W}$-invariant
(with $\mathcal{W}$ acting as constant coefficient $q$-difference
reflection operators) and that $H_{\ddagger}1\in\mathcal{M}$
is $\mathcal{W}$-invariant (with $\mathcal{W}$ now acting as
constant coefficient $q^{-1}$-difference reflection operators).
So in terms of exponential coordinates 
\[x=e^{2\pi i w}=\bigl(e^{2\pi i w_1},e^{2\pi i w_2},\ldots,e^{2\pi i
w_n}\bigr),
\]
we obtain meromorphic functions
\[ w\mapsto H(1)\bigl(e^{2\pi i w}\bigr),\qquad
w\mapsto H_{\ddagger}(1)\bigl(e^{2\pi i w}\bigr)
\]
on the compact torus $\mathbb{C}^n/\Gamma_q$, with
$\Gamma_q=\mathbb{Z}^n+\mathbb{Z}^n\,v$  and $v$ an element in the
upper half plane satisfying $q=e^{2\pi i v}$.
Now both the maps $w\mapsto H(1)\bigl(e^{2\pi i w}\bigr)
G_\tau(e^{2\pi iw})^{-1}$ and $w\mapsto  
H_{\ddagger}(1)\bigl(e^{2\pi i w}\bigr)G_\tau(e^{2\pi iw})^{-1}$
are analytic in $w\in \mathbb{C}^n$. Since 
\[G_\tau(e^{2\pi i w})=
\prod_{j=1}^n\frac{1}{\bigl(-q^{\frac{1}{2}}u_0t_0^{-1}e^{2\pi i w_j},
-q^{\frac{1}{2}}u_0t_0^{-1}e^{-2\pi i w_j};q\bigr)_{\infty}},
\]
it now follows that the function $w_j\mapsto H(1)(e^{2\pi i w})$ and 
$w_j\mapsto H_{\ddagger}(1)(e^{2\pi i
w})$ for fixed, regular $w_k\in\mathbb{C}$ ($k\not=j$)
is elliptic with period lattice $\mathbb{Z}+\mathbb{Z}\,v$,
whose possible poles are at most simple and located at
\begin{equation}\label{polesform}
\bigl(u+\mathbb{Z}+\mathbb{N}\,v\bigr)\cup
\bigl(-u+\mathbb{Z}-\mathbb{N}\,v\bigr),
\end{equation}
where $u\in\mathbb{C}$ is chosen such that
$e^{2\pi i u}=-q^{-\frac{1}{2}}u_0^{-1}t_0$. Standard elliptic
function theory now implies that the functions
$w_j\mapsto H(1)(e^{2\pi i w})$ and 
$w_j\mapsto H_{\ddagger}(1)(e^{2\pi i
w})$ are constant. Consequently, the meromorphic functions
$H1, H_{\ddagger}1\in\mathcal{M}$, regarded as
analytic functions on the open, dense set of elements in
$(\mathbb{C}^\times)^n$ which do not belong to the zero set of
$G_\tau^{-1}\in\mathcal{O}$, are constant. Hence 
$H1, H_{\ddagger}1\in\mathcal{O}$, and they are constant on 
$(\mathbb{C}^\times)^n$. 

To complete the proof, we note that
\begin{equation*}
\begin{split}
H(1)(s_0^{\ddagger\sigma})&=G_\tau(s_0^{\sigma})\sum_{s\in\mathcal{S}_\tau}
\mu_\tau(s)\\
&=G_{\sigma\tau}(s_0)^{-1}
\mathfrak{E}(s_0^{\ddagger},s_0^{\ddagger\sigma})
=G_{\sigma\tau}(s_0)^{-1}
\end{split}
\end{equation*}
by the normalization of the Mac\-do\-nald-Koorn\-win\-der 
polynomials and of the
Che\-red\-nik kernel. Similarly one can show that $H_\ddagger(1)(s_0^\sigma)=
G_{\sigma\tau}(s_0)^{-1}$.
\end{proof}

The following corollary is immediate from Lemma \ref{transfH} 
and from (the proof of) Lemma \ref{H1}.
\begin{cor}\label{uniqueness}
The linear maps $H,H_\ddagger: \mathcal{A}\rightarrow \mathcal{O}G_{\tau}$
are, up to a multiplicative constant, 
uniquely determined by the intertwining properties
\eqref{auxH1} and \eqref{auxH2} respectively under the action of 
the double affine Hecke algebra $\mathcal{H}_{\tau\sigma\tau}$.
\end{cor}

Continuing the same line of arguments as for auxiliary transforms, 
we can now state the following consequence of Lemma \ref{H1}.
\begin{cor}\label{generalizedeigenvalue2}
{\bf a)} The transform $H$ maps into $\mathcal{A}$.
Furthermore,
\[ H\bigl(E_{\tau\sigma\tau}(v;\cdot)\bigr)=e(v)E(v;\cdot),\qquad
\forall\,v\in \mathcal{S}=\mathcal{S}_{\tau\sigma\tau}
\]
for some $e\in \mathcal{F}(\mathcal{S})$.\\
{\bf b)} The transform $H_{\ddagger}$ maps into
$\mathcal{A}$. Furthermore,
\[ H_{\ddagger}\bigl(E_{\tau\sigma\tau}(v^{-1};\cdot)\bigr)=
e_{\ddagger}(v)E_{\ddagger}(v;\cdot),\qquad \forall\,
v\in\mathcal{S}_\ddagger
\]
for some $e_{\ddagger}\in \mathcal{F}(\mathcal{S}_\ddagger)$.
\end{cor}
\begin{proof}
By Lemma \ref{transfH} and Lemma \ref{H1}, we clearly have
$H(\mathcal{A})\subseteq \mathcal{A}$ and
$H_{\ddagger}(\mathcal{A})\subseteq \mathcal{A}$.
Combining this fact with the formulas
\begin{equation*}
\begin{split}
\bigl(\dagger_{\ddagger}\circ
\kappa_{\tau\sigma\tau}^I\circ\iota_{\tau\sigma\tau}\bigr)
(Y_i^{\tau\sigma\tau})&=Y_i,\\
\bigl(\kappa_{\ddagger\tau\sigma\tau}\circ\ddagger_{\tau\sigma\tau}\bigr)
(Y_i^{\tau\sigma\tau})&=Y_i^{\ddagger}{}^{-1}
\end{split}
\end{equation*}
for $i=1,\ldots,n$, and using Lemma \ref{transfH}
and the results of Subsection 4.2, we directly obtain the
second statement of the corollary (compare with the proof of 
Corollary \ref{auxcor}).
\end{proof}

Corollary \ref{generalizedeigenvalue2} combined with
\eqref{formulatransform1} and \eqref{formulatransform2} already
prove that $\mathfrak{E}(s^{-1},\cdot)$ for $s\in\mathcal{S}$ (respectively 
$\mathfrak{E}_\ddagger(v^{-1},\cdot)$ for $v\in\mathcal{S}_\ddagger$)
is a constant multiple of $E(s;\cdot)$ (respectively
$E_\ddagger(v;\cdot)$). The next aim is to prove that the constant
multiple is in fact always equal to one.
For the proof of this result we first need to compute the functions 
$e\in \mathcal{F}(\mathcal{S})$
and $e_\ddagger\in \mathcal{F}(\mathcal{S}_\ddagger)$,
see Corollary \ref{generalizedeigenvalue2}. Again,
the method is similar to the computation of the generalized
eigenvalues of the auxiliary transforms $L$ and $L_\ddagger$,
respectively. 
\begin{lem}\label{generalizedeigenvalue3}
The maps $e\in \mathcal{F}(\mathcal{S})$ and 
$e_\ddagger\in \mathcal{F}(\mathcal{S}_\ddagger)$
are given explicitly by
\[ e(\cdot)=G_{\sigma\tau}(\cdot)^{-1}|_{\mathcal{S}},\qquad
e_\ddagger(\cdot)=G_{\sigma\tau}(\cdot)^{-1}|_{\mathcal{S}_\ddagger}.
\]
\end{lem} 
\begin{proof}
Let $\widetilde{e}\in \mathcal{F}(\mathcal{S})$ and let
$\widetilde{H}:\mathcal{A}\rightarrow \mathcal{A}$ be the linear
map defined by
\[
\widetilde{H}(E_{\tau\sigma\tau}(v;\,\cdot)\bigr)=
\widetilde{e}(v)E(v;\,\cdot),\qquad
\forall\,v\in\mathcal{S}.
\]
Suppose furthermore that $\widetilde{e}$, considered as
multiplication operator on $\mathcal{F}(\mathcal{S})=
\mathcal{F}(\mathcal{S}_{\sigma\tau\sigma})$,
induces an algebra isomorphism $\widetilde{\nu}_{\ddagger\sigma\tau}:
\mathcal{H}_{\ddagger\sigma\tau}\rightarrow
\mathcal{H}_{\ddagger\sigma}$ by the formula
\[ X\circ \widetilde{e}=\widetilde{e}
\circ \widetilde{\nu}_{\ddagger\sigma\tau}(X),\qquad \forall\,X\in
\mathcal{H}_{\ddagger\sigma\tau}
\]
in $\hbox{End}_{\mathbb{C}}(\mathcal{F}(\mathcal{S}))$, where 
$\mathcal{H}_{\ddagger\sigma\tau}$ and $\mathcal{H}_{\ddagger\sigma}$ 
act via the dot-action on $\mathcal{F}(\mathcal{S})=
\mathcal{F}(\mathcal{S}_{\sigma\tau\sigma})$.

By repeating the arguments of the proof of Lemma \ref{eigenconj},
we deduce that if the isomorphism 
$\widetilde{\nu}_{\ddagger\sigma\tau}$ satisfies
\begin{equation}\label{ehio}
\psi_\sigma\circ\dagger_{\ddagger\sigma}\circ
\widetilde{\nu}_{\ddagger\sigma\tau}
\circ\dagger_{\sigma\tau}\circ\psi_{\sigma\tau\sigma}=
\dagger_{\ddagger}\circ\kappa_{\tau\sigma\tau}^I\circ\iota_{\tau\sigma\tau},
\end{equation}
then the linear map $\widetilde{H}$ has the intertwining property
\begin{equation}\label{dhio}
\widetilde{H}\circ
X=\bigl(\dagger_{\ddagger}\circ\kappa_{\tau\sigma\tau}^I\circ
\iota_{\tau\sigma\tau}\bigr)(X)\circ \widetilde{H}
\end{equation}
for all
$X\in\mathcal{H}_{\sigma\tau\sigma}=\mathcal{H}_{\tau\sigma\tau}$.
Note now that \eqref{ehio} is equivalent to
\begin{equation*}
\begin{split}
\widetilde{\nu}_{\ddagger\sigma\tau}&=
\dagger_\sigma\circ\psi\circ\dagger_{\ddagger}\circ\kappa_{\tau\sigma\tau}^I
\circ\iota_{\tau\sigma\tau}\circ\psi_{\sigma\tau}\circ
\dagger_{\ddagger\sigma\tau}\\
&=\tau_{\dagger\sigma}^{-1},
\end{split}
\end{equation*}
where the last equality follows from computing both sides explicitly
on suitable algebraic generators of
$\mathcal{H}_{\ddagger\sigma\tau}$. By Proposition \ref{Gaussianconj},
Lemma \ref{actioncompatibleH} and \eqref{Gaussianinversion}, 
we conclude that if
$\widetilde{e}$ equals $G_{\ddagger\sigma}(\cdot)|_{\mathcal{S}}=
G_{\sigma\tau}(\cdot)^{-1}|_{\mathcal{S}}$
up to a multiplicative constant, then $\widetilde{e}$ induces the
isomorphism $\nu_{\ddagger\sigma\tau}$, and hence the corresponding
linear map $\widetilde{H}$ has the intertwining property \eqref{dhio}.
By the uniqueness of such linear maps (see Corollary
\ref{uniqueness}), we conclude that these are in fact the only possibilities
for $\widetilde{e}\in\mathcal{F}(\mathcal{S})$ 
for which the associated linear maps
$\widetilde{H}$ have the intertwining property \eqref{dhio}. 
Applying this result to $H$ and using
the normalization $H1=G_{\sigma\tau}(s_0)^{-1}1$, we 
conclude that $e(\cdot)=G_{\sigma\tau}(\cdot)^{-1}|_{\mathcal{S}}$.

The proof for $e_\ddagger$ is similar; the arguments lead to the
transformation behaviour 
\begin{equation*}
\begin{split}
X\circ e_{\ddagger}&=e_{\ddagger}\circ\bigl(\dagger_{\ddagger\sigma}
\circ\psi_{\ddagger}\circ\kappa_{\ddagger\tau\sigma\tau}\circ
\ddagger_{\tau\sigma\tau}\circ\psi_{\sigma\tau}\bigr)(X)\circ
e_{\ddagger}\\
&=e_{\ddagger}\circ
\tau_{\sigma\tau}(X),\qquad\qquad\qquad\qquad\qquad
\qquad\qquad\forall\,X\in\mathcal{H}_{\sigma\tau}
\end{split}
\end{equation*}
in $\hbox{End}_{\mathbb{C}}(\mathcal{F}(\mathcal{S}_\ddagger))$, where
the function $e_\ddagger$ is considered as a multiplication operator
in the endomorphism space
$\hbox{End}_{\mathbb{C}}(\mathcal{F}(\mathcal{S}_\ddagger))=
\hbox{End}_{\mathbb{C}}(\mathcal{F}(\mathcal{S}_{\ddagger\sigma\tau\sigma}))$.
It follows from this that
$e_{\ddagger}(\cdot)=G_{\sigma\tau}(\cdot)^{-1}|_{\mathcal{S}_\ddagger}$.
\end{proof}

We are now in the position to 
prove that the normalized Cherednik kernel  
meromorphically extend the
Maconald-Koornwinder polynomials in their degrees.

\begin{thm}\label{polreduc}
{\bf a)} Let $\mathfrak{E}$ be the normalized Cherednik kernel
associated to $\alpha$.
For all $v\in \mathcal{S}$,
\[\mathfrak{E}(v^{-1},\,\cdot)=E(v;\,\cdot),
\]
where $E(v;\,\cdot)=E_\alpha(v;\,\cdot)$ 
is the normalized Macdonald-Koornwinder polynomial
corresponding to the spectral point $v\in\mathcal{S}$.\\
{\bf b)} Let $\mathfrak{E}_\ddagger$ be the normalized Cherednik
kernel associated to $\alpha_\ddagger$.
For all $v\in \mathcal{S}_\ddagger$,
\[\mathfrak{E}_\ddagger(v^{-1},\,\cdot)=E_\ddagger(v;\,\cdot),
\]
where $E_\ddagger(v;\,\cdot)=E_{\alpha_\ddagger}(v;\,\cdot)$ 
is the normalized Macdonald-Koornwinder polynomial
corresponding to the spectral point $v\in\mathcal{S}_\ddagger$.
\end{thm}
\begin{proof}
By formula \eqref{formulatransform1}, Corollary
\ref{generalizedeigenvalue2} and Lemma \ref{generalizedeigenvalue3},
we have 
\begin{equation*}
\begin{split}
\mathfrak{E}(v^{-1},\,\cdot)
&=G_{\sigma\tau}(v)H\bigl(E_{\tau\sigma\tau}(v;\cdot)\bigr)\\
&=G_{\sigma\tau}(v)e(v)E(v;\cdot)\\
&=E(v;\cdot)
\end{split}
\end{equation*}
for $v\in\mathcal{S}$.
The proof for $\mathfrak{E}_\ddagger$ is similar.
\end{proof}

Several direct consequences can be derived by combining the 
explicit series expansion \eqref{Eexpansions} of the normalized
Cherednik kernel with its polynomial reduction (see Theorem \ref{polreduc}). 

\begin{cor}\label{consequenceab}
{\bf a)} The expansion of the inverse Gaussian
$G_\tau^{-1}\in\mathcal{O}$ in terms of Macdonald-Koornwinder 
polynomials is given by
\begin{equation*}
\begin{split}
G_\tau(x)^{-1}&=G_{\sigma\tau}(s_0)
\sum_{s\in\mathcal{S}_\tau}\mu_\tau(s)E_\tau(s;x)\\
&=G_{\sigma\tau}(s_0)
\sum_{s\in\mathcal{S}_\tau}\mu_\tau(s)E_{\ddagger\tau}(s^{-1};x).
\end{split}
\end{equation*}
Here the series converge absolutely and uniformly 
for $x$ in compacta of $(\mathbb{C}^\times)^n$.\\
{\bf b)} For all $s\in\mathcal{S}=\mathcal{S}_{\sigma\tau\sigma}$, we
have the identities 
\begin{equation*}
\begin{split}
E_{\sigma\tau\sigma}(s;Y^\tau{}^{-1})\bigl(G_\tau^{-1}\bigr)&=
\frac{G_{\sigma\tau}(s_0)}{G_{\sigma\tau}(s)}E(s;\cdot)G_\tau^{-1},\\
E_{\sigma\tau\sigma}(s;Y^{\ddagger\tau})\bigl(G_\tau^{-1}\bigr)&=
\frac{G_{\sigma\tau}(s_0)}{G_{\sigma\tau}(s)}
E_\ddagger(s^{-1};\cdot)G_\tau^{-1}
\end{split}
\end{equation*}
in $\mathcal{O}$.
\end{cor}
\begin{proof}
{\bf a)} By Theorem \ref{polreduc}, the normalized Cherednik kernels
$\mathfrak{E}$ and $\mathfrak{E}_\ddagger$ satisfy
\begin{equation}\label{nul}
\mathfrak{E}(s_0^{-1},\cdot)=1=\mathfrak{E}_\ddagger(s_0,\cdot).
\end{equation}
The identities now follow by substituting the explicit series 
expansions \eqref{Eexpansions} for 
$\mathfrak{E}$ and $\mathfrak{E}_\ddagger$ into \eqref{nul}.\\
{\bf b)} We prove the first equality, the second is similar.
By {\bf a)}, by the polynomial duality
$E_{\sigma\tau\sigma}(s;v^{-1})=E_{\sigma\tau}(v;s^{-1})$ for
$v\in\mathcal{S}_\tau=\mathcal{S}_{\sigma\tau}$ and
$s\in\mathcal{S}=\mathcal{S}_{\sigma\tau\sigma}$,
and by the explicit expansion \eqref{Eexpansions}
for the Cherednik kernel $\mathfrak{E}$,
\begin{equation*}
\begin{split}
E_{\sigma\tau\sigma}(s;Y^\tau{}^{-1})\bigl(G_\tau^{-1}\bigr)&=
G_{\sigma\tau}(s_0)\sum_{v\in\mathcal{S}_\tau}\mu_\tau(v)
E_{\sigma\tau\sigma}(s;Y^\tau{}^{-1})\bigl(E_\tau(v;\cdot)\bigr)\\
&=G_{\sigma\tau}(s_0)\sum_{v\in\mathcal{S}_\tau}\mu_\tau(v)
E_{\sigma\tau\sigma}(s;v^{-1})E_\tau(v;\cdot)\\
&=G_{\sigma\tau}(s_0)\sum_{v\in\mathcal{S}_\tau}\mu_\tau(v)
E_{\sigma\tau}(v;s^{-1})E_\tau(v;\cdot)\\
&=\frac{G_{\sigma\tau}(s_0)}{G_{\sigma\tau}(s)}\mathfrak{E}(s^{-1},\cdot)
G_\tau^{-1}.
\end{split}
\end{equation*}
The formula then follows from the polynomial reduction
$\mathfrak{E}(s^{-1},\cdot)=E(s;\cdot)$ of the Cherednik kernel
$\mathfrak{E}$, see Theorem \ref{polreduc}.
\end{proof}
 
\begin{rem}
The series expansion of $G_\tau^{-1}\in\mathcal{O}$ in terms of 
Macdonald-Koornwinder polynomials (see Corollary
\ref{consequenceab}{\bf a)}) 
may be regard as a generalization of the Jacobi triple product
identity, cf. Remark \ref{remabc}{\bf c)}.
\end{rem}


\subsection{Symmetrization}

Recall that the symmetrizer $C_+\in H_0\subset \mathcal{H}$ was used
in Subsection 4.5 to derive $W_0$-invarariant versions of all main results on
Macdonald-Koornwinder polynomials. In this subsection we 
consider the action of the symmetrizer $C_+$
on the normalized Cherednik kernel $\mathfrak{E}$.

\begin{Def}
Let $\mathfrak{E}$ be the normalized Cherednik kernel associated to 
$\alpha$. The meromorphic kernel 
$\mathfrak{E}^+(\cdot,\cdot)=\mathfrak{E}_\alpha^+(\cdot,\cdot)\in
\mathcal{M}\bigl((\mathbb{C}^\times)^n\times (\mathbb{C}^\times)^n\bigr)$
defined by
\[\mathfrak{E}^+(\gamma,x)=\bigl(C_+\mathfrak{E}(\gamma,\cdot)\bigr)(x)
\]
is called the symmetric Cherednik kernel associated to $\alpha$.
\end{Def}

In the following lemma we give some elementary properties of the 
symmetric Cherednik kernel. 
\begin{lem}\label{elementary}
Let $\mathfrak{E}^+$ be the symmetric Cherednik kernel
associated to $\alpha$.\\
{\bf a)} $\mathfrak{E}^+(\gamma,x)$ is
$W_0$-invariant in $x$.\\
{\bf b)} $\mathfrak{E}^+(\gamma,x)=
\bigl(C_+^{\sigma}\mathfrak{E}(\cdot,x)\bigr)(\gamma)$.\\
{\bf c)} $\mathfrak{E}^+(\gamma,x)$ is
$W_0$-invariant in $\gamma$.\\
{\bf d)} $\mathfrak{E}^+(s_0,s_0^{\sigma})=1$ 
\textup{(}normalization\textup{)}.
\end{lem}
\begin{proof}
{\bf a)} This follows from the identities $T_iC_+=t_iC_+$
for $i=1,\ldots,n$ in $\mathcal{H}$ and the explicit form of the
reflection operators $T_i$, compare with the proof of Lemma \ref{L}.\\
{\bf b)} This follows from \eqref{stableC} and from
the transformation behaviour of the Cherednik kernel
under the action of the double 
affine Hecke algebra $\mathcal{H}$, see Theorem \ref{Cherednikkernels}.\\
{\bf c)} This is immediately clear from {\bf b)} and from the proof of
{\bf a)}.\\
{\bf d)} For any $g\in \mathcal{O}G_{\sigma\tau}$, we have 
\[\bigl(T_i^{\sigma}g\bigr)(s_0^\ddagger)=t_ig(s_0^\ddagger),\qquad
i=1,\ldots,n
\]
since the rational function $c_{a_i}^\sigma\in\mathbb{C}(x)$
occurring in the definition of $T_i^\sigma$ vanishes at
$s_0^\ddagger$, see Lemma \ref{actioncompatibleW}{\bf b)}.
In particular, $(C_+^{\sigma}g)(s_0^\ddagger)=g(s_0^\ddagger)$
for $g\in\mathcal{O}G_{\sigma\tau}$. It follows that for generic
$x\in (\mathbb{C}^\times)^n$ (in particular, for $x=s_0^\sigma$),
\[
\mathfrak{E}^+(s_0^\ddagger,x)=
\bigl(C_+^{\sigma}\mathfrak{E}(\cdot,x)\bigr)(s_0^\ddagger)
=\mathfrak{E}(s_0^{\ddagger},x)=1,
\]
where the last equality follows from the polynomial reduction
of $\mathfrak{E}$, see Theorem \ref{polreduc}. In particular,
by part {\bf c)} of the lemma,
$\mathfrak{E}^+(s_0,s_0^\sigma)=\mathfrak{E}^+(s_0^\ddagger,s_0^\sigma)=1$.
\end{proof}

We denote 
\begin{equation}\label{muexplicit}
\mu^+(s)=C_0\frac{G_{\sigma\tau}(s)N_{\sigma}^+(s^{-1})}
{G_{\sigma\tau}(s_0)N_{\sigma}^+(s_0^{\ddagger})},\qquad
s\in\mathcal{S}^+,
\end{equation}
with $C_0=C_0^\alpha\in\mathbb{C}$ the normalization constant 
\eqref{C0plus}.
Note that the normalization constant is chosen in such a way that
\[\sum_{s\in\mathcal{S}^+}\mu^+(s)=\sum_{s\in\mathcal{S}}\mu(s)=
\frac{1}{G_{\tau\sigma\tau}(s_0^\tau)G(s_0^{\tau\sigma})},
\]
see \eqref{sumplus} for the first equality and \eqref{summutau} for
the second equality.
Applying the results of Proposition
\ref{symmpol} now
leads to the following result.

\begin{prop}\label{Eexpansionplus}
The symmetric Cherednik kernel $\mathfrak{E}^+$
associated to $\alpha$ is given by the series expansion
\[\mathfrak{E}^+(\gamma,x)=G_{\sigma\tau}(\gamma)G_\tau(x)
\sum_{s\in\mathcal{S}_\tau^+}\mu_\tau^+(s)E_\tau^+(s;x)
E_{\sigma\tau}^+(s;\gamma),
\]
with the series converging absolutely and uniformly
on compacta of $(\mathbb{C}^\times)^n\times (\mathbb{C}^\times)^n$.
\end{prop}
\begin{proof} 
Let the symmetrizer $C_+$ respectively $C_+^\sigma$ 
act on the $x$-variables respectively
the $\gamma$-variables within the series expansion 
\eqref{Eexpansions}
of $\mathfrak{E}(\gamma,x)$,
and apply \eqref{stableC} and Proposition \ref{symmpol}. We obtain
\[\mathfrak{E}^+(\gamma,x)=
G_{\sigma\tau}(\gamma)G_\tau(x)\sum_{\lambda\in\Lambda^+}
\Bigl(\sum_{\mu\in W_0\cdot\lambda}\mu_\tau(s_\mu^\tau)
\Bigr)E_{\tau}^+(s_\lambda^\tau;x)
E_{\sigma\tau}^+(s_\lambda^{\sigma\tau};\gamma).
\]
By \eqref{Ndecomposition}, \eqref{symmistrivialdiscrete} and
\eqref{WnulGaussian}, and by the explicit expressions \eqref{mutau}
and \eqref{muexplicit} for $\mu_\tau$ and $\mu_\tau^+$ respectively, 
we obtain
\[\sum_{\mu\in W_0\cdot\lambda}\mu_\tau(s_\mu^\tau)=
\mu_\tau^+(s_\lambda^\tau),\qquad \forall\,\lambda\in\Lambda^+,\]
from which the explicit series expansion for the symmetric Cherednik
function now immediately follows.
The convergence properties of the series are clear
from e.g. (the proof of) 
Corollary \ref{convcor} and Remark \ref{symmbound}.
\end{proof}

Define the spherical sub-algebra $\mathcal{H}^+=\mathcal{H}_\alpha^+$
of the double affine Hecke algebra $\mathcal{H}$ by
\[\mathcal{H}^+=C_+\mathcal{H}C_+=
\{X\in \mathcal{H} \, | \, C_+X=X=XC_+\}.
\]
Observe that $X\in\mathcal{H}$ is an element of
$\mathcal{H}^+$ if and only if $T_iX=t_iX=XT_i$ for all
$i=1,\ldots,n$. In particular, any element $X\in\mathcal{H}^+$,
considered as endomorphism of $\mathcal{M}$, 
maps into the field $\mathcal{M}_+$ of $W_0$-invariant
meromorphic functions, and factorizes through
the projection $C_+:\mathcal{M}\rightarrow \mathcal{M}_+$.
Note furthermore that the duality anti-isomorphism
$\psi$ restricts to an anti-isomorphism
$\psi: \mathcal{H}^+\rightarrow \mathcal{H}_\sigma^+$
since $\psi(C_+)=C_+^\sigma$.

Symmetrizing the main properties
of the normalized Cherednik kernels $\mathfrak{E}$ and
$\mathfrak{E}_\ddagger$ 
leads to the following
main result of this subsection.

\begin{thm}\label{symmmain} 
{\bf a)} $\mathfrak{E}^+(\gamma,x)=
\bigl(C_+^\ddagger\mathfrak{E}_\ddagger(\gamma,\cdot)\bigr)(x)$
\textup{(}inversion-invariance\textup{)}.\\
{\bf b)} For $X\in\mathcal{H}^+$, 
$\bigl(X\mathfrak{E}^+(\gamma,\cdot)\bigr)(x)=
\bigl(\psi(X)\mathfrak{E}^+(\cdot,x)\bigr)(\gamma)$ 
\textup{(}transformation behaviour\textup{)}.\\
{\bf c)} $\mathfrak{E}^+(\gamma,x)=\mathfrak{E}_\sigma^+(x,\gamma)$
\textup{(}duality\textup{)}.\\
{\bf d)} $\mathfrak{E}^+(s;x)=E^+(s;x)$ for $s\in\mathcal{S}^+$
\textup{(}polynomial reduction\textup{)}.
\end{thm}
\begin{proof}
{\bf a)} Similarly as for the normalized Cherednik kernel 
$\mathfrak{E}$ (see Lemma \ref{elementary}),
we have $\bigl(C_+^\ddagger\mathfrak{E}_\ddagger(\gamma,\cdot)\bigr)(x)=
\bigl(C_+^{\ddagger\sigma}\mathfrak{E}_\ddagger(\cdot,x)\bigr)(\gamma)$.
Using Proposition
\ref{symmpol}, the series expansion \eqref{Eexpansions}
for $\mathfrak{E}_\ddagger$, \eqref{Gaussianinversion}
and \eqref{stableC}, we can repeat now the arguments of the
proof of Proposition \ref{Eexpansionplus} to show that
\begin{equation*}
\begin{split}
\bigl(C_+^\ddagger\mathfrak{E}_\ddagger(\gamma,\cdot)\bigr)(x)&=
G_{\sigma\tau}(\gamma)G_\tau(x)\sum_{s\in\mathcal{S}_\tau^+}
\mu_\tau^+(s)E_{\tau}^+(s;x)E_{\sigma\tau}^+(s;\gamma)\\
&=\mathfrak{E}^+(\gamma,x).
\end{split}
\end{equation*}
{\bf b)} For $X\in \mathcal{H}^+$ we compute, using the transformation
behaviour of the normalized Cherednik kernel $\mathfrak{E}$ under the action of
$\mathcal{H}$ (see Theorem \ref{Cherednikkernels}),
\begin{equation*}
\begin{split}
\bigl(X\mathfrak{E}^+(\gamma,\cdot)\bigr)(x)&=
\bigl(X\mathfrak{E}(\gamma,\cdot)\bigr)(x)\\
&=\bigl(\psi(X)\mathfrak{E}(\cdot,x)\bigr)(\gamma)\\
&=\bigl(\psi(X)\mathfrak{E}^+(\cdot,x)\bigr)(\gamma),
\end{split}
\end{equation*}
where we used Lemma \ref{elementary}{\bf b)} and 
$\psi(X)\in\mathcal{H}_\sigma^+$ in the last equality.\\
{\bf c)} This follows from the duality of the normalized Cherednik kernel
$\mathfrak{E}$ (see Theorem \ref{dualitytheorem})
and Lemma \ref{elementary}{\bf b)}, 
\[\mathfrak{E}^+(\gamma,x)=\bigl(C_+\mathfrak{E}(\gamma,\cdot)\bigr)(x)
=\bigl(C_+\mathfrak{E}_\sigma(\cdot,\gamma)\bigr)(x)=
\mathfrak{E}_\sigma^+(x,\gamma).
\]
{\bf d)} This follows from the polynomial reduction of the
normalized Cherednik kernel (see Theorem \ref{polreduc}) and Lemma
\ref{elementary}{\bf c)},
\begin{equation*}
\begin{split}
\mathfrak{E}^+(s,\cdot)&=\mathfrak{E}^+(s^{-1},\cdot)\\
&=C_+\mathfrak{E}(s^{-1},\cdot)\\
&=C_+E(s;\cdot)=E^+(s;\cdot)
\end{split}
\end{equation*}
for all $s\in\mathcal{S}^+$.
\end{proof}

\begin{rem}
Note that $\mathcal{A}_+\subset \mathcal{H}^+$
(considered as multiplication operators), as well as
$\mathcal{A}_+(Y)\subset \mathcal{H}^+$, where
$\mathcal{A}_+(Y)$ is the sub-algebra of elements
$p(Y)$ ($p\in\mathcal{A}_+$) in $\mathcal{H}$ (e.g. by applying
the duality anti-isomorphism $\psi_\sigma$ to the inclusion
$\mathcal{A}_+\subset \mathcal{H}_\sigma^+$). In view of
Theorem \ref{symmmain}{\bf b)} we thus conclude that 
\begin{equation}\label{spectralsymmetric}
p(Y)\mathfrak{E}^+(\gamma,\cdot)=p(\gamma)\mathfrak{E}^+(\gamma,\cdot),
\qquad p(Y_\sigma)\mathfrak{E}^+(\cdot,x)=p(x)\mathfrak{E}^+(\cdot,x)
\end{equation}
for all $p\in\mathcal{A}_+$. The algebra of commuting $q$-difference
reflection operators $\mathcal{A}_+(Y)$, considered as endomorphisms
of $\mathcal{M}_+$, can be identified with an algebra
generated by $n$ algebraically independent, commuting $q$-difference
operators (see van Diejen \cite{vD0} for an explicit description of
these $q$-difference operators). 
One of these $q$-difference operators may be taken to be
Koornwinder's \cite{K} multivariable
extension of the Askey-Wilson second order $q$-difference
operator, see \cite{N} and Subsection 4.5.
The formulas \eqref{spectralsymmetric}
thus show that the symmetric Cherednik function 
$\mathfrak{E}^+(\gamma,\cdot)\in\mathcal{M}_+$ is a common
eigenfunction for these $q$-difference operators. 
\end{rem}

We note that Corollary \ref{consequenceab} in symmetrized
form reads as follows.
\begin{cor}\label{consequenceabsym}
{\bf a)} The expansion of the inverse Gaussian
$G_\tau^{-1}\in\mathcal{O}$ in terms of symmetric Macdonald-Koornwinder 
polynomials is given by
\[
G_\tau(x)^{-1}=
G_{\sigma\tau}(s_0)\sum_{s\in\mathcal{S}_\tau^+}\mu_\tau^+(s)E_\tau^+(s;x).
\]
Here the series converge absolutely and uniformly on compacta
of $(\mathbb{C}^\times)^n$.\\
{\bf b)} For all $s\in\mathcal{S}^+=\mathcal{S}_{\sigma\tau\sigma}^+$, we
have the identity
\[
E_{\sigma\tau\sigma}^+(s;Y^\tau)\bigl(G_\tau^{-1}\bigr)=
\frac{G_{\sigma\tau}(s_0)}{G_{\sigma\tau}(s)}E^+(s;\cdot)G_\tau^{-1}
\]
in $\mathcal{O}$.
\end{cor}
\begin{proof}
The proof is completely analogous to the proof of Corollary
\ref{consequenceab}.
\end{proof}
The explicit relations of the symmetric Macdonald-Koornwinder
polynomials under permutations of the parameters $\{a,b,c,d\}$
(see Proposition \ref{symmetrypar}) can be lifted to the symmetric
Cherednik kernel. We give one explicit example. 
\begin{prop}\label{Cherednikpar}
With the notations as in Proposition \ref{symmetrypar},
\begin{equation*}
\mathfrak{E}_\beta^+(\gamma,x)=
\left(\prod_{i=1}^n
\frac{\bigl(act^{2(n-i)}, qt^{2(i-n)}/bd;q\bigr)_{\infty}}
{\bigl(bct^{2(n-i)},qt^{2(i-n)}/ad;
q\bigr)_{\infty}}\right)\frac{G_{\beta_{\sigma\tau}}(\gamma)}
{G_{\alpha_{\sigma\tau}}(\gamma)}
\mathfrak{E}_\alpha^+(\gamma,x).
\end{equation*}
\end{prop}
\begin{proof}
We use the explicit expansion of $\mathfrak{E}^+$ in terms of 
symmetric Macdonald-Koornwinder polynomials (see Proposition 
\ref{Eexpansionplus}) to prove the proposition. We consider what
happens to each term in this expansion sum when replacing the
difference multiplicity function $\alpha$
by $\beta$. 

Observe that the spectral points 
$s_\lambda^\tau\in\mathcal{S}_\tau^+$ are invariant under replacement
of $\alpha$ by $\beta$.
By Proposition \ref{symmetrypar} the Macdonald-Koornwinder polynomials
in the expansion sum transform as
\begin{equation*}
\begin{split}
E_{\beta_{\sigma\tau}}^+(s_\lambda^\tau;\gamma)&=
E_{\alpha_{\sigma\tau}}^+(s_\lambda^\tau;\gamma),\\
E_{\beta_{\tau}}^+(s_\lambda^\tau;x)&=
\left(\prod_{i=1}^n
\frac{\bigl(act^{2(n-i)},qat^{2(n-i)}/d;q\bigr)_{\lambda_i}}
{\bigl(bct^{2(n-i)}, qbt^{2(n-i)}/d;q\bigr)_{\lambda_i}}
\left(\frac{b}{a}\right)^{\lambda_i}\right)E_{\alpha_\tau}(s_\lambda^\tau;x)
\end{split}
\end{equation*}
for $\lambda\in\Lambda^+$.
The weights $\mu_\tau^+(s_\lambda^\tau)$ ($\lambda\in\Lambda^+$)
can be expressed in terms of $q$-shifted factorials
in view of the formulas  \eqref{muexplicit}, \eqref{C0plus},
\eqref{Nplus} and \eqref{Gaussform},
compare with 
the formula \eqref{lhs} in the proof of Proposition \ref{evaluationprop}.
{}From this it follows by a straightforward computation that
\[\frac{\mu_{\beta_\tau}^+(s_\lambda^\tau)}
{\mu_{\alpha_\tau}^+(s_\lambda^\tau)}
=\prod_{i=1}^n\frac{\bigl(act^{2(n-i)}, qt^{2(i-n)}/bd;q\bigr)_{\infty}
\bigl(bct^{2(n-i)},qbt^{2(n-i)}/d;q\bigr)_{\lambda_i}}
{\bigl(bct^{2(n-i)},qt^{2(i-n)}/ad;q\bigr)_{\infty}
\bigl(act^{2(n-i)},qat^{2(n-i)}/d;q\bigr)_{\lambda_i}}
\left(\frac{a}{b}\right)^{\lambda_i}.
\] 
Since $G_{\beta_\tau}(x)=G_{\alpha_\tau}(x)$, the proposition now follows
by a direct computation using the explicit
series expansion for $\mathfrak{E}^+$ (see Proposition \ref{Eexpansionplus})
\end{proof}

\begin{rem}
The behaviour of the symmetric Mac\-do\-nald-Koorn\-win\-der 
polynomials and of the symmetric
Cherednik kernel under permutations of the Askey-Wilson parameters 
$\{a,b,c,d\}$ (see
Proposition \ref{symmetrypar} and Proposition \ref{Cherednikpar},
respectively) can be extended to the nonsymmetric level. For the
Macdonald-Koornwinder polynomials, one now uses 
the evaluation formula for nonsymmetric
Macdonald-Koornwinder polynomials as proven in \cite[Thm. 9.3]{St1}.
We do not give the formulas explicitly, since we do not need 
them in the present paper.
\end{rem}

We end this subsection by relating the one variable set-up
($n=1$) to the theory of basic hypergeometric series. Note first that
the roots of medium length in $R_{nr}$ are no longer present for rank
one. In particular, the parameter $t$ in the difference multiplicity
function $\alpha$ disappears. Furthermore, the rank one double
affine Hecke algebra involves two difference reflection operators, 
namely $T_0$ and $T_n=T_1$. The associated 
$Y$-operator is $Y=Y_1=T_1T_0$, see e.g.
\cite{NS} for details. 

In rank one, the action of $Y+Y^{-1}\in\mathcal{A}_+(Y)
\subset \mathcal{H}^+$
on $W_0$-invariant meromorphic
functions on $\mathbb{C}^\times$
(i.e. meromorphic functions $f$ satisfying $f(x^{-1})=f(x)$),
essentially 
coincides with the action of
the standard Askey-Wilson second-order $q$-difference operator
involving four parameters $\{a,b,c,d\}$, which are related to 
the difference multiplicity function $\alpha$ via formula \eqref{AWpar},
see e.g. \cite[Prop. 5.8]{NS}.
In particular, $\mathfrak{E}^+(\gamma,\cdot)\in\mathcal{M}_+$
is a meromorphic eigenfunction of the Askey-Wilson second
order $q$-difference operator, which admits an explicit
series expansion in terms of the Askey-Wilson polynomials, see
Proposition \ref{Eexpansionplus}. On the other hand, a basis
of eigenfunctions for the Askey-Wilson second order $q$-difference
operator was given explicitly by Ismail and Rahman \cite{IR}
in terms of very-well-poised ${}_8\phi_7$ basic hypergeometric
series, see 
also Suslov \cite{Sus}. Here the very-well-poised ${}_8\phi_7$
basic hypergeometric series, denoted by ${}_8W_7$, is defined by
\[
{}_8W_7\bigl(a;b,c,d,e,f;q,z\bigr)=
\sum_{k=0}^{\infty}
\frac{1-aq^{2k}}{1-a}
\frac{\bigl(a,b,c,d,e,f;q\bigr)_kz^k}
{\bigl(q,qa/b,qa/c,qa/d,qa/e,qa/f;q\bigr)_k}, 
\] 
see \cite{GR} for details. To relate these solutions with $\mathfrak{E}^+$, 
we thus need to evaluate the explicit
series expansion for $\mathfrak{E}^+$ (see Proposition 
\ref{Eexpansionplus}) in terms of ${}_8\phi_7$ basic hypergeometric series.
This was done in \cite{St3}. 
The result is as follows, see \cite[Thm. 4.2]{St3} 
and use \eqref{AWpolynomial}.
\begin{thm}\label{linkhyp}
In the one variable set-up \textup{(}$n=1$\textup{)},
the symmetric Cherednik kernel $\mathfrak{E}^+$ can be written as
\begin{equation*}
\begin{split}
\mathfrak{E}^+(\gamma,x)=
&\frac{\bigl(qax\gamma/\tilde{d},
qa\gamma/\tilde{d}x,q/ad,qa/d;q\bigr)_{\infty}}
{\bigl(\tilde{a}\tilde{b}\tilde{c}\gamma,q\gamma/\tilde{d},
qx/d,q/dx;q\bigr)_{\infty}}\\
&\qquad\qquad\times
{}_8W_7\bigl(\tilde{a}\tilde{b}\tilde{c}\gamma/q;
ax,a/x,\tilde{a}\gamma,\tilde{b}\gamma,\tilde{c}\gamma;
q,q/\tilde{d}\gamma\bigr)
\end{split}
\end{equation*}
when $|q/\tilde{d}\gamma|<1$, where we have denoted 
$\{a,b,c,d\}$ and
$\{\tilde{a},\tilde{b},\tilde{c},\tilde{d}\}$
for the Askey-Wilson parameters \eqref{AWpar} associated to
$\alpha$ and $\alpha_\sigma$, respectively.
\end{thm}
\begin{rem}\label{linkhypremark}
{\bf a)}
The meromorphic continuation for the expression of $\mathfrak{E}^+$
as a ${}_8W_7$ series can be written explicitly as a sum
of two balanced ${}_4\phi_3$'s using Bailey's formula 
\cite[(2.10.10)]{GR}, see \cite[(3.2)]{St3}.

{\bf b)} Theorem \ref{linkhyp} shows that the Cherednik kernel 
$\mathfrak{E}^+$ in rank one ($n=1$) is the special eigenfunction
of the Askey-Wilson second-order $q$-difference operator
named the {\it Askey-Wilson function} in \cite{KS2} and \cite{St3}
(which was defined in these papers in terms of the
${}_8W_7$ series).
In \cite{KS1} the Askey-Wilson function was interpreted as
spherical function on the noncompact quantum group 
$\hbox{SU}_q(1,1)$. In particular, the kernel $\mathfrak{E}^+$ in rank one
may be regarded as the natural $q$-analogue of the Jacobi function
(see also \cite{KS3}).
\end{rem}


\section{An extension of the Macdonald-Koornwinder transform}

In this section we study a difference Fourier transform which is closely
related to the Macdonald-Koornwinder transform (see Section 4). 
In fact, one essentially 
replaces in the Macdonald-Koornwinder transform $F_{\mathcal{A}}$
the kernel $\mathfrak{E}_{\mathcal{A},\ddagger}$ by the normalized Cherednik
kernel $\mathfrak{E}_\ddagger$, and the cyclic module $\mathcal{A}$ by the
cyclic $\mathcal{H}$-module $\mathcal{A}G^{-1}$. We follow closely
the general line of arguments for difference Fourier transforms
as explained in Section 3.
The results in this section generalize results of Cherednik \cite{C1},
\cite{Cnew} to the nonreduced setup, as well as results in \cite{St3} 
to the multivariable setup.

Throughout this section we keep the same assumptions 
on the difference multiplicity function
$\alpha=(\mathbf{t},q^{\frac{1}{2}})=(t_0,u_0,t_n,u_n,t,q^{\frac{1}{2}})$
as in Section 5.


\subsection{The transform}

Recall the definition of the contour $\mathcal{T}=\mathcal{T}_\alpha$
in Subsection 4.3. We fix in this section such a contour $\mathcal{T}$,
satisfying the additional requirement that 
the parameter $d/q=-q^{-\frac{1}{2}}t_0u_0^{-1}$ 
is in the exterior of $\mathcal{T}$. It is convenient to fix the
contour $\mathcal{T}$ once and for all, although the main definitions in this
section are easily seen to be independent of this choice. 

The extra assumption on the fixed contour
$\mathcal{T}=\mathcal{T}_\alpha$ implies that the pairings
$\bigl(\cdot,\cdot\bigr)_{\mathcal{A}}=
\bigl(\cdot,\cdot\bigr)_{\mathcal{A},\alpha}$ and 
$\bigl(\cdot,\cdot\bigr)_{\mathcal{A},\tau}$ on $\mathcal{A}$
can be written as integrals over the same deformed torus $\mathcal{T}^n$,
\begin{equation}\label{formulaspairing}
\begin{split}
\bigl(p,r\bigr)_{\mathcal{A}}&=\frac{1}{(2\pi i)^n}
\underset{\mathcal{T}^n}{\iint}p(x)r(x^{-1})\Delta(x)\frac{dx}{x},\\
\bigl(p,r\bigr)_{\mathcal{A},\tau}&=\frac{1}{(2\pi i)^n}
\underset{\mathcal{T}^n}{\iint}p(x)r(x^{-1})\Delta_\tau(x)\frac{dx}{x}
\end{split}
\end{equation}
for $p,r\in\mathcal{A}$. We furthermore use the 
formulas \eqref{formulaspairing} as the definitions of   
$\bigl(p,r)_{\mathcal{A}}$ and $\bigl(p,r\bigr)_{\mathcal{A},\tau}$
for those meromorphic functions $p$ and $r$ which are regular on 
$\mathcal{T}^n$.

Consider now the subspace $V=V_\alpha=\mathcal{A}G^{-1}\subset
\mathcal{O}$. For all $X\in \mathcal{H}=\mathcal{H}_\alpha$ we have
\[ X\bigl(pG^{-1}\bigr)=\bigl(\tau(X)p\bigr)G^{-1},\qquad
p\in\mathcal{A},
\]
hence $\mathcal{A}G^{-1}$ is a cyclic $\mathcal{H}$-module, with cyclic vector
$G^{-1}$. Let $\mathfrak{E}_\ddagger$ be the normalized Cherednik kernel
associated to $\alpha$. Since 
$G_{\sigma\tau}(\gamma)^{-1}G_\tau(x)^{-1}\mathfrak{E}_\ddagger(\gamma,x)$
is analytic at $(\gamma,x)\in (\mathbb{C}^\times)^n\times 
(\mathbb{C}^\times)^n$ and $G_\tau\in\mathcal{M}$ is regular on
$\mathcal{T}^n$, we may define the transform 
$F=F_\alpha: \mathcal{A}G^{-1}\rightarrow
\mathcal{O}G_{\sigma\tau}$ by
\begin{equation}
\bigl(Fg\bigr)(\gamma)=\bigl(g,\mathfrak{E}_\ddagger(\gamma^{-1},
\cdot)\bigr)_{\mathcal{A}},\qquad g\in \mathcal{A}G^{-1}.
\end{equation}
We collect some elementary properties of $F$.
\begin{lem}
{\bf a)} For all $p\in\mathcal{A}$,
\begin{equation}\label{analyticform}
F\bigl(pG^{-1}\bigr)(\gamma)=\bigl(p,\mathfrak{E}_{\ddagger}(\gamma^{-1},
\cdot)G_{\tau}^{-1}\bigr)_{\mathcal{A},\tau}.
\end{equation}
{\bf b)} $F:\mathcal{A}G^{-1}\rightarrow \mathcal{M}$ 
is a Fourier transform associated with $\sigma$.
\end{lem}
\begin{proof}
{\bf a)} By a direct computation using the 
explicit expression of the weight function
$\Delta$ (see \eqref{decomposition} and \eqref{Deltaform})
and of the Gaussian $G$, we have
\begin{equation}\label{relationweights}
G(x)^{-1}G_\tau(x)\Delta(x)=\Delta_\tau(x).
\end{equation}
The claim is now immediate in view of the special choice of contour
$\mathcal{T}$.\\
{\bf b)} The advantage of the expression \eqref{analyticform} is that the
kernel  
\[
\mathfrak{E}_\ddagger(\gamma^{-1},\cdot)G_\tau^{-1}=
\mathfrak{E}_\ddagger(\gamma^{-1},\cdot)G_{\ddagger}
\]
is analytic on $(\mathbb{C}^\times)^n$ for generic 
$\gamma\in (\mathbb{C}^\times)^n$, hence
the adjoint of $X\in\mathcal{H}_\tau$ 
acting on $p$ in the pairing \eqref{analyticform} is  
$\ddagger_\tau(X)$ (see Remark \ref{adjoint1remark}). 
We now compute for $X\in\mathcal{H}$,
\begin{equation*}
\begin{split}
F\bigl(X(pG^{-1})\bigr)(\gamma)&=F\bigl((\tau(X)p)G^{-1}\bigr)(\gamma)\\
&=\bigl(\tau(X)p,
\mathfrak{E}_{\ddagger}(\gamma^{-1},\cdot)G_{\ddagger}
\bigr)_{\mathcal{A},\tau}\\
&=\bigl(p,(\ddagger_\tau\circ\tau)(X)
\bigl(\mathfrak{E}_{\ddagger}(\gamma^{-1},\cdot)G_{\ddagger}
\bigr)\bigr)_{\mathcal{A},\tau}\\
&=\bigl(Z\bigl(F(pG^{-1})\bigr)\bigr)(\gamma),
\end{split}
\end{equation*}
with $Z\in\mathcal{H}_\sigma$ given by
\[Z=(\dagger_{\ddagger\sigma}\circ\psi_\ddagger\circ\tau_{\ddagger}^{-1}
\circ\ddagger_\tau\circ\tau)(X)=\sigma(X).
\]
Here the last equality follows by computing the left hand side and
the right hand side explicitly on a set of 
algebraic generators of $\mathcal{H}$.
\end{proof}

Consider the sub-space $W_\sigma=W_{\alpha_\sigma}=
\mathcal{A}G_{\sigma\tau}\subset \mathcal{M}$. For all
$X\in\mathcal{H}_\sigma$ we have
\[X(pG_{\sigma\tau})=(\tau_{\sigma\tau}^{-1}(X)p)G_{\sigma\tau},\qquad
p\in\mathcal{A},\]
hence $\mathcal{A}G_{\sigma\tau}$ is a cyclic
$\mathcal{H}_{\sigma}$-module with cyclic vector $G_{\sigma\tau}$.
The expansion formula \eqref{Eexpansions} 
for the normalized Cherednik kernel $\mathfrak{E}_\ddagger$
leads now to the following result. 
\begin{prop}\label{Fexplicit}
The difference Fourier transform $F$ 
defines a linear bijection $F: \mathcal{A}G^{-1}\rightarrow 
\mathcal{A}G_{\sigma\tau}$. Explicitly,
we have for all $s\in \mathcal{S}_\tau=\mathcal{S}_{\sigma\tau}$,
\begin{equation}\label{formulades}
F\bigl(E_\tau(s;\cdot)G^{-1}\bigr)(\gamma)=
D_0\,G_{\tau\sigma\tau}(s)E_{\sigma\tau}(s;\gamma)G_{\sigma\tau}(\gamma),
\end{equation}
with $D_0$ the constant
\begin{equation}\label{D0}
D_0=\mathcal{C}(s_0^{\ddagger\sigma})
\prod_{i=1}^n\frac{\bigl(t^2,bct^{2(n-i)},dt^{2(i-n)}/a,qt^{2(i-n)}/ad;
q\bigr)_{\infty}}{\bigl(q,t^{2(n-i+1)},abt^{2(n-i)},act^{2(n-i)};
q\bigr)_{\infty}},
\end{equation}
where we used the Askey-Wilson parametrization \eqref{AWpar} for part
of the difference multiplicity function $\alpha$.
\end{prop}
\begin{proof}
We use the series expansion \eqref{Eexpansions} for the kernel
$\mathfrak{E}_\ddagger$
together with \eqref{analyticform} and the orthogonality relations 
\eqref{ortho} for the Macdonald-Koornwinder polynomials. Then we obtain for
$s\in\mathcal{S}_\tau$,
\begin{equation*}
\begin{split}
F\bigl(E_\tau(s;\cdot)G^{-1}\bigr)(\gamma)&=
G_{\sigma\tau}(\gamma)\sum_{v\in\mathcal{S}_\tau}\mu_\tau(v)
\bigl(E_{\tau}(s;\cdot),
E_{\ddagger\tau}(v^{-1};\cdot)\bigr)_{\mathcal{A},\tau}
E_{\sigma\tau}(v;\gamma)\\
&=\mu_\tau(s)\bigl(E_\tau(s;\cdot),
E_{\ddagger\tau}(s^{-1};\cdot)\bigr)_{\mathcal{A},\tau}
E_{\sigma\tau}(s;\gamma)G_{\sigma\tau}(\gamma).
\end{split}
\end{equation*}
By the explicit expression \eqref{mutau} for $\mu_\tau$
and by Theorem \ref{corpoly}{\bf b)}, this simplifies to
\[F\bigl(E_\tau(s;\cdot)G^{-1}\bigr)(\gamma)= 
D_0\,G_{\tau\sigma\tau}(s)E_{\sigma\tau}(s;\gamma)G_{\sigma\tau}(\gamma),
\]
with the constant $D_0$ given by
\[D_0=\frac{C_0^\tau\bigl(1,1\bigr)_{\mathcal{A},\tau}}
{G_{\tau\sigma\tau}(s_0^\tau)}.
\]
Now $D_0$ can be evaluated explicitly using \eqref{symmcontinuous},
\eqref{Gustafson}, Remark \ref{Gaussianevaluation} and \eqref{C0plus}.
Using furthermore the fact that
$\mathcal{C}_{\tau}(s_0^{\ddagger\tau\sigma})=
\mathcal{C}(s_0^{\ddagger\sigma})$, we obtain 
the explicit expression \eqref{D0} for the constant $D_0$.
\end{proof}


\subsection{The inverse transform}

The next step is to invert the difference Fourier transform
$F: \mathcal{A}G^{-1}\rightarrow \mathcal{A}G_{\sigma\tau}$. 
We define a transform $J_\sigma=J_{\alpha_\sigma}$ by the formula
\begin{equation}\label{inversedef}
\bigl(J_{\sigma}g\bigr)(x):=\lbrack
g,\mathfrak{E}(\cdot,x)\rbrack_{\mathcal{A},\sigma}
=\sum_{s\in \mathcal{S}_\ddagger}g(s)\mathfrak{E}(s,x)N_\sigma(s)
\end{equation}
for $g\in\mathcal{A}G_{\sigma\tau}$. The defining sum
\eqref{inversedef} converges absolutely and uniformly for $x$ in
compacta of $(\mathbb{C}^\times)^n$. This follows from the alternative
expression 
\begin{equation}\label{alternativeJ}
\bigl(J_{\sigma}g\bigr)(x)=\lbrack
g,\mathfrak{E}_{\mathcal{A}}(\cdot,x)\rbrack_{\mathcal{A},\sigma},\qquad
g\in \mathcal{A}G_{\sigma\tau}
\end{equation}
(see Theorem \ref{polreduc}), combined with 
the bounds for the Gaussian and
for the Mac\-do\-nald-Koorn\-win\-der polynomials 
$\mathfrak{E}_{\mathcal{A}}(s,x)=E(s^{-1};x)$ 
($s\in\mathcal{S}_\ddagger$), see \eqref{Gaussianbound} and
Proposition \ref{MKbounds} respectively.
In particular, $J_{\sigma}$ defines a linear map $J_\sigma:
\mathcal{A}G_{\sigma\tau}\rightarrow \mathcal{O}$. 
\begin{lem}
$J_\sigma: \mathcal{A}G_{\sigma\tau}\rightarrow \mathcal{O}$ 
is a Fourier transform associated with $\sigma^{-1}$.
\end{lem}
\begin{proof}
This follows from
\eqref{alternativeJ}, Lemma \ref{actioncompatibleH}, 
Remark \ref{adjoint2remark}, Proposition
\ref{allrelations} and 
the fact that $\sigma^{-1}=\psi_\sigma\circ\iota_\sigma$ as unital
algebra isomorphisms from $\mathcal{H}_\sigma$ to $\mathcal{H}$.
\end{proof}

\begin{prop}\label{Jexplicit}
The difference Fourier transform $J_\sigma$ defines a
bijection $J_\sigma: \mathcal{A}G_{\sigma\tau}\rightarrow
\mathcal{A}G^{-1}$. Explicitly, we have for
$s\in\mathcal{S}_\tau=\mathcal{S}_{\sigma\tau}$,
\[J_\sigma\bigl(E_{\sigma\tau}(s;\cdot)G_{\sigma\tau}\bigr)(x)=
E_0\,G_{\tau\sigma\tau}(s)^{-1}E_\tau(s;x)G(x)^{-1},
\]
with the constant $E_0\in\mathbb{C}$ given by
\begin{equation}\label{E0}
\begin{split}
E_0=N_\sigma(s_0^\ddagger)
\prod_{i=1}^n&\left\{\frac{\bigl(abcdt^{2(2n-i-1)};q\bigr)_{\infty}}
{\bigl(adt^{2(n-i)},bdt^{2(n-i)},cdt^{2(n-i)};q\bigr)_{\infty}}\right.\\
&\left.\,\,\,\times\frac{1}{\bigl(bct^{2(n-i)},bct^{2(n-i)},
dt^{2(i-n)}/a,qt^{2(i-n)}/ad;q\bigr)_{\infty}}\right\},
\end{split}
\end{equation}
where we used the Askey-Wilson parametrization \eqref{AWpar} for
part of the difference multiplicity function $\alpha$.
\end{prop}
\begin{proof}
In view of \eqref{alternativeJ},
\eqref{Eexpansions}, \eqref{polduality} and 
the definition \eqref{mutau} for the weight $\mu\in\mathcal{F}(\mathcal{S})$,
we easily see that
\begin{equation*}
\begin{split}
J_\sigma\bigl(E_{\sigma\tau}(s;\cdot)G_{\sigma\tau}\bigr)(x)&=
\frac{G_{\sigma\tau}(s_0)N_\sigma(s_0^\ddagger)}{C_0}\,
G_{\tau\sigma\tau}(s)^{-1}\mathfrak{E}_\tau(s^{-1},x)G(x)^{-1}\\
&=\frac{G_{\sigma\tau}(s_0)N_\sigma(s_0^\ddagger)}{C_0}\,
G_{\tau\sigma\tau}(s)^{-1}E_\tau(s;x)G(x)^{-1},\\
\end{split}
\end{equation*}
where we used the polynomial reduction in the second equality
(see Theorem \ref{polreduc}). By Remark \ref{Gaussianevaluation}
and \eqref{C0plus}, the coefficient 
$G_{\sigma\tau}(s_0)N_\sigma(s_0^\ddagger)/C_0$ is easily
seen to be equal to the constant $E_0$ as defined in \eqref{E0}.
\end{proof}
Recall the notation
$c_{\mathcal{A}}=\bigl(1,1\bigr)_{\mathcal{A}}N_\sigma(s_0^{-1})$
from Theorem \ref{corpoly}.
\begin{cor}\label{inversionextension}
The difference Fourier transform $F: \mathcal{A}G^{-1}\rightarrow 
\mathcal{A}G_{\sigma\tau}$ 
is a linear bijection, with inverse 
$c_{\mathcal{A}}^{-1}J_\sigma: \mathcal{A}G_{\sigma\tau}\rightarrow 
\mathcal{A}G^{-1}$.
\end{cor}
\begin{proof}
Proposition \ref{Fexplicit} and Proposition \ref{Jexplicit}
show that
$F\circ J_\sigma=D_0E_0\,\hbox{Id}$ on $\mathcal{A}G_{\sigma\tau}$ 
and $J_\sigma\circ F=D_0E_0\,\hbox{Id}$ on $\mathcal{A}G^{-1}$.
It thus suffices to show that
$D_0E_0=c_{\mathcal{A}}$.
Now by \eqref{Gustafson} and by
the explicit expressions \eqref{D0}, \eqref{E0}
for $D_0$ and $E_0$, it follows that
\begin{equation}\label{DE}
D_0E_0=N_\sigma(s_0^\ddagger)\frac{\mathcal{C}(s_0^{\ddagger\sigma})}{2^nn!}
\bigl(1,1\bigr)_{\mathcal{A},+}=
N_\sigma(s_0^\ddagger)\bigl(1,1\bigr)_{\mathcal{A}}
=c_{\mathcal{A}},
\end{equation}
where the second equality follows from \eqref{symmcontinuous}.
\end{proof}


\subsection{Plancherel-type formulas}
In this subsection we prove Plancherel-type formulas for the transform
$F$ and its inverse $c_{\mathcal{A}}^{-1}J_{\sigma}$.
For this, we introduce two new 
transforms $\widetilde{F}: \mathcal{A}G^{-1}\rightarrow
\mathcal{O}G_{\sigma\tau}$ and $\widetilde{J}_\sigma:
\mathcal{A}G_{\sigma\tau}\rightarrow \mathcal{O}$ by
\begin{equation}
\bigl(\widetilde{F}g\bigr)(\gamma)=
\bigl(\mathfrak{E}(\gamma,\cdot),g\bigr)_{\mathcal{A}},
\qquad
\bigl(\widetilde{J}_{\sigma}h\bigr)(x)=
\lbrack I\mathfrak{E}_\ddagger(\cdot,x),
h\rbrack_{\mathcal{A},\sigma}
\end{equation}
for $g\in \mathcal{A}G^{-1}$ and $h\in\mathcal{A}G_{\sigma\tau}$, with
$I$ the inversion operator $(Ig)(\gamma)=g(\gamma^{-1})$.
Repeating the arguments of the previous subsections lead to the
following result.
\begin{prop}\label{explicittilde}
 The transform $\widetilde{F}$ defines a linear bijection
$\widetilde{F}: \mathcal{A}G^{-1}\rightarrow
\mathcal{A}G_{\sigma\tau}$, whose inverse is given by
$c_{\mathcal{A}}^{-1}\widetilde{J}_\sigma: \mathcal{A}G_{\sigma\tau}
\rightarrow \mathcal{A}G^{-1}$. Explicitly, we have
\begin{equation*}
\begin{split}
\widetilde{F}\bigl(E_{\ddagger\tau}(s^{-1};\cdot)G^{-1}\bigr)(\gamma)&=
D_0\,G_{\tau\sigma\tau}(s)E_{\sigma\tau}(s;\gamma)G_{\sigma\tau}(\gamma),\\
\widetilde{J}_\sigma\bigl(E_{\sigma\tau}(s;\cdot)G_{\sigma\tau}\bigr)(x)&=
E_0\,G_{\tau\sigma\tau}(s)^{-1}E_{\ddagger\tau}(s^{-1};x)G(x)^{-1}
\end{split}
\end{equation*}
for $s\in\mathcal{S}_\tau=\mathcal{S}_{\sigma\tau}$, with the
constants $D_0\in\mathbb{C}$ and $E_0\in\mathbb{C}$ given by
\eqref{D0} and \eqref{E0}, respectively.
\end{prop}

The transforms $\widetilde{F}$ and $\widetilde{J}_\sigma$ are defined
in such a way that
\[
\lbrack Fg,h\rbrack_{\mathcal{A},\sigma}=
\bigl(g,\widetilde{J}_{\sigma}h\bigr)_{\mathcal{A}},
\qquad \bigl(J_{\sigma}h,g\bigr)_{\mathcal{A}}=\lbrack
h,\widetilde{F}g\rbrack_{\mathcal{A},\sigma}
\]
for all $g\in \mathcal{A}G^{-1}$ and $h\in\mathcal{A}G_{\sigma\tau}$.
This can be proven by interchanging integration and summation using
Fubini's Theorem, which is justified by the polynomial
reduction of the Cherednik kernels (Theorem \ref{polreduc})
and by the bounds for the
Macdonald-Koornwinder polynomials (Proposition \ref{MKbounds}).
This leads now immediately to the following Plancherel-type
formulas.
\begin{prop}\label{Plancherelnew}
{\bf a)} Let $g,h\in\mathcal{A}G^{-1}$, then
\[ \lbrack Fg, \widetilde{F}h\rbrack_{\mathcal{A},\sigma}=
c_{\mathcal{A}}\bigl(g,h\bigr)_{\mathcal{A}}.
\]
{\bf b)} Let $g,h\in\mathcal{A}G_{\sigma\tau}$, then
\[ \bigl(J_{\sigma}g,\widetilde{J}_{\sigma}h\bigr)_{\mathcal{A}}=
c_{\mathcal{A}}\lbrack g,h\rbrack_{\mathcal{A},\sigma}.
\]
\end{prop}
Combining Proposition \ref{Plancherelnew} with Proposition
\ref{Fexplicit}, Proposition
\ref{explicittilde} and \eqref{DE} leads to the following 
formulas involving Macdonald-Koornwinder polynomials.
\begin{cor}\label{formulasnew}
Let $s,v\in\mathcal{S}_\tau=\mathcal{S}_{\sigma\tau}$, then
\begin{equation*}
\begin{split}
D_0G_{\tau\sigma\tau}(s)G_{\tau\sigma\tau}(v)\,
&\lbrack E_{\sigma\tau}(s;\cdot)G_{\sigma\tau},
E_{\sigma\tau}(v;\cdot)G_{\sigma\tau}\rbrack_{\mathcal{A},\sigma}=\\
&\quad=E_0\bigl(E_\tau(s;\cdot)G^{-1},
E_{\ddagger\tau}(v^{-1};\cdot)G^{-1}\bigr)_{\mathcal{A}},
\end{split}
\end{equation*}
with $D_0$ and $E_0$ given by \eqref{D0} and \eqref{E0},
respectively.
\end{cor}
The explicit formulas 
of Corollary \ref{formulasnew} are the
analogues of the orthogonality relations \eqref{ortho}
and quadratic
norm evaluations (Theorem \ref{corpoly}{\bf b)})
for the Mac\-do\-nald-Koorn\-win\-der polynomials.

With the main results for the difference Fourier transform
$F:\mathcal{A}G^{-1}\rightarrow \mathcal{A}G_{\sigma\tau}$
now established, we can make the connection with the
Macdonald-Koornwinder transform $F_{\mathcal{A}}$ and its inverse
$c_{\mathcal{A}}^{-1}J_{\mathcal{A},\sigma}$ 
more explicit. Since the inverse transform $J_\sigma$ is a discrete
transform supported on the polynomial spectrum
$\mathcal{S}_\ddagger$, we may as well consider $F$ and $J_\sigma$ as maps
\begin{equation}\label{discreteview}
F_{res}: \mathcal{A}G^{-1}\rightarrow 
\bigl(\mathcal{A}G_{\sigma\tau}\bigr)|_{\mathcal{S}_\ddagger},
\qquad J_{res,\sigma}:
\bigl(\mathcal{A}G_{\sigma\tau}\bigr)|_{\mathcal{S}_\ddagger}
\rightarrow \mathcal{A}G^{-1}
\end{equation}
by restriction of the spectral variable $\gamma$ to $\mathcal{S}_\ddagger$.
Here the space
$\bigl(\mathcal{A}G_{\sigma\tau}\bigr)|_{\mathcal{S}_\ddagger}$
is again a (cyclic) $\mathcal{H}_\sigma$-submodule of
$\mathcal{F}(\mathcal{S}_\ddagger)$ by Lemma \ref{actioncompatibleH}.
The transforms $F_{res}$ and $J_{res,\sigma}$
are then Fourier transforms associated to
$\sigma$ and $\sigma^{-1}$ respectively, and 
$c_{\mathcal{A}}^{-1}J_{res,\sigma}$ is the inverse of the transform
$F_{res}$. Moreover, the transforms
$F_{res}$ and $J_{res,\sigma}$ can be expressed as
\begin{equation}\label{discreteviewexplicit}
\begin{split}
\bigl(F_{res}g\bigr)(s)&=
\bigl(g,\mathfrak{E}_{\mathcal{A},\ddagger}(s^{-1},\cdot)
\bigr)_{\mathcal{A}},\qquad s\in\mathcal{S}_\ddagger,\\
\bigl(J_{res,\sigma}h\bigr)(x)&=
\lbrack h,\mathfrak{E}_{\mathcal{A}}(\cdot,x)\rbrack_{\mathcal{A},\sigma}
\end{split}
\end{equation}
for $g\in \mathcal{A}G^{-1}$
and $h\in
\bigl(\mathcal{A}G_{\sigma\tau}\bigr)|_{\mathcal{S}_\ddagger}$,
which coincide with the defining formulas for the 
Macdonald-Koornwinder transform
$F_{\mathcal{A}}:\mathcal{A}\rightarrow
\mathcal{F}_0(\mathcal{S}_\ddagger)$ and the discrete transform
$J_{\mathcal{A},\sigma}: \mathcal{F}_0(\mathcal{S}_\ddagger)\rightarrow
\mathcal{A}$.

In particular, the Macdonald-Koornwinder transform and 
its extension to the cyclic $\mathcal{H}$-module 
$\mathcal{A}G^{-1}$ can be treated together in a uniform manner
by considering difference Fourier transforms on the 
$\mathcal{H}$-module $V=V_\alpha$ defined by
\[ V=\mathcal{A}\oplus \mathcal{A}G^{-1}\subset \mathcal{M},
\]
and on the $\mathcal{H}_\sigma$-submodule $W_\sigma=
W_{\alpha_\sigma}$ defined by
\[ W_\sigma=\mathcal{F}_0(\mathcal{S}_\ddagger)\oplus 
\bigl(\mathcal{A}G_{\sigma\tau}\bigr)_{|\mathcal{S}_\ddagger}\subset
\mathcal{F}(\mathcal{S}_\ddagger)
\]
(clearly, the sum is direct in both cases). The transforms
$F_{res}: V\rightarrow W_\sigma$ and $J_{res,\sigma}:
W_\sigma\rightarrow V$ are
then defined by \eqref{discreteviewexplicit}, now with
$g\in V$ and $h\in W_\sigma$. These extended transforms
$F_{res}: V\rightarrow W_\sigma$ and 
$J_{res,\sigma}: W_\sigma\rightarrow V$
are Fourier transforms associated to $\sigma$ and $\sigma^{-1}$
respectively, and $c_{\mathcal{A}}^{-1}J_{res,\sigma}$ is the inverse
of $F_{res}$. Furthermore, applying Fubini's Theorem we 
have the Plancherel-type formulas
\begin{equation}\label{newformPlancherel}
\begin{split}
\lbrack F_{res}g, \widetilde{F}_{res}h\rbrack_{\mathcal{A},\sigma}&=
c_{\mathcal{A}}\bigl(g,h\bigr)_{\mathcal{A}},\qquad g,h\in V,\\
\bigl(J_{res,\sigma}g, \widetilde{J}_{res,\sigma}h\bigr)_{\mathcal{A}}&=
c_{\mathcal{A}}\lbrack g,h\rbrack_{\mathcal{A},\sigma},\qquad g,h\in W_\sigma,
\end{split}
\end{equation}
with $\widetilde{F}_{res}: V\rightarrow W_\sigma$ and 
$\widetilde{J}_{res,\sigma}: W_\sigma\rightarrow V$ the transforms
\[
\bigl(\widetilde{F}_{res}g\bigr)(s)=\bigl(\mathfrak{E}_{\mathcal{A}}(s,\cdot),
g\bigr)_{\mathcal{A}},\qquad
\bigl(\widetilde{J}_{res,\sigma}h\bigr)(x)=\lbrack
I\mathfrak{E}_{\mathcal{A},\ddagger}(\cdot,x),h\rbrack_{\mathcal{A},\sigma}
\]
for $s\in\mathcal{S}_\ddagger$, $g\in V$ and $h\in W_\sigma$, where
$I$ is the inversion operator $(Ig)(v)=g(v^{-1})$.
By the explicit expressions for the
images of suitable
bases of $V$ and $W_\sigma$ under $F_{res}$, $\widetilde{F}_{res}$ and
$J_{res,\sigma}$, $\widetilde{J}_{res,\sigma}$ respectively
(see Section 4 as well as this section), the 
Plancherel-type formulas \eqref{newformPlancherel} lead to 
the orthogonality relations and quadratic norm evalutions
of the Macdonald-Koornwinder polynomials and to the 
formulas of Corollary \ref{formulasnew}.
It also leads to ``mixed identities'', for which $g$ and $h$
in \eqref{newformPlancherel} are taken from different summands
in $V$ (respectively $W_\sigma$). These mixed identities 
are completely covered by the following integral formulas for
Macdonald-Koornwinder polynomials.
\begin{prop}\label{mixed}
For $v\in\mathcal{S}=\mathcal{S}_{\sigma\tau\sigma}$ 
and $s\in\mathcal{S}_\tau$,
\begin{equation*}
\bigl(E(v;\cdot),E_{\ddagger\tau}(s^{-1};\cdot)G^{-1}
\bigr)_{\mathcal{A}}=
D_0G_{\sigma\tau}(v)G_{\tau\sigma\tau}(s)
E_{\sigma\tau}(s;v^{-1}).
\end{equation*}
\end{prop}
\begin{proof} 
Let $v\in\mathcal{S}$
and $s\in\mathcal{S}_\tau$. Let
$\delta_{v^{-1}}\in\mathcal{F}_0(\mathcal{S}_\ddagger)$ be the
function which is one at $v^{-1}\in\mathcal{S}_\ddagger$ and zero
otherwise. Then \eqref{ortho}
and Theorem \ref{corpoly}{\bf b)} imply that
\[ F_{res}(E(v;\cdot))=\bigl(E(v;\cdot),E_\ddagger(v^{-1};\cdot)
\bigr)_{\mathcal{A}}\delta_{v^{-1}}=
\frac{c_{\mathcal{A}}}{N_\sigma(v^{-1})}\,\delta_{v^{-1}}.
\]
On the other hand,
\[\widetilde{F}_{res}(E_{\ddagger\tau}(s^{-1};\cdot)G^{-1})=
D_0G_{\tau\sigma\tau}(s)E_{\sigma\tau}(s;\cdot)G_{\sigma\tau}.
\]
Combining these two formulas with the first Plancherel-type
formula in \eqref{newformPlancherel} leads to the desired identity.
\end{proof}


\subsection{The symmetric theory}
The results on the ex\-ten\-ded Mac\-do\-nald-Koorn\-win\-der transform
$F: \mathcal{A}G^{-1}\rightarrow \mathcal{A}G_{\sigma\tau}$
and its inverse $J_\sigma: \mathcal{A}G_{\sigma\tau}\rightarrow
\mathcal{A}G^{-1}$ can be symmetrized in the usual manner
by applying the symmetrizer $C_+\in H_0\subset \mathcal{H}$.
We collect the main formulas in this subsection.

We define the symmetric transforms $F^+=F_{\alpha}^+:
\mathcal{A}_+G^{-1}\rightarrow \mathcal{M}$ 
and $J^+_\sigma=J_{\alpha_\sigma}^+: \mathcal{A}_+G_{\sigma\tau}\rightarrow
\mathcal{M}$ by the formulas
\begin{equation}\label{Fplus}
\begin{split}
\bigl(F^+g\bigr)(\gamma)&=\bigl(g,\mathfrak{E}^+(\gamma,\cdot)
\bigr)_{\mathcal{A},+},\\
\bigl(J^+_{\sigma}h\bigr)(x)&=\lbrack
h,\mathfrak{E}^+(\cdot,x)\rbrack_{\mathcal{A},+,\sigma}=
\lbrack
h,\mathfrak{E}^+_\sigma(x,\cdot)\rbrack_{\mathcal{A},+,\sigma},
\end{split}
\end{equation}
for $g\in \mathcal{A}_+G^{-1}$ and $h\in\mathcal{A}_+G_{\sigma\tau}$,
where the integration for the pairing
$\bigl(\cdot,\cdot\bigr)_{\mathcal{A},+}$ (see \eqref{pairingsymm})
is over the deformed torus $\mathcal{T}^n$ with $\mathcal{T}$ as in
the previous subsection, and with
$\lbrack\cdot,\cdot\rbrack_{\mathcal{A},+}=
\lbrack \cdot,\cdot\rbrack_{\mathcal{A},+,\alpha}$ defined by
\[\lbrack f,g\rbrack_{\mathcal{A},+}=
\sum_{s\in\mathcal{S}_{\sigma}^+}f(s)g(s)N^+(s^{-1})
\]
for functions $f$ and $g$ such that the sum is absolutely convergent.
By standard arguments (cf. Subsection 4.5) we obtain 
\begin{equation}\label{linksymmetric}
\begin{split}
\bigl(Fg\bigr)(\gamma)&=\bigl(\widetilde{F}g\bigr)(\gamma)=
\frac{\mathcal{C}(s_0^{\ddagger\sigma})}{2^nn!}\bigl(F^+g\bigr)(\gamma),
\qquad g\in \mathcal{A}_+G^{-1},\\
\bigl(J_{\sigma}h\bigr)(x)&=\bigl(\widetilde{J}_{\sigma}h\bigr)(x)=
\mathcal{C}_\sigma(s_0^\ddagger)\bigl(J^+_{\sigma}h\bigr)(x),\qquad
h\in \mathcal{A}_+G_{\sigma\tau}.
\end{split}
\end{equation}
We obtain from these formulas the following result.
\begin{thm}\label{FFFF}
The transform $F^+$ defines a linear bijection
$F^+: \mathcal{A}_+G^{-1}\rightarrow 
\mathcal{A}_+G_{\sigma\tau}$, whose inverse is given by 
$(c_{\mathcal{A}}^+)^{-1}J^+_\sigma: \mathcal{A}_+G_{\sigma\tau}\rightarrow
\mathcal{A}_+G^{-1}$, where the constant
$c_{\mathcal{A}}^+$ is given by
\begin{equation}\label{cAplus}
 c_{\mathcal{A}}^+
=N^+_\sigma(s_0^{-1})\bigl(1,1\bigr)_{\mathcal{A},+}.
\end{equation}
Furthermore, we have the Plancherel formulas 
\begin{equation*}
\begin{split}
\lbrack F^+g,
F^+h\rbrack_{\mathcal{A},+,\sigma}&=
c_{\mathcal{A}}^+\,\bigl(g,h\bigr)_{\mathcal{A},+},\qquad
g,h\in\mathcal{A}_+G^{-1},\\
\bigl(J^+_{\sigma}g, J^+_{\sigma}h\bigr)_{\mathcal{A},+}&=
c_{\mathcal{A}}^+\,\lbrack g,h\rbrack_{\mathcal{A},+,\sigma},\qquad
g,h\in \mathcal{A}_+G_{\sigma\tau}.
\end{split}
\end{equation*}
\end{thm}
\begin{proof}
First note that for functions
$g,h\in\mathcal{F}(\mathcal{S}_\ddagger)$ which are $W_0$-invariant
under the dot-action,
\[\lbrack g,h\rbrack_{\mathcal{A},\sigma}=
\mathcal{C}_\sigma(s_0^\ddagger)\lbrack
g,h\rbrack_{\mathcal{A},+,\sigma},
\]
provided that the sums absolutely converge, cf. \eqref{JAsym}.
An analogous statement holds true for
$\bigl(\cdot,\cdot\bigr)_{\mathcal{A}}$, see \eqref{symmcontinuous}.
By the inversion formula and Plancherel formula for $F$ and by
\eqref{linksymmetric}, 
it then suffices to note that
\begin{equation}\label{DEplus}
\frac{2^nn!}{\mathcal{C}(s_0^{\ddagger\sigma})
\mathcal{C}_\sigma(s_0^\ddagger)}\,c_{\mathcal{A}}
=N^+_\sigma(s_0^{-1})\bigl(1,1\bigr)_{\mathcal{A},+},
\end{equation}
which follows from \eqref{Ndecomposition} and \eqref{symmcontinuous}.
\end{proof}

We finish the subsection by symmetrizing 
the explicit formulas in Proposition \ref{Fexplicit},
Proposition \ref{Jexplicit}, Corollary \ref{formulasnew}
and Proposition \ref{mixed}. We define constants $D_0^+$ and
$E_0^+$ by
\begin{equation}\label{D0plus}
D_0^+=2^nn!\frac{D_0}{\mathcal{C}(s_0^{\ddagger\sigma})},\qquad
E_0^+=\frac{E_0}{\mathcal{C}_\sigma(s_0^\ddagger)}
\end{equation}
with $D_0$ and $E_0$ given by \eqref{D0} and \eqref{E0},
respectively. Note that $D_0^+E_0^+=c_{\mathcal{A}}^+$
by \eqref{DE} and \eqref{DEplus}. 
Using standard symmetrization techniques, one obtains
the following proposition.
\begin{prop}\label{Fplusexplicit}
Let $s,u\in\mathcal{S}_\tau^+=\mathcal{S}_{\sigma\tau}^+$
and $v\in \mathcal{S}^+=\mathcal{S}_{\sigma\tau\sigma}^+$.\\
{\bf a)} The transform $F^+$ satisfies
\[
F^+\bigl(E_\tau^+(s;\cdot)G^{-1}\bigr)(\gamma)=
D_0^+\,G_{\tau\sigma\tau}(s)E_{\sigma\tau}^+(s;\gamma)G_{\sigma\tau}(\gamma).
\]
{\bf b)} The transform $J^+_\sigma$ satisfies 
\[J^+_\sigma\bigl(E_{\sigma\tau}^+(s;\cdot)G_{\sigma\tau}\bigr)(x)=
E_0^+\,G_{\tau\sigma\tau}(s)^{-1}E_\tau^+(s;x)G(x)^{-1}.
\]
{\bf c)} The following identities are valid:
\begin{equation*}
\begin{split}
D_0^+G_{\tau\sigma\tau}(s)G_{\tau\sigma\tau}(u)\,
&\lbrack E_{\sigma\tau}^+(s;\cdot)G_{\sigma\tau},
E_{\sigma\tau}^+(u;\cdot)G_{\sigma\tau}\rbrack_{\mathcal{A},+,\sigma}=\\
&\quad=E_0^+\bigl(E_\tau^+(s;\cdot)G^{-1},
E_{\tau}^+(u;\cdot)G^{-1}\bigr)_{\mathcal{A},+}.
\end{split}
\end{equation*}
{\bf d)} We have the integral evaluations
\[
\bigl(E^+(v;\cdot),E_{\tau}^+(s;\cdot)G^{-1}
\bigr)_{\mathcal{A},+}=
D_0^+G_{\sigma\tau}(v)G_{\tau\sigma\tau}(s)
E_{\sigma\tau}^+(s;v).
\]
\end{prop}


\section{The (non)symmetric Askey-Wilson function transform}
In this section we restrict attention to rank one ($n=1$). We
study nonsymmetric analogues of the 
(spherical) Fourier transform on the noncompact quantum
$\hbox{SU}(1,1)$ group, following the general philosophy of
Section 3.

Recall 
that the Jacobi function transform is a generalized Fourier transform
with kernel given by the Jacobi function. It has an
interpretation (for certain discrete parameter values) as the
spherical Fourier transform on $\hbox{SU}(1,1)$, see e.g. \cite{Koover}.
In recent papers of Koelink and the author, see \cite{KS1}, \cite{KS2} and
\cite{KS3}, a Fourier transform was defined and studied which admits an
interpretation as (spherical) Fourier transform on the noncompact
quantum group $\hbox{SU}_q(1,1)$. The transform, named the
Askey-Wilson function transform, is an integral transform with
kernel given by the so-called Askey-Wilson function. 

By Theorem \ref{linkhyp}
and Remark \ref{linkhypremark}{\bf b)}, the Askey-Wilson function
is precisely the symmetric rank one Cherednik kernel
$\mathfrak{E}^+$. This observation naturally leads to  
nonsymmetric variants
of the Askey-Wilson function transform, defined as integral 
transforms involving the normalized rank one Cherednik kernels
$\mathfrak{E}$ and $\mathfrak{E}_\ddagger$. These
transforms qualify as difference Fourier transforms in the sense of
Section 3. The underlying spaces are given explicitly as 
direct sums of two cyclic $\mathcal{H}$-modules.

In Subsection 8.1 we define the bilinear forms and the
cyclic $\mathcal{H}$-modules. In Subsection 8.2 we define
the (non)symmetric Askey-Wilson function transforms. In Subsection 8.3
we analyze the transforms on the first, ``classical'' cyclic
$\mathcal{H}$-module, which reduces to (the rank one case) of the extended
Macdonald-Koornwinder transform 
as discussed in the previous section.
In Subsection 8.4 we compute the image of the cyclic vector of the
second, ``strange'' $\mathcal{H}$-module under the transforms.
In Subsection 8.5 we prove algebraic Plancherel and inversion formulas
for the (non)symmetric Askey-Wilson function transform. In Subsection
8.6 we show how these results can be extended to the $L^2$-level
for the symmetric Askey-Wilson function transform, yielding new proofs
for the main results of \cite{KS2}.

In this section we keep the same generic 
assumptions on the difference multiplicity
function
$\alpha=(\mathbf{t},q^{\frac{1}{2}})=(t_0,u_0,t_1,u_1,q^{\frac{1}{2}})$
as in Section 5. Recall that in the rank one setup ($n=1$), 
the roots in $R_{nr}$ of medium length have disappeared, 
hence the associated 
parameter $t$ in the  multiplicity function $\mathbf{t}$ 
disappears (see \cite{NS} for the detailed treatment of the
polynomial theory in the rank one setup).  
In addition, an extra parameter $e\in\mathbb{C}^\times$ enters in the
definition of the
nonsymmetric Askey-Wilson function transform, which we only 
assume to be generic (unless stated explicitly otherwise).


\subsection{The bilinear forms}
It is well known from the theory of elliptic functions that there
exists, up to a multiplicative constant, a unique meromorphic
function
$\mathcal{P}_e=\mathcal{P}_{e}^{\alpha}\in\mathcal{M}=
\mathcal{M}(\mathbb{C}^\times)$
satisfying the invariance properties
\[\mathcal{P}_e(x^{-1})=\mathcal{P}_e(x),\qquad
\mathcal{P}_e(qx)=\mathcal{P}_e(x),\]
and with divisor on the elliptic curve
$\mathbb{C}^\times/q^{\mathbb{Z}}$ (written multiplicatively) given by
\[ \hbox{Div}(\mathcal{P}_e)=(d)+ (d^{-1})-(e)-(e^{-1}).
\]
Here and in the remainder of this section
we use the Askey-Wilson parametrization
\begin{equation}\label{AWparone}
\{a,b,c,d\}=\{t_1u_1,t_1u_1^{-1},q^{\frac{1}{2}}t_0u_0,
-q^{\frac{1}{2}}t_0u_0^{-1}\}
\end{equation}
for the difference multiplicity function $\alpha$, cf. \eqref{AWpar}.
The function $\mathcal{P}_e$ is closely related to the Weierstrass
$\mathcal{P}$-function. We fix $\mathcal{P}_e$ here by defining it
as a quotient of Jacobi theta-functions. For this we introduce the notation
$\theta(y_1,y_2,\ldots,y_n)=\theta(y_1)\theta(y_2)\cdots \theta(y_n)$
for products of the renormalized Jacobi theta function
\begin{equation}\label{Jacobi}
\theta(x)=\bigl(x,q/x;q\bigr)_{\infty}.
\end{equation}
Then we fix $\mathcal{P}_e=\mathcal{P}_e^{\alpha}
\in\mathcal{M}$ uniquely as the quotient
\begin{equation}\label{P}
\mathcal{P}_e(x)=\frac{\theta(dx,dx^{-1})}
{\theta(ex,ex^{-1})}.
\end{equation}
Since $\mathcal{P}_e\in\mathcal{M}$ is $\mathcal{W}$-invariant,
the associated multiplication
operator in $\hbox{End}_{\mathbb{C}}(\mathcal{M})$ commutes with the action 
of the double affine Hecke algebra $\mathcal{H}$ on $\mathcal{M}$.

We now define new weight functions $W_e=W_{e}^{\alpha}\in\mathcal{M}$
and $W_e^+=W_{e}^{+,\alpha}\in\mathcal{M}$
by
\begin{equation}\label{W}
W_e(x)=\mathcal{P}_e(x)\Delta(x),\qquad W_e^+(x)=\mathcal{P}_e(x)\Delta^+(x)
\end{equation}
with $\Delta=\Delta_\alpha$ the weight function for the
rank one Macdonald-Koornwinder polynomials (see \eqref{weight}) and
$\Delta^+=\Delta_\alpha^+$ the weight function for the symmetric
rank one Macdonald-Koornwinder polynomials (see \eqref{CDform}). 
Observe that the weight function
$W_e^+$ is $W_0$-invariant, i.e. $W_e^+(x^{-1})=W_e^+(x)$. In terms
of $q$-shifted factorials, the weight function $\Delta$ is given by
\[\Delta(x)=\frac{\bigl(x^2,qx^{-2};q\bigr)_{\infty}}
{\bigl(ax,qax^{-1},bx,qbx^{-1},cx,cx^{-1},dx,dx^{-1};
q\bigr)_{\infty}},
\]
see \eqref{decomposition}, \eqref{CDform} and \eqref{Deltaform},
whence $W_e\in\mathcal{M}$ can be expressed as 
\[W_e(x)= \frac{\bigl(x^2,qx^{-2},qx/d,q/dx;q\bigr)_{\infty}}
{\bigl(ax,qax^{-1},bx,qbx^{-1},cx,cx^{-1};q\bigr)_{\infty}
\theta(ex,ex^{-1})}.
\]
Similarly the weight function $W_e^+$ is given by
\[W_e^+(x)= \frac{\bigl(x^2,x^{-2},qx/d,q/dx;q\bigr)_{\infty}}
{\bigl(ax,ax^{-1},bx,bx^{-1},cx,cx^{-1};q\bigr)_{\infty}
\theta(ex,ex^{-1})}.
\]
For $\epsilon>0$ sufficiently small, let $C_\epsilon\subset
\mathbb{C}$ be a closed, counterclockwise oriented rectifiable Jordan
curve around the origin $0\in\mathbb{C}$ satisfying 
$C_\epsilon^{-1}=C_\epsilon$ (set-theoretically), and containing
the sequences
\begin{equation}\label{polesW}
\{aq^m,bq^m,cq^m,q^{1+m}/d \,\, | \,\, m\in\mathbb{N}\}
\cup \{eq^m \,\, | \,\, m\in\mathbb{Z},\,\, |eq^m|<\epsilon^{-1}\}
\end{equation}
in its interior, respectively the $eq^m$ ($m\in\mathbb{Z}$) with
$|eq^m|\geq \epsilon^{-1}$ in its exterior.  
In the special case that $a,b,c$ and $q/d$ have moduli $\leq 1$, we may
choose the contour $C_\epsilon$ to be the unit circle
in the complex plane with two deformations, one to include the
poles $eq^m$ with moduli $<\epsilon^{-1}$ and
one to exclude the poles $e^{-1}q^m$ with moduli $>\epsilon$.

For such $\epsilon>0$ sufficiently small, we define now two pairings
as follows. 
For meromorphic functions $g,h\in\mathcal{O}G_\tau$
we define the pairing
$\langle g,h\rangle_e^\epsilon=
\langle g,h\rangle_{e}^{\epsilon,\alpha}$ by the formula
\begin{equation}\label{bilinearepsilon}
\langle g,h\rangle_e^\epsilon=\frac{1}{2\pi i}
\int_{C_\epsilon}g(x)h(x^{-1})W_e(x)\frac{dx}{x}.
\end{equation}
Similarly, we define the pairing $\langle g,h\rangle_{e,+}^\epsilon=
\langle g,h\rangle_{e,+}^{\epsilon,\alpha}$ for meromorphic
functions $g,h\in\mathcal{O}G_\tau$ by
\begin{equation}
\langle g,h\rangle_{e,+}^\epsilon=\frac{1}{2\pi i}
\int_{C_\epsilon}g(x)h(x)W_e^+(x)\frac{dx}{x}.
\end{equation}
By Cauchy's Theorem, the pairings $\langle g,h\rangle_e^\epsilon$
and $\langle g,h\rangle_{e,+}^\epsilon$ are 
independent of the choice of contour $C_\epsilon$ satisfying the
specific defining conditions as stated above. 

We now define two subspaces $V^{cl}_e=V^{cl,\alpha}_{e}$
and $V^{str}=V^{str,\alpha}$ of $\mathcal{O}G_{\tau}$ 
by
\begin{equation}\label{modules}
V^{cl}_e=\mathcal{A}\mathcal{P}_e^{-1}G^{-1},\qquad
V^{str}=\mathcal{A}G_\tau.
\end{equation}
Furthermore, we write
$V^{cl}_{e,+}=V^{cl,\alpha}_{e,+}$ and $V^{str}_+=V^{str,\alpha}_{+}$
for the associated subspaces of $W_0$-invariant functions,
\begin{equation}\label{modulessym}
V^{cl}_{e,+}=\mathcal{A}_+\mathcal{P}_e^{-1}G^{-1},\qquad
V^{str}_+=\mathcal{A}_+G_\tau.
\end{equation}
The superscripts 
{\it{cl}} and {\it{str}} stand for ``classical'' and
``strange'' respectively. This terminology is motivated by the fact
that 
$V^{cl}_e$ (respectively $V^{str}$) covers the contributions
of the Plancherel measure of the spherical Fourier transform on the
quantum group $\hbox{SU}_q(1,1)$ which arise
from the unitary principal series
representations (respectively strange series representations),
see e.g. \cite{KS1} and \cite{St3}.

Observe that $V^{cl}_e$ and 
$V^{str}$ are cyclic $\mathcal{H}$-submodules of $\mathcal{M}$
with corresponding cyclic vectors $\mathcal{P}_e^{-1}G^{-1}$ and
$G_\tau$, respectively.
In fact, for any $p\in\mathcal{A}$ and $X\in\mathcal{H}$ we have
\begin{equation}
\begin{split}
X\bigl(p\mathcal{P}_e^{-1}G^{-1}\bigr)&=\bigl(\tau(X)p\bigr)
\mathcal{P}_e^{-1}G^{-1},\\
X\bigl(pG_\tau\bigr)&=\bigl(\tau_\tau^{-1}(X)p\bigr)G_\tau
\end{split}
\end{equation}
and the action of $\mathcal{H}_\tau$ preserves
the subspace $\mathcal{A}$ of Laurent polynomials. 
Any function
$g\in V^{cl}_e$ vanishes at points $x\in
(eq^{\mathbb{Z}})^{\pm 1}$, hence
$V^{cl}_e\cap V^{str}=\{0\}$. We define 
the $\mathcal{H}$-module $M_e=M_{e}^{\alpha}$ by
\begin{equation}
M_e=V^{cl}_e\oplus V^{str}\subset \mathcal{O}G_{\tau}.
\end{equation}
The subspace $M_e^+=M_{e}^{+,\alpha}\subset M_e$
consisting of $W_0$-invariant functions in $M_e$ is given by
the direct sum
\[ M_e^+=V_{e,+}^{cl}\oplus V_{+}^{str}.
\]
Observe that the identity $G_{\ddagger\tau}=G^{-1}$ implies that $M_e$ is also
$\mathcal{H}_\ddagger$-stable under the natural action of
$\mathcal{H}_\ddagger$ on $\mathcal{M}$ as $q^{-1}$-difference
reflection operators. The explicit formulas are given by
\begin{equation*}
\begin{split}
X\bigl(p\mathcal{P}_e^{-1}G^{-1}\bigr)&=
\bigl(\tau_{\ddagger\tau}^{-1}(X)p\bigr)\mathcal{P}_e^{-1}G^{-1},\\
X\bigl(pG_\tau\bigr)&=\bigl(\tau_\ddagger(X)p\bigr)G_\tau
\end{split}
\end{equation*}
for $p\in\mathcal{A}$ and for $X\in\mathcal{H}_\ddagger$.
\begin{lem}\label{existence} 
For $g,h\in M_e$, the limits 
\[\langle g,h\rangle_e=\lim_{\epsilon\searrow 0}\,
\langle g,h\rangle_e^\epsilon,\qquad\quad
\langle g,h\rangle_{e,+}=\lim_{\epsilon\searrow 0}\,
\langle g,h\rangle_{e,+}^\epsilon
\]
exist.
\end{lem}
\begin{proof}
We prove the lemma 
for the pairing $\langle \cdot,\cdot\rangle_{e,+}$ (the proof for
$\langle \cdot,\cdot\rangle_{e}$ is similar).
Let $g,h\in M_e$. Let $\epsilon>0$ be sufficiently
small. We may rewrite the integral over $C_\epsilon$ in the definition
of $\langle g,h\rangle_{e,+}^\epsilon$ by an integral over an 
($\epsilon$-independent) closed, rectifiable Jordan curve $C$ on the
cost of picking up residues at points of the form $eq^l$ and
$e^{-1}q^{-l}$ for $m_\epsilon\leq l\leq l_0$, where $m_\epsilon\in
\mathbb{Z}$ is the smallest integer such that 
$|eq^{m_\epsilon}|<\epsilon^{-1}$,
and with $l_0\in\mathbb{Z}$ some suitably chosen, fixed integer.
Hence convergence of the limit will follow from the
convergence of the series 
\begin{equation}\label{convsum}
\sum_{l\leq l_0}g(y_l^{\pm 1})h(y_l^{\pm 1})\underset{x=y_l^{\pm 1}}
{\hbox{Res}}\left(
\frac{W_e^+(x)}{x}\right),
\end{equation}
with $y_l=eq^l$. In case that $g\in V_e^{cl}$ or
$h\in V_e^{cl}$ these sums vanish. When 
$g,h\in V^{str}=\mathcal{A}G_\tau$, then 
\begin{equation}\label{gGaussianbound}
|g(y_l^{\pm 1})|\leq c_0c_1^{|l|}q^{l^2/2},\qquad \forall\,l\leq l_0
\end{equation}
for some constants $c_0,c_1>0$, and similarly for $h$. Furthermore,
\[
\left|\underset{x=y_l^{\pm 1}}
{\hbox{Res}}\left(\frac{W_e^+(x)}{x}\right)\right|\leq d_0d_1^{|l|},\qquad
\forall\,l\leq l_0
\]
for some constants $d_0,d_1>0$ (which e.g. follows from the explicit expression
\cite[(5.8)]{KS2} for the residue of the weight function $W_e^+$ at
$y_l^{\pm 1}$).
Hence the Gaussian contributions $q^{l^2/2}$ in the asymptotics of
$g$ and $h$ force the convergence of \eqref{convsum}. This completes
the proof of the lemma.
\end{proof}

\begin{prop}\label{adjointAW}
The bilinear form $\langle \cdot,\cdot\rangle_e$ on 
$M_e$ induces the anti-isomorphism
$\ddagger$ on $\mathcal{H}$. In other words,
\[\langle Xg,h\rangle_e=\langle g,\ddagger(X)h\rangle_e,
\]
for $X\in\mathcal{H}$ and $g,h\in M_e$.
\end{prop}
\begin{proof}
It suffices to prove the proposition for $X=T_0$ and
$X=T_1$.
By the explicit form of the operators
$T_j$, we have for $j=0$ and for $j=1$,
\[
\langle T_jg,h\rangle_e^\epsilon-\langle g,\ddagger(T_j)h\rangle_e^\epsilon
=\frac{t_j^{-1}}{2\pi i}\int_{C_\epsilon}
\left\{(r_jg)(x)\widetilde{h}(x)-g(x)(r_j\widetilde{h})(x)\right\}
c_{a_j}(x)W_e(x)\frac{dx}{x},
\]
with $\widetilde{h}(x)=h(x^{-1})$ and with the $r_j$'s acting as
constant coefficient $q$-difference reflection operators.
For $j=1$, the function
$x\mapsto c_{a_1}(x)W_e(x)$ is $r_1$-invariant,
hence the right hand side vanishes by the inversion-invariance of
the contour $C_\epsilon$. Taking the limit $\epsilon\searrow 0$ 
and using Lemma \ref{existence} then
yields $\langle T_1g,h\rangle_e=\langle g,\ddagger(T_1)h\rangle_e$.

For $j=0$, the function
$x\mapsto c_{a_0}(x)W_e(x)$ is $r_0$-invariant, hence
\[\langle T_0g,h\rangle_e^\epsilon-\langle g,
\ddagger(T_0)h\rangle_e^\epsilon
=\frac{-t_0^{-1}}{2\pi i}\int_{C_\epsilon-qC_\epsilon}
g(x)\widetilde{h}(qx^{-1})c_{a_0}(x)W_e(x)\frac{dx}{x}.
\]
A straightforward computation shows that the poles in $\mathbb{C}^\times$ 
of the integrand are simple and contained in the zero set of
\[x\mapsto
\bigl(ax,qax^{-1},bx,qbx^{-1},cx,qcx^{-1},qx/d,q^2/dx;q\bigr)_{\infty}
\theta(ex,ex^{-1}).
\]
Let now $m\in\mathbb{Z}$ be the smallest integer satisfying
$|eq^m|<\epsilon^{-1}$. By the assumptions on the contour 
$C_\epsilon$ and by Cauchy's Theorem, we pick up poles at $eq^m$
and at $e^{-1}q^{1-m}$ 
when shifting $C_\epsilon$ to $qC_\epsilon$, whence
\begin{equation*}
\begin{split}
\langle &T_0g,h\rangle_e^\epsilon-\langle g,\ddagger(T_0)h\rangle_e^\epsilon=\\
&=t_0^{-1}\bigl(g(q^{1-m}e^{-1})\widetilde{h}(eq^m)-
g(eq^m)\widetilde{h}(q^{1-m}e^{-1})\bigr)
\underset{x=eq^m}{\hbox{Res}}\left(\frac{c_{a_0}(x)W_e(x)}{x}\right).
\end{split}
\end{equation*}
Now take the limit $\epsilon\searrow 0$ and use 
the bounds derived in the proof
of Lemma \ref{existence}. This implies that
$\langle T_0g,h\rangle_e=\langle g,\ddagger(T_0)h\rangle_e$,
as desired.
\end{proof}

\begin{rem}\label{betterstat}
Observe that the conditions on the functions $g$ and $h$ in
Lemma \ref{existence} and Proposition \ref{adjointAW} may be
relaxed. For instance, Lemma \ref{existence} and Proposition \ref{adjointAW}
hold true when $g\in M_e$ and $h\in \mathcal{O}G_\tau$ with $h$ satisfying 
\begin{equation}\label{boundanalytic}
|h(y_l^{\pm 1})|\leq c_0c_1^{|l|},\qquad \forall\,l\leq l_0
\end{equation}
for some constants $c_0,c_1>0$, where $l_0\in\mathbb{Z}$ is some
arbitrary, fixed integer and $y_l=eq^l$. If
$h=h_\gamma\in\mathcal{O}G_\tau$ furthermore depends analytically on
an additional parameter $\gamma\in\mathbb{C}^\times$ and
\eqref{boundanalytic} holds 
true uniformly for $\gamma$ in compacta of $\mathbb{C}^\times$, then
the proof of Lemma \ref{existence} in addition implies that $\langle
g,h_\gamma\rangle_e$ depends analytically on
$\gamma\in\mathbb{C}^\times$, for all $g\in M_e$.  
\end{rem}


\subsection{The transforms} 
In order to define the (non)symmetric Askey-Wilson function
transform we need to establish certain bounds for the normalized
Cherednik kernels $\mathfrak{E}$ and $\mathfrak{E}_\ddagger$
associated to $\alpha$ and $\alpha_\ddagger$ respectively, as well
as for the normalized symmetric Cherednik kernel $\mathfrak{E}^+$. 
In fact, we prove bounds for the analytic functions
$G_{\sigma\tau}^{-1}\mathfrak{E}(\cdot,y_l^{\pm 1})$,
$G_{\sigma\tau}^{-1}\mathfrak{E}_\ddagger(\cdot,y_l^{\pm 1})$
and $G_{\sigma\tau}^{-1}\mathfrak{E}^+(\cdot,y_l)$
in $l\in\mathbb{Z}$, where $y_l=eq^l$.
\begin{lem}\label{Chbounds}
For any compacta $K\subset (\mathbb{C}^\times)^n$, there exist 
constants $C,D>0$ \textup{(}depending on $K$\textup{)}, such that
\begin{equation*}
\begin{split}
|G_{\sigma\tau}(\gamma)^{-1}\mathfrak{E}(\gamma,y_l^{\pm 1})|
&\leq CD^{|l|},\\
|G_{\sigma\tau}(\gamma)^{-1}\mathfrak{E}_\ddagger(\gamma,y_l^{\pm 1})|
&\leq CD^{|l|},\\
|G_{\sigma\tau}(\gamma)^{-1}\mathfrak{E}^+(\gamma,y_l)|
&\leq CD^{|l|},
\end{split}
\end{equation*}
for all $l\in\mathbb{Z}$ and all $\gamma\in K$, where
$y_l=eq^l$.
\end{lem}
\begin{proof}
Recurrence relations for the normalized Cherednik kernels
$\mathfrak{E}$ and $\mathfrak{E}_\ddagger$ are essentially the same as
Pieri formulas for the polynomial kernel
$\mathfrak{E}_{\mathcal{A}}$, since the kernels satisfy the same
transformation behaviour under the action of the double affine
Hecke algebra (see Proposition \ref{allrelations} and Theorem 
\ref{Cherednikkernels}). Hence the proof of the bounds for
the Macdonald-Koornwinder polynomials (see Proposition \ref{MKbounds}
and the appendix) can be easily adjusted to obtain the desired bounds for the
analytic functions $G_{\sigma\tau}^{-1}\mathfrak{E}(\cdot,y_l^{\pm 1})$
and $G_{\sigma\tau}^{-1}\mathfrak{E}_\ddagger(\cdot,y_l^{\pm 1})$.
We leave the precise details to the reader. The bounds for
$\mathfrak{E}^+$ follows easily from the bounds for $\mathfrak{E}$, 
using e.g. the formula
$\mathfrak{E}^+(\gamma,x)=\bigl(C_+\mathfrak{E}(\gamma,\cdot)\bigr)(x)$.
\end{proof}

\begin{prop}\label{85} 
{\bf a)} The assignment
\[\bigl(\mathcal{F}_eg\bigr)(\gamma)=
\langle g,\mathfrak{E}_\ddagger(\gamma^{-1},\cdot)\rangle_e
\]
for $g\in M_e$ defines
a linear map 
$\mathcal{F}_e=\mathcal{F}_{e}^{\alpha}: 
M_e\rightarrow \mathcal{O}G_{\sigma\tau}$. Furthermore,
$\mathcal{F}_e$ is a Fourier transform associated with $\sigma$.\\
{\bf b)} The assignment
\[\bigl(\mathcal{J}_eg\bigr)(\gamma)=
\langle g,I\mathfrak{E}(\gamma,\cdot)\rangle_e
\]
for $g\in M_e$, with $(Ih)(x)=h(x^{-1})$, defines a linear map 
$\mathcal{J}_e=\mathcal{J}_{e}^{\alpha}: 
M_e\rightarrow \mathcal{O}G_{\sigma\tau}$. Furthermore,
$\mathcal{J}_e$ is a Fourier transform associated with
$\sigma_\sigma^{-1}$.
\end{prop}
\begin{proof}
{\bf a)} Remark \ref{betterstat} and
Lemma \ref{Chbounds} imply that the map
\[\gamma\mapsto \langle g,
G_{\sigma\tau}(\gamma)^{-1}\mathfrak{E}_\ddagger(\gamma^{-1},\cdot)\rangle_e 
\]
for fixed $g\in M_e$ defines an analytic function on
$\mathbb{C}^\times$. Hence $\mathcal{F}_e$ defines
a linear map $\mathcal{F}_e: M_e\rightarrow \mathcal{O}G_{\sigma\tau}$.
By Proposition \ref{adjointAW}, Remark \ref{betterstat} and 
the general arguments of Section 3 it is clear that $\mathcal{F}_e$
is a Fourier transform associated with $\sigma$. 
The proof of {\bf b)} is similar. 
\end{proof}

\begin{cor}\label{symmAW}
The assignment
\[\bigl(\mathcal{F}_e^+g\bigr)(\gamma)=
\langle g,\mathfrak{E}^+(\gamma,\cdot)\rangle_{e,+}
\]
for $g\in M_e^+$ defines a linear map 
$\mathcal{F}_e^+=\mathcal{F}_{e}^{+,\alpha}: 
M_e^+\rightarrow \mathcal{O}G_{\sigma\tau}$. 
Furthermore,
\[\mathcal{F}_eg=\mathcal{J}_eg=
\frac{\mathcal{C}(s_0^{\ddagger\sigma})}{2}\mathcal{F}_e^+g
\]
for all $g\in M_e^+$.
\end{cor}
\begin{proof}
By similar arguments as in the proof of Proposition
\ref{85}, we have that
$\mathcal{F}_e^+$ defines a linear map 
$\mathcal{F}_e^+: M_e^+\rightarrow \mathcal{O}G_{\sigma\tau}$.
By \eqref{stableC}, Lemma \ref{elementary}, Proposition
\ref{adjointAW}, Remark \ref{betterstat} 
and Theorem \ref{symmmain}{\bf a)}, we obtain
\[\bigl(\mathcal{F}_eg\bigr)(\gamma)=\bigl(\mathcal{J}_eg\bigr)(\gamma)=
\langle g,\mathfrak{E}^+(\gamma,\cdot)\rangle_e
\] 
for $g\in M_e^+$. Symmetrizing the integral 
using \eqref{symmistrivial} and using the decomposition
$W_e(x)=\mathcal{C}(x)W_e^+(x)$ gives
\[\langle g,\mathfrak{E}^+(\gamma,\cdot)\rangle_e=
\frac{\mathcal{C}(s_0^{\ddagger\sigma})}{2}
\langle g,\mathfrak{E}^+(\gamma,\cdot)\rangle_{e,+}
\]
for $g\in M_e^+$, which completes the proof.
\end{proof}
By Theorem \ref{linkhyp} the 
transform $\mathcal{F}_e^+$ is, up to a multiplicative constant,
precisely the Askey-Wilson function transform as defined and studied
in \cite{KS2}. The only difference
is that the transform $\mathcal{F}_e^+$ in \cite{KS2} is (initially) defined 
on a certain space of compactly supported functions, while  
$\mathcal{F}_e^+$ in this subsection is defined
on the space $M_e^+$ (which is a more natural subspace
from the viewpoint of
double affine Hecke algebras). It turns out though that
for restricted parameter values, 
their continuous extensions to the $L^2$-level do coincide
(see \cite{KS2}, \cite{St3} and Subsection 8.6).

In order to distinguish the two transforms $\mathcal{F}_e$ 
and $\mathcal{F}_e^+$, we use the following terminology throughout
the remainder of the paper.
\begin{Def}
{\bf a)} The transform $\mathcal{F}_e: M_e\rightarrow
\mathcal{O}G_{\sigma\tau}$ is called the Askey-Wilson
function transform.\\
{\bf b)} The transform $\mathcal{F}_e^+: M_e^+\rightarrow
\mathcal{O}G_{\sigma\tau}$ is called the symmetric Askey-Wilson
function transform. 
\end{Def}


\subsection{The classical part of the transforms}

In this subsection we consider the difference Fourier transforms
$\mathcal{F}_e$ and $\mathcal{J}_e$ on
the cyclic $\mathcal{H}$-submodule $V_e^{cl}\subset M_e$. This
is related to the extended (rank one) Macdonald-Koornwinder
transforms $F$ and $\widetilde{F}$ of the previous section in the
following way. 
\begin{lem}\label{reduction2}
{\bf a)} For $p\in\mathcal{A}$,
\[\mathcal{F}_e\bigl(p\mathcal{P}_e^{-1}G^{-1}\bigr)(\gamma)=
F\bigl(pG^{-1}\bigr)(\gamma),\qquad
\mathcal{J}_e\bigl(p\mathcal{P}_e^{-1}G^{-1}\bigr)(\gamma)=
\widetilde{F}\bigl((Ip)G^{-1}\bigr)(\gamma),
\]
where $I$ is the inversion operator $(Ig)(x)=g(x^{-1})$.\\
{\bf b)} For $p\in\mathcal{A}_+$, 
\[\mathcal{F}_e^+\bigl(p\mathcal{P}_e^{-1}G^{-1}\bigr)(\gamma)=
F^+\bigl(pG^{-1}\bigr)(\gamma),
\] 
with $F^+$ the extended \textup{(}rank one\textup{)} 
Macdonald-Koornwinder transform
defined by \eqref{Fplus}.
\end{lem}
\begin{proof}
{\bf a)} Let $\mathcal{T}=\mathcal{T}_\alpha$ be a deformed circle as
defined in Subsection 7.1.
For $p\in\mathcal{A}$ we have
\begin{equation*}
\begin{split}
\mathcal{F}_e\bigl(p\mathcal{P}_e^{-1}G^{-1}\bigr)(\gamma)&=
\frac{1}{2\pi i}\int_{\mathcal{T}_\tau}p(x)
\mathfrak{E}_\ddagger(\gamma^{-1},x^{-1})G(x)^{-1}\Delta(x)\frac{dx}{x}\\
&=\bigl(pG^{-1},\mathfrak{E}_\ddagger(\gamma^{-1},\cdot)\bigr)_{\mathcal{A}}\\
&=F\bigl(pG^{-1}\bigr)(\gamma).
\end{split}
\end{equation*}
Here the first equality holds by Cauchy's Theorem since
$\mathcal{P}_e^{-1}W_e=\Delta$ is regular at the points
$x\in (eq^{\mathbb{Z}})^{\pm 1}$. The proof for $\mathcal{J}_e$
is similar.\\
{\bf b)} The proof is similar to the proof of {\bf a)}.
\end{proof}

Invoking the main results on the extended Macdonald-Koornwinder
transform,
we arrive at the following result.
\begin{cor}\label{classicalpart}
{\bf a)} The Askey-Wilson function transform $\mathcal{F}_e$ 
restricts to a linear bijection
$\mathcal{F}_e: V_e^{cl}\rightarrow V^{str,\sigma}$. Explicitly,
we have 
\begin{equation}\label{ee}
\mathcal{F}_e\bigl(E_\tau(s;\cdot)\mathcal{P}_e^{-1}G^{-1}\bigr)(\gamma)=
D_0G_{\tau\sigma\tau}(s)E_{\sigma\tau}(s;\gamma)G_{\sigma\tau}(\gamma)
\end{equation}
for all $s\in\mathcal{S}_\tau=\mathcal{S}_{\sigma\tau}$, 
with $D_0=D_0^\alpha$ the constant 
\begin{equation}\label{Dnul}
D_0=\frac{\bigl(bc,d/a,q/ad;q\bigr)_{\infty}}
{\bigl(q,qab,ac;q\bigr)_{\infty}}.
\end{equation}
{\bf b)} The difference Fourier transform 
$\mathcal{J}_e$ restricts to a linear bijection
$\mathcal{J}_e: V_e^{cl}\rightarrow V^{str,\sigma}$. Explicitly,
we have 
\[\mathcal{J}_e\bigl(I(E_{\ddagger\tau}(s^{-1};\cdot))\mathcal{P}_e^{-1}
G^{-1}\bigr)(\gamma)=D_0G_{\tau\sigma\tau}(s)E_{\sigma\tau}(s;\gamma)
G_{\sigma\tau}(\gamma)
\]
for all $s\in\mathcal{S}_\tau$, with $D_0=D_0^\alpha$ the constant 
\eqref{Dnul}.\\
{\bf c)} The symmetric Askey-Wilson function transform
$\mathcal{F}_e^+$ restricts to a linear bijection
$\mathcal{F}_e^+: V_{e,+}^{cl}\rightarrow V_{+}^{str,\sigma}$. 
Explicitly, we have 
\[\mathcal{F}_e^+\bigl(E_\tau^+(s;\cdot)\mathcal{P}_e^{-1}G^{-1}\bigr)(\gamma)=
D_0^+G_{\tau\sigma\tau}(s)E_{\sigma\tau}^+(s;\gamma)G_{\sigma\tau}(\gamma)
\]
for all $s\in\mathcal{S}_\tau^+=\mathcal{S}_{\sigma\tau}^+$, 
with $D_0^+=D_0^{+,\alpha}$ the constant 
\begin{equation}\label{Dnulplus}
D_0^+=\frac{2\bigl(bc,d/a,q/ad;q\bigr)_{\infty}}
{\bigl(q,ab,ac;q\bigr)_{\infty}}.
\end{equation}
\end{cor}
\begin{proof}
{\bf a)} By Lemma \ref{reduction2} and Proposition \ref{Fexplicit},
formula \eqref{ee} holds with the constant $D_0$ given by \eqref{D0}. 
Now for rank one ($n=1$),
$\mathcal{C}(s_0^{\ddagger\sigma})=1-ab$, hence \eqref{D0} reduces
to \eqref{Dnul}. The proof of {\bf b)} and {\bf c)} are similar,
now using Proposition \ref{explicittilde} and Proposition
\ref{Fplusexplicit}, respectively. 
\end{proof}


\subsection{The Fourier transform of the Gaussian}
We consider the transforms $\mathcal{F}_e$ and $\mathcal{J}_e$ on
the cyclic $\mathcal{H}$-submodule $V^{str}=\mathcal{A}\,G_{\tau}$. 
Since $\mathcal{F}_e$ 
and $\mathcal{J}_e$ are Fourier transforms associated with $\sigma$
and $\sigma_\sigma^{-1}$ respectively, it suffices to evaluate
the image of the cyclic vector $G_\tau$ under $\mathcal{F}_e$
and under $\mathcal{J}_e$. This amounts to the same thing, since
Corollary \ref{symmAW} implies
\begin{equation}\label{samething}
\bigl(\mathcal{F}_eG_\tau\bigr)(\gamma)=
\bigl(\mathcal{J}_eG_\tau\bigr)(\gamma)=
\frac{\mathcal{C}(s_0^{\ddagger\sigma})}{2}
\bigl(\mathcal{F}_e^+G_\tau\bigr)(\gamma).
\end{equation}
We start with
the following one-variable $q$-analogue of the Macdonald-Mehta integral.
\begin{lem}\label{MM}
We have the explicit evaluation
\[\bigl(\mathcal{F}_e^+G_\tau\bigr)(s_0)=
\frac{2}{\bigl(q,ab,ac,bc;q\bigr)_{\infty}}
\frac{\theta(abce)}{\theta(ae,be,ce)}.
\]
\end{lem}
\begin{proof}
By the polynomial reduction, we have $\mathfrak{E}^+(s_0,\cdot)=1$,
hence
\begin{equation}\label{SSSS}
\bigl(\mathcal{F}_e^+G_\tau\bigr)(s_0)=
\lim_{\epsilon\searrow 0}\frac{1}{2\pi i}\int_{C_\epsilon}
G_\tau(x)W_e^+(x)\frac{dx}{x}.
\end{equation}
Since 
\begin{equation}\label{integrandexplicit}
G_\tau(x)W_e^+(x)=\frac{\bigl(x^2,1/x^{2};q\bigr)_{\infty}}
{\bigl(ax,a/x,bx,b/x,cx,c/x;q\bigr)_{\infty}\theta(ex,e/x)},
\end{equation}
the right hand side of \eqref{SSSS}
can be easily matched with the one variable 
$q$-Macdonald-Mehta integral \cite[(5.9)]{St3}, with the parameter
$u$ in \cite[(5.9)]{St3} taken to be $e$
(observe in particular that $\bigl(\mathcal{F}_e^+G_\tau\bigr)(s_0)$
does not depend on the parameter $d$). 
Its evaluation (see \cite[Thm. 5.5]{St3}) 
yields the desired result. 
\end{proof}
\begin{rem}
One of the goals in this section is to prove the main results on
the symmetric
Askey-Wilson function transform (see \cite{KS2} and \cite{St3}) 
using only affine Hecke algebra techniques and some explicit 
``constant term'' evaluations. 
It is therefore noteworthy to mention that
two proofs of the evaluation of the
above Macdonald-Mehta type integral are available 
 which only use some direct basic
hypergeometric series manipulations (see \cite[Appendix B]{St3}).
\end{rem}

\begin{prop}\label{FourierGaussian}
The image of the Gaussian $G_\tau$ under the symmetric Askey-Wilson
function transform $\mathcal{F}_e^+$ is given by
\[\bigl(\mathcal{F}_e^+G_\tau\bigr)(\gamma)=
\frac{2\bigl(q/ad;q\bigr)_{\infty}}{\bigl(q,ab,ac;q\bigr)_{\infty}
\theta(ae,be,ce,qe/d)}
\mathcal{P}_{e_\sigma}^\sigma(\gamma)^{-1}G_\sigma(\gamma)^{-1},
\]
with the parameter $e_\sigma\in\mathbb{C}$ defined by
\[e_\sigma=-\frac{q^{\frac{1}{2}}}{u_0t_1e}.
\]
\end{prop}
\begin{proof}
We consider the $W_0$-invariant meromorphic function
$f\in\mathcal{O}G_{\sigma}G_{\sigma\tau}$ defined by 
\begin{equation}\label{f}
 f(\gamma)=\frac{2}{\mathcal{C}(s_0^{\ddagger\sigma})}\,
G_\sigma(\gamma)\bigl(\mathcal{F}_eG_\tau\bigr)(\gamma)
=G_\sigma(\gamma)\bigl(\mathcal{F}_e^+G_\tau\bigr)(\gamma).
\end{equation}
We first show that $f$ is invariant under the action of
$r_0\in\mathcal{W}$, i.e. that $f(q\gamma^{-1})=f(\gamma)$. We compute
\begin{equation*}
\begin{split}
(T_0^{\sigma\tau}f)(\gamma)&=\frac{\mathcal{C}(s_0^{\ddagger\sigma})}{2}\,
G_\sigma(\gamma)\mathcal{F}_e\bigl(G_\tau
(\tau_\tau^{-1}\circ\sigma^{-1}\circ\tau_{\sigma}^{-1})(T_0^{\sigma\tau})(1)
\bigr)(\gamma)\\
&=\frac{\mathcal{C}(s_0^{\ddagger\sigma})}{2}\,
G_\sigma(\gamma)\mathcal{F}_e\bigl(G_\tau
(Y_1^{\tau}T_1^{\tau}{}^{-1})(1)\bigr)(\gamma)\\
&=u_0f(\gamma),
\end{split}
\end{equation*}
since $Y_1^{\tau}T_1^{\tau}{}^{-1}$ acts on $1\in\mathcal{A}$
as multiplication by the constant $u_0t_1t_1^{-1}=u_0$.
By the explicit form of the $q$-difference reflection
operator $T_0^{\sigma\tau}$, we conclude that $f$ is $r_0$-invariant. 
Hence $f\in\mathcal{M}$ is $\mathcal{W}$-invariant. 
So $f$ may be regarded as a meromorphic function
on the elliptic curve $T=\mathbb{C}^\times/q^{\mathbb{Z}}$ (written
multiplicatively). The only possible poles of $f$ on $T$ 
are at most simple and located at
$-q^{\frac{1}{2}+\mathbb{Z}}u_1u_0^{-1}$ and
$-q^{-\frac{1}{2}+\mathbb{Z}}u_1^{-1}u_0$, since $f\in 
\mathcal{O}G_{\sigma}G_{\sigma\tau}$ and
\[
G_{\sigma}(\gamma)G_{\sigma\tau}(\gamma)=
\frac{1}{\theta\bigl(-q^{\frac{1}{2}}u_1u_0^{-1}\gamma,
-q^{\frac{1}{2}}u_1u_0^{-1}\gamma^{-1}\bigr)}.
\]
By standard elliptic function theory it follows that
\begin{equation}\label{fgoal}
f(\gamma)=C\mathcal{P}_{\tilde{e}}^\sigma(\gamma)^{-1}=
C\,\frac{\theta\bigl(\tilde{e}\gamma,\tilde{e}\gamma^{-1}\bigr)}
{\theta\bigl(-q^{\frac{1}{2}}u_1u_0^{-1}\gamma,
-q^{\frac{1}{2}}u_1u_0^{-1}\gamma^{-1}\bigr)}
\end{equation}
for some $\tilde{e}\in \mathbb{C}^\times$ 
and some $C\in\mathbb{C}$. 
The choice of $\tilde{e}\in\mathbb{C}^\times$ in \eqref{fgoal}
is not unique.
In fact, if \eqref{fgoal} is valid for $\tilde{e}=u$ for some choice of
constant $C$, then all other possible choices for $\tilde{e}$ are given
by $(uq^{\mathbb{Z}})^{\pm 1}$.

We fix now a pair $(\tilde{e},C)$ such that equation \eqref{fgoal}
is satisfied. The next step is to 
derive explicit identities for $C$ and $\tilde{e}$
using the evaluation of the Macdonald-Mehta type integral
(see Lemma \ref{MM}) and using symmetries in the four Askey-Wilson
parameters $a,b,c$ and $d$. These additional identities in $C$ and
$\tilde{e}$ can be solved explicitly and lead to the
explicit expressions for $C$ and $\tilde{e}$ as given in the statement
of the proposition.

We first evaluate $f(s_0)=f(t_0t_1)$ in two different ways.
Since $G_\sigma(t_0t_1)=\bigl(ad,q/bc;q\bigr)_{\infty}^{-1}$,
we have by \eqref{f} and by Lemma \ref{MM},
\[f(t_0t_1)=\frac{2}{\bigl(q,ab,ac,ad;q\bigr)_{\infty}\theta(bc)}
\frac{\theta(abce)}{\theta(ae,be,ce)}.
\]
In particular we have $f(t_0t_1)\not=0$, which implies that
$\tilde{e}\not\in (t_0t_1q^{\mathbb{Z}})^{\pm 1}$ and $C\not=0$.

On the other hand, by \eqref{fgoal},
\[ f(t_0t_1)=C\frac{\theta(t_0t_1\tilde{e},t_0^{-1}t_1^{-1}\tilde{e})}
{\theta(ad,bc)}.
\]
Combining the two identities, we obtain an expression for the constant
$C$ in terms of $\tilde{e}$:
\begin{equation}\label{Ce}
C=\frac{2\bigl(q/ad;q\bigr)_{\infty}}{\bigl(q,ab,ac;q\bigr)_{\infty}}
\frac{\theta(abce)}{\theta(ae,be,ce,t_0t_1\tilde{e},
t_0^{-1}t_1^{-1}\tilde{e})}.
\end{equation}

In order to find $\tilde{e}$ explicitly,
we mimic the previous approach by evaluating 
$f(-t_0^{-1}t_1)$ in two different ways. In order to do so,
we consider a new multiplicity function $\beta$ by
\[\beta=(t_0,u_0,t_1,-u_1^{-1},q^{\frac{1}{2}}),
\]
i.e. the value $u_1$ of the multiplicity function $\mathbf{t}$ 
at the orbit $\mathcal{W}a_1^\vee$ is replaced by $-u_1^{-1}$.
In terms of the Askey-Wilson parametrization \eqref{AWparone},
this amounts to interchanging the role of $a$ and $b$.  
Note that $\beta_\sigma=(-u_1^{-1},u_0,t_1,t_0,q^{\frac{1}{2}})$,
so that $s_0^{\beta_\sigma}=-u_1^{-1}t_1$. Using the duality and the
polynomial reduction of the symmetric Cherednik kernel
$\mathfrak{E}^+$ (see Theorem
\ref{symmmain}), as well as the symmetry of $\mathfrak{E}^+$
when interchanging the role of $a$ and $b$ (see Proposition
\ref{Cherednikpar}), we derive
\begin{equation*}
\begin{split}
\mathfrak{E}_{\alpha_\sigma}^+(-u_1^{-1}t_1,\gamma)&=
\mathfrak{E}_\alpha^+(\gamma,-u_1^{-1}t_1)\\
&=
\frac{\bigl(bc,q/ad;q\bigr)_{\infty}}
{\bigl(ac,q/bd;q\bigr)_{\infty}}
\frac{G_{\alpha_{\sigma\tau}}(\gamma)}{G_{\beta_{\sigma\tau}}(\gamma)}
\mathfrak{E}_\beta^+(\gamma,-u_1^{-1}t_1)\\
&=\frac{\bigl(bc,q/ad;q\bigr)_{\infty}}
{\bigl(ac,q/bd;q\bigr)_{\infty}}
\frac{G_{\alpha_{\sigma\tau}}(\gamma)}{G_{\beta_{\sigma\tau}}(\gamma)}.
\end{split}
\end{equation*}
Interchanging the role of $t_0$ and $u_1$ (i.e. replacing $\alpha$
by $\alpha_\sigma$) then gives the explicit
evaluation formula
\[
\mathfrak{E}^+(-t_0^{-1}t_1,x)=
\frac{\bigl(qa/d,q/ad;q\bigr)_{\infty}}
{\bigl(ac,c/a;q\bigr)_{\infty}}\frac{\bigl(cx,c/x;q\bigr)_{\infty}}
{\bigl(qx/d,q/dx;q\bigr)_{\infty}}.
\]
Hence we obtain 
\begin{equation}\label{othercase}
\bigl(\mathcal{F}_e^+G_\tau\bigr)(-t_0^{-1}t_1)=
\frac{\bigl(qa/d,q/ad;q\bigr)_{\infty}}
{\bigl(ac,c/a;q\bigr)_{\infty}}\,K,
\end{equation}
with $K$ given by
\begin{equation*}
K=\lim_{\epsilon\searrow 0}\frac{1}{2\pi i}\int_{C_\epsilon}
G_\tau(x)\frac{\bigl(cx,c/x;q\bigr)_{\infty}}
{\bigl(qx/d,q/dx;q\bigr)_{\infty}}W_e^+(x)\frac{dx}{x}.
\end{equation*}
It follows from the explicit
form \eqref{integrandexplicit}
of the integrand $G_{\tau}W_e^+$ 
that $K$ is the Macdonald-Mehta type integral 
$(\mathcal{F}_e^+G_\tau)(s_0)$ in which
the parameter $c$ is replaced by
the parameter $q/d$. In particular, 
Lemma \ref{MM} yields an explicit evaluation for $K$.
Formula \eqref{f}
combined with \eqref{othercase} then shows that
\begin{equation}\label{outcome1}
\begin{split}
f(-t_0^{-1}t_1)&=\frac{1}{\bigl(d/b,qa/c;q\bigr)_{\infty}}
\bigl(\mathcal{F}_e^+G_\tau\bigr)(-t_0^{-1}t_1)\\
&=\frac{\bigl(qa/d,q/ad;q\bigr)_{\infty}}
{\bigl(d/b,ac;q\bigr)_{\infty}\theta(c/a)}\,K\\
&=\frac{2\bigl(q/ad;q\bigr)_{\infty}}
{\bigl(q,ab,ac;q\bigr)_{\infty}\theta(d/b,c/a)}
\frac{\theta(qabe/d)}{\theta(ae,be,qe/d)}.
\end{split}
\end{equation}
On the other hand, \eqref{fgoal} and \eqref{Ce} show that
\begin{equation}\label{outcome2}
\begin{split}
f(-t_0^{-1}t_1)&=C\frac{\theta(-t_0t_1^{-1}\tilde{e},-t_0^{-1}t_1\tilde{e})}
{\theta(d/b,c/a)}\\
&=\frac{2\bigl(q/ad;q\bigr)_{\infty}}
{\bigl(q,ab,ac;q\bigr)_{\infty}\theta(d/b,c/a)}
\frac{\theta(abce,-t_0t_1^{-1}\tilde{e},-t_0^{-1}t_1\tilde{e})}
{\theta(ae,be,ce,t_0t_1\tilde{e},t_0^{-1}t_1^{-1}\tilde{e})}.
\end{split}
\end{equation}
Comparing \eqref{outcome1} and \eqref{outcome2} leads to the identity
\begin{equation}\label{sole}
\frac{\theta(-t_0t_1^{-1}\tilde{e},-t_0^{-1}t_1\tilde{e})}
{\theta(t_0t_1\tilde{e},t_0^{-1}t_1^{-1}\tilde{e})}=
\frac{\theta(qabe/d,ce)}{\theta(abce,qe/d)}.
\end{equation}
In other words, if  
\eqref{fgoal} holds true for the pair $(\tilde{e},C)$, 
then $\tilde{e}$ is necessarily a solution of \eqref{sole}.

It is easy to verify that 
$\tilde{e}:=e_\sigma=-q^{\frac{1}{2}}/u_0t_1e$ is a solution
of \eqref{sole}, as well as $\tilde{e}=e_\sigma^{-1}$. Furthermore,  
the left hand side of \eqref{sole} is an elliptic function in 
$\tilde{e}$ on $T=\mathbb{C}^\times/q^{\mathbb{Z}}$, hence
standard elliptic function theory shows that 
$\bigl(e_{\sigma}q^{\mathbb{Z}}\bigr)^{\pm 1}$ are all possible
solutions for $\tilde{e}$ of equation \eqref{sole}.
Hence we conclude that
$\bigl(e_{\sigma}q^{\mathbb{Z}}\bigr)^{\pm 1}$
are all the possible values for $\tilde{e}$ such that
\eqref{fgoal} holds true for some constant $C$.

For the choice $\tilde{e}=e_\sigma$ in the equation
$(\mathcal{F}_e^+G_\tau)(\gamma)=
C\mathcal{P}_{\tilde{e}}^\sigma(\gamma)^{-1}G_\sigma(\gamma)^{-1}$, 
the corresponding constant $C$ is given by
\[C=\frac{2\bigl(q/ad;q\bigr)_{\infty}}
{\bigl(q,ab,ac;q\bigr)_{\infty}\theta(ae,be,ce,qe/d)}
\]
in view of \eqref{Ce}. This completes the proof of the proposition.
\end{proof}


\subsection{Algebraic inversion and Plancherel formulas}

The results of the previous subsections easily lead to an
algebraic inversion formula for the Askey-Wilson function transform. 
For the formulation
we introduce normalization constants $K_e=K_e^\alpha$
and $K_e^+=K_e^{+,\alpha}$ by
\begin{equation}\label{Ke}
\begin{split}
K_e&=\frac{\bigl(q,qab,ac;q\bigr)_{\infty}}
{\bigl(q/ad;q\bigr)_{\infty}}
\sqrt{\theta(ae,be,ce,qe/d)},\\
K_e^+&=\frac{\bigl(q,ab,ac;q\bigr)_{\infty}}
{2\bigl(q/ad;q\bigr)_{\infty}}
\sqrt{\theta(ae,be,ce,qe/d)},
\end{split}
\end{equation}
where we choose an arbitrary branch of the square root.
Observe that the constants $K_e$ and $K_e^+$ are self-dual, i.e.
$K_{e_\sigma}^\sigma=K_e$ and $K_e^+=K_{e_\sigma}^{+,\sigma}$.
We define normalized difference Fourier transforms 
$\widetilde{\mathcal{F}}_e=\widetilde{\mathcal{F}}_e^\alpha$
and $\widetilde{\mathcal{J}}_e=\widetilde{\mathcal{J}}_e^{\alpha}$
on $M_e$ by
\begin{equation}\label{normAWfunction}
\widetilde{\mathcal{F}}_e=K_e\mathcal{F}_e,\qquad
\widetilde{\mathcal{J}}_e=K_e\mathcal{J}_e.
\end{equation}
Similarly, we define the normalized symmetric Askey-Wilson function
transform $\widetilde{\mathcal{F}}_e^+=
\widetilde{\mathcal{F}}_e^{+,\alpha}$ on $M_e^+$ by
\begin{equation}\label{normsymmAWfunction}
\widetilde{\mathcal{F}}_e^+=K_e^+\mathcal{F}_e^+.
\end{equation}
The algebraic inversion formulas can now be formulated as follows.
\begin{thm}\label{inversionAW} 
{\bf a)}
The normalized Askey-Wilson function transform
$\widetilde{\mathcal{F}}_e$ defines a linear bijection
$\widetilde{\mathcal{F}}_e: M_e\rightarrow M_{e_\sigma}^\sigma$. Its inverse
is given by $\widetilde{\mathcal{J}}_{e_\sigma}^\sigma: M_{e_\sigma}^\sigma
\rightarrow M_e$.\\
{\bf b)} The normalized symmetric Askey-Wilson function transform
$\widetilde{\mathcal{F}}_e^+$ defines a linear bijection
$\widetilde{\mathcal{F}}_e^+: M_e^+\rightarrow
M_{e_\sigma}^{+,\sigma}$. Its inverse is given by
$\widetilde{\mathcal{F}}_{e_\sigma}^{+,\sigma}:
M_{e_\sigma}^{+,\sigma}\rightarrow M_e^+$.
\end{thm}
\begin{proof}
{\bf a)} By Corollary \ref{classicalpart},
\begin{equation}\label{ye1}
\mathcal{F}_e\bigl(\mathcal{P}_e^{-1}G^{-1}\bigr)(\gamma)=
\mathcal{J}_e\bigl(\mathcal{P}_e^{-1}G^{-1}\bigr)(\gamma)
=\frac{\bigl(q/ad;q\bigr)_{\infty}}{\bigl(q,qab,ac;q\bigr)_{\infty}}
G_{\sigma\tau}(\gamma)
\end{equation}
in view of the explicit expression \eqref{Dnul} for $D_0$
and Remark \ref{Gaussianevaluation}.
Furthermore,
by \eqref{samething} and Proposition \ref{FourierGaussian},
\begin{equation}\label{ye2}
\begin{split}
\bigl(\mathcal{F}_e G_\tau\bigr)(\gamma)&=
\bigl(\mathcal{J}_e G_\tau\bigr)(\gamma)\\
&=\frac{\bigl(q/ad;q\bigr)_{\infty}}
{\bigl(q,qab,ac;q\bigr)_{\infty}\theta(ae,be,ce,qe/d)}\,
\mathcal{P}_{e_\sigma}^\sigma(\gamma)^{-1}G_\sigma(\gamma)^{-1}.
\end{split}
\end{equation}
It follows from \eqref{ye1} and \eqref{ye2} that $\mathcal{F}_e$
and $\mathcal{J}_e$ map the cyclic vector of $V_e^{cl}$ (respectively
$V^{str}$) to a multiple of the cyclic vector of $V^{str,\sigma}$
(respectively $V_{e_\sigma}^{cl,\sigma}$). Since $\mathcal{F}_e$
and $\mathcal{J}_e$ are Fourier transforms associated with $\sigma$
and $\sigma_\sigma^{-1}$ respectively, $\mathcal{F}_e$ and
$\mathcal{J}_e$ thus define linear mappings
\begin{equation}\label{domran}
\mathcal{F}_e,\mathcal{J}_e: M_e\rightarrow M_{e_\sigma}^\sigma.
\end{equation}
Observe that the extended parameter map
\[(\alpha,e)\mapsto (\alpha_\sigma,e_\sigma)
\]
is an involution. By the definition \eqref{Ke}
of the constant $K_e$ and by \eqref{ye1}, \eqref{ye2} and \eqref{domran},
we obtain
\begin{equation*}
\begin{split}
\bigl(\mathcal{J}_{e_\sigma}^\sigma\circ\mathcal{F}_e\bigr)
(\mathcal{P}_e^{-1}G^{-1})&=K_e^{-2}
\mathcal{P}_e^{-1}G^{-1}=
\bigl(\mathcal{F}_{e_\sigma}^\sigma\circ\mathcal{J}_e\bigr)
(\mathcal{P}_e^{-1}G^{-1}),\\
\bigl(\mathcal{J}_{e_\sigma}^\sigma\circ\mathcal{F}_e\bigr)(G_\tau)&=
K_e^{-2}G_\tau=
\bigl(\mathcal{F}_{e_\sigma}^\sigma\circ\mathcal{J}_e\bigr)(G_\tau).
\end{split}
\end{equation*}
Since $\mathcal{F}_e$ and $\mathcal{J}_e$ are Fourier transforms
associated with $\sigma$ and $\sigma_\sigma^{-1}$ respectively, and 
$G_\tau, \mathcal{P}_e^{-1}G^{-1}\in M_e$ generate $M_e$
as $\mathcal{H}$-module, we conclude that 
\[\mathcal{J}_{e_\sigma}^\sigma\circ \mathcal{F}_e=
K_e^{-2}\hbox{Id}|_{M_e}=\mathcal{F}_{e_\sigma}^\sigma\circ\mathcal{J}_e.
\]
Combined with \eqref{domran}, this completes the proof of part {\bf a)}.\\
{\bf b)} This follows from {\bf a)} and from Corollary \ref{symmAW}.
\end{proof}

We already noticed that $\widetilde{\mathcal{F}}_e$ and 
$\widetilde{\mathcal{J}}_e$
restrict to linear bijections $\widetilde{\mathcal{F}}_e,
\widetilde{\mathcal{J}}_e:
V_e^{cl}\rightarrow V^{str,\sigma}$, see Corollary \ref{classicalpart}. 
By (the proof of) the above theorem, 
the difference Fourier transforms $\widetilde{\mathcal{F}}_e$ 
and $\widetilde{\mathcal{J}}_e$ also restricts
to linear bijections 
$\widetilde{\mathcal{F}}_e,\widetilde{\mathcal{J}}_e: V^{str}\rightarrow
V_{e_\sigma}^{cl,\sigma}$. In a similar fashion it follows that 
the normalized symmetric Askey-Wilson
function transform $\widetilde{\mathcal{F}}_e^+$
restricts to a linear bijection
$\widetilde{\mathcal{F}}_e^+: 
V^{str,+}\rightarrow V_{e_\sigma}^{cl,+,\sigma}$. 
In the following proposition we 
describe these transforms on a suitable basis of $V^{str}$ (respectively
of $V^{str,+}$).

\begin{prop}\label{strangepart}
{\bf a)} For $s\in\mathcal{S}_\tau$, we have
\begin{equation*}
\begin{split}
\widetilde{\mathcal{F}}_e\bigl(E_\tau(s;\cdot)G_\tau\bigr)(\gamma)&=
\frac{1}{D_0^{\sigma}K_eG_{\tau\sigma}(s)}
E_{\ddagger\sigma\tau}(s^{-1};\gamma^{-1})
\mathcal{P}_{e_\sigma}^\sigma(\gamma)^{-1}G_\sigma(\gamma)^{-1},\\
\widetilde{\mathcal{J}}_e\bigl(E_\tau(s;\cdot)G_\tau\bigr)(\gamma)&=
\frac{1}{D_0^{\sigma}K_eG_{\tau\sigma}(s)}
E_{\sigma\tau}(s;\gamma)\mathcal{P}_{e_\sigma}^\sigma(\gamma)^{-1}
G_\sigma(\gamma)^{-1}
\end{split}
\end{equation*}
with $D_0$ and $K_e$ the explicit constants \eqref{Dnul}
and \eqref{Ke} respectively.\\
{\bf b)} For $s\in\mathcal{S}_\tau^+$, we have
\[
\widetilde{\mathcal{F}}_e^+\bigl(E_\tau^+(s;\cdot)G_\tau\bigr)(\gamma)=
\frac{1}{D_0^{+,\sigma}K_e^+G_{\tau\sigma}(s)}
E_{\sigma\tau}^+(s;\gamma)\mathcal{P}_{e_\sigma}^\sigma(\gamma)^{-1}
G_{\sigma}(\gamma)^{-1},
\]
with $D_0^+$ and $K_e^+$ the explicit constants \eqref{Dnulplus}
and \eqref{Ke} respectively.
\end{prop}
\begin{proof}
We derive the formulas for $\mathcal{J}_e$,
the other cases are derived in a similar fashion.
Let $s\in\mathcal{S}_\tau$. Since $\sigma\tau\sigma=
\tau\sigma\tau$ when acting on the difference
multiplicity function $\alpha$, we have by Corollary \ref{classicalpart}
and Theorem \ref{inversionAW},
\begin{equation*}
\begin{split}
\widetilde{\mathcal{J}}_e\bigl(E_\tau(s;\cdot)G_\tau\bigr)&=
\frac{1}{D_0^{\sigma}K_eG_{\tau\sigma}(s)}
\bigl(\widetilde{\mathcal{J}}_e\circ
\widetilde{\mathcal{F}}_{e_\sigma}^\sigma\bigr)
\bigl(E_{\sigma\tau}(s;\cdot)\mathcal{P}_{e_\sigma}^\sigma{}^{-1}
G_\sigma^{-1}\bigr)\\
&=\frac{1}{D_0^{\sigma}K_eG_{\tau\sigma}(s)}
E_{\sigma\tau}(s;\cdot)\mathcal{P}_{e_\sigma}^{\sigma}{}^{-1}
G_\sigma^{-1},
\end{split}
\end{equation*}
which is the desired result.
\end{proof}
\begin{rem}
Corollary \ref{classicalpart} and Proposition \ref{strangepart}
give an explicit description of the image of a suitable basis
of $M_e$ (respectively $M_e^+$) under the (non)symmetric Askey-Wilson
function transform. For the symmetric Askey-Wilson function 
transform,  a different proof for these formulas 
was derived in \cite[Thm. 5.2]{St3}. 
\end{rem}
We end this subsection with the following
algebraic Plan\-che\-rel type formulas for the 
(non)symmetric Askey-Wilson function transform.
\begin{thm}\label{PlAW}
{\bf a)} For all $g,h\in M_e$,
\[\langle \widetilde{\mathcal{F}}_eg, 
I\bigl(\widetilde{\mathcal{J}}_eh\bigr)\rangle_{e_\sigma}^{\sigma}=
\langle g,Ih\rangle_e,
\]
with $I$ the inversion operator $(If)(\gamma)=f(\gamma^{-1})$.\\
{\bf b)} For all $g,h\in M_e^+$,
\[\langle \widetilde{\mathcal{F}}_e^+g,\widetilde{\mathcal{F}}_e^+h
\rangle_{e_\sigma,+}^{\sigma}
=\langle g,h\rangle_{e,+}.
\]
\end{thm} 
\begin{proof}
{\bf a)} By a formal computation using the duality 
of the normalized rank one Cherednik kernels
$\mathfrak{E}$ and $\mathfrak{E}_\ddagger$  
(see Theorem \ref{dualitytheorem}), we obtain
\begin{equation}\label{tussenstap}
 \langle \widetilde{\mathcal{F}}_eg,Ih\rangle_{e_\sigma}^{\sigma}=
\langle g,I\bigl(\widetilde{\mathcal{F}}_{e_{\sigma}}^{\sigma}h\bigr)
\rangle_e,\qquad
\langle \widetilde{\mathcal{J}}_eg,Ih\rangle_{e_\sigma}^\sigma=
\langle g,I\bigl(\widetilde{\mathcal{J}}_{e_\sigma}^{\sigma}h\bigr)
\rangle_e
\end{equation}
for $g\in M_e$ and $h\in M_{e_\sigma}^{\sigma}$. For a 
rigorous proof of formula
\eqref{tussenstap} we need to justify that
limits and integrations may be interchanged in the formal computation.
In case $g\in V_e^{cl}$, 
this is easily justified by Lemma \ref{Chbounds} and by Fubini's
Theorem, since $g$ vanishes on $(eq^{\mathbb{Z}})^{\pm 1}$. 
In a similar fashion, \eqref{tussenstap} is proven to be correct
when $h\in V_{e_\sigma}^{cl,\sigma}$.

Let now $g\in V^{str}$ and $h\in V^{str,\sigma}$. 
We write $h=\widetilde{\mathcal{J}}_ef$ with
$f\in V_e^{cl}$, then Theorem \ref{inversionAW}
shows that $f=\widetilde{\mathcal{F}}_{e_\sigma}^{\sigma}h$.
We derive that
\begin{equation*}
\begin{split}
\langle \widetilde{\mathcal{F}}_eg, Ih\rangle_{e_\sigma}^\sigma&=
\langle \widetilde{\mathcal{F}}_eg,
I\bigl(\widetilde{\mathcal{J}}_ef\bigr)\rangle_{e_\sigma}^\sigma\\
&=\langle \bigl(\widetilde{\mathcal{J}}_{e_\sigma}^\sigma\circ
\widetilde{\mathcal{F}}_e\bigr)(g), If\rangle_{e}\\
&=\langle g,If\rangle_{e}\\
&=\langle g,I\bigl(\widetilde{\mathcal{F}}_{e_\sigma}^{\sigma}h\bigr)
\rangle_{e}.
\end{split}
\end{equation*}
Here the second equality is allowed since $f\in V_e^{cl}$,
and the third equality follows from Theorem \ref{inversionAW}.
The second identity in \eqref{tussenstap} for $g\in V^{str}$
and $h\in V^{str,\sigma}$ follows by a similar computation.

Formula \eqref{tussenstap} combined with Theorem \ref{inversionAW}
now shows that
\[ \langle \mathcal{F}_eg,I\bigl(\mathcal{J}_eh\bigr)
\rangle_{e_\sigma}^{\sigma}=\langle g,Ih\rangle_e
\]
for $g,h\in M_e$.\\
{\bf b)} This follows by similar arguments as for the nonsymmetric
Askey-Wilson function transform (see part {\bf a)}).
\end{proof}
 

\subsection{The $L^2$-theory}

The results thusfar obtained for the symmetric Askey-Wilson
function transform $\mathcal{F}_e^+$ are sufficient to 
derive ``analytic'' Plancherel and inversion formulas
as follows. 

We need to
restrict the parameter domain first in order to obtain positive
measures. As usual, we assume throughout this subsection
that $0<q^{\frac{1}{2}}<1$.
We assume furthermore that the parameter $e$ and the
difference multiplicity function $\alpha$ satisfy the conditions 
\begin{equation}\label{parconditions}
\begin{split}
&e<0,\qquad\qquad\qquad 0<b,c\leq a<d/q,\\
&bd,cd\geq q,\qquad\qquad ab,ac<1,
\end{split}
\end{equation}
where we used the Askey-Wilson parametrization \eqref{AWparone}
for the difference multiplicity function $\alpha$.
These conditions are invariant under the extended
involution $(\alpha,e)\mapsto (\alpha_\sigma,e_\sigma)$.
For generic parameters satisfying these conditions, we define
a measure $m_e(\cdot)=m_{e}^{\alpha}(\cdot)$ by
\[\int f(x)dm_e(x)=
\frac{1}{2\pi i}\int_{x\in \mathbb{T}}f(x)W_e^+(x)\frac{dx}{x}
+\sum_{x\in D_e}\bigl(f(x)+f(x^{-1})\bigr)
\underset{y=x}{\hbox{Res}}\left(\frac{W_e^+(y)}{y}\right),
\]
where $D_e=D_e^+\cup D_e^-$ is the discrete set
\begin{equation*}
\begin{split}
D_e^+&=\{aq^k \, | \, k\in\mathbb{N}\,:\,\, aq^k>1 \},\\
D_e^-&=\{eq^k \, | \, k\in\mathbb{Z}\,:\,\, eq^k<-1 \}.
\end{split}
\end{equation*}
By continuous extension in the parameters, we obtain a positive
measure $m_{e}^{\alpha}$ for all parameters $\alpha,e$ satisfying
\eqref{parconditions}, see
\cite{KS2} and \cite{St3} for details. The support of 
$m_{e}^{\alpha}$ is $\mathbb{T}\cup D_e\cup D_e^{-1}$.

For the remainder of this subsection we fix parameters $\alpha$ and
$e$ satisfying the conditions \eqref{parconditions}.
Let $L_+^2(m_e)$ be the Hilbert space consisting of $L^2$-functions
$f$ with respect to the measure $m_e$ satisfying $f(x)=f(x^{-1})$
$m_e$-a.e. Clearly, the space $M_e^+$ may be considered
as subspace of $L_+^2(m_e)$. The following result from
\cite[Prop. 6.7]{St3} tells us exactly
when $M_e^+\subseteq L_+^2(m_e)$ is dense.
\begin{prop}\label{density}
Let $k\in\mathbb{Z}$
be the unique integer 
such that $1<eq^k\leq q^{-1}$. Then $M_e^+\subseteq L_+^2(m_e)$ is dense
if and only if $|e_\sigma^{-1}q^k|\geq 1$.
\end{prop} 
\begin{proof}
Using the explicit results on 
the extended, symmetric rank one Macdonald-Koorn\-win\-der transform $F^+$
as derived in Subsection 7.4 (see also \cite{St3}), the proof 
can be reduced to proving the density of some polynomial
space in an explicit $L^2$-space of square integrable functions
with respect to a compactly supported measure.
The proof then follows from standard density results.
For a detailed proof we refer
the reader to \cite[Prop. 6.7]{St3}. 
\end{proof}
As a consequence of Proposition \ref{density}
we arrive at a new proof for the analytic 
inversion and Plancherel formula for the symmetric Askey-Wilson
function transform, see \cite[Thm. 1]{KS2} for the classical approach. 
\begin{thm}
Let $k\in\mathbb{Z}$
be the unique integer 
such that $1<eq^k\leq q^{-1}$ and assume that $|e_\sigma^{-1}q^k|\geq 1$.
The normalized
symmetric Askey-Wilson function transform
$\widetilde{\mathcal{F}}_e^+$ \textup{(}see 
\eqref{normsymmAWfunction}\textup{)}
uniquely extends by continuity
to a surjective isometric isomorphism 
\[\widetilde{\mathcal{F}}_e^+: L_+^2(m_e)\rightarrow 
L_+^2(m_{e_\sigma}^{\sigma}),
\]
with inverse given by
\[\widetilde{\mathcal{F}}_{e_\sigma}^{+,\sigma}:
L_+^2(m_{e_{\sigma}}^{\sigma})\rightarrow L_+^2(m_e).
\]
\end{thm}
\begin{proof}
The conditions
on the parameters $\alpha,e$ 
are invariant under the involution  
$(\alpha,e)\mapsto (\alpha_\sigma,e_\sigma)$.
Hence $L_+^2(m_{e_\sigma}^{\sigma})$ is defined, and
$M_e^+\subseteq L_+^2(m_e)$,  
$M_{e_\sigma}^{+,\sigma}\subseteq L_+^2(m_{e_\sigma}^{\sigma})$ are both
dense by Proposition \ref{density}.
The result now follows immediately from
the algebraic inversion and
Plan\-che\-rel formula for $\widetilde{\mathcal{F}}_e^+$,
see Theorem \ref{inversionAW} and Theorem \ref{PlAW}.
\end{proof}


\section{Appendix}
In this section we prove the bounds for the Macdonald-Koornwinder
polynomials as stated in Proposition \ref{MKbounds}.

Let $\Lambda_i=\{ \lambda\in\Lambda \, | \, r_i\cdot
\lambda\not=\lambda\}$
and denote  
\[\mathcal{S}_i=\{s_{\lambda} \, | \, \lambda\in\Lambda_i\}
\subseteq \mathcal{S}
\]
for $i=0,\ldots,n$. Let
$\mathcal{F}(\mathcal{S}_i)$ be the space of functions
$g:\mathcal{S}_i\rightarrow \mathbb{C}$.
Recall the definition of $c_i:=c_{a_i}\in \mathbb{C}(x)$
for $i=0,\ldots,n$, see \eqref{cf1} and \eqref{cf2}.
We first observe the following elementary fact.
\begin{lem}\label{boundedcoeff}
The functions $|c_i^{\ddagger\sigma}(\cdot)|\in
\mathcal{F}(\mathcal{S})$ and 
$|c_i^{\ddagger\sigma}(\cdot)|^{-1}\in 
\mathcal{F}(\mathcal{S}_i)$
are bounded for $i=0,\ldots,n$.
\end{lem}
\begin{proof}
This follows immediately from 
Lemma \ref{actioncompatibleW}{\bf b)}
and from the explicit expression of the spectral
points $s_\lambda\in\mathcal{S}$ (see Subsection 4.1), 
since the parameters are
assumed to be generic.
\end{proof}
If $\lambda=\sum_i\lambda_i\epsilon_i\in\Lambda$ then we call
$\lambda_i$ the $i$th coordinate of $\lambda$.
For $\lambda\in\Lambda\setminus \{0\}$ we denote
$u_{\lambda}\in W_0$ for the unique element of minimal
length such that the first coordinate of $u_\lambda^{-1}\cdot\lambda$ 
is strictly negative. There are essentially two cases to
consider here; if some of the coordinates of $\lambda$ are
strictly negative, then $u_\lambda=r_{i-1}r_{i-2}\cdots r_1$
with $i\geq 1$ the smallest such that $\lambda_i<0$ (and
$u_\lambda=1\in W_0$ the identity element when $i=1$).
If all coordinates of $\lambda$ are nonnegative and
$\lambda_j>0$, $\lambda_{j+1}=\cdots=\lambda_n=0$, then
$u_\lambda=r_j\cdots r_{n-1}r_nr_{n-1}\cdots r_1$. 
Finally, for $0\in \Lambda$ the zero element, we denote
$u_0=r_1\cdots r_{n-1}r_nr_{n-1}\cdots r_1\in W_0$. Observe that
\[T_{u_0}=T_1\cdots T_{n-1}T_nT_{n-1}\cdots T_1=U\in\mathcal{H}
\]
(see \eqref{U}), since $u_0=r_1\cdots r_{n-1}r_nr_{n-1}\cdots r_1$
is a reduced expression in $W_0$.
We prove the following refinement of Proposition \ref{MKbounds}.
\begin{lem}\label{tussenstapp}
Let $K\subset \bigl(\mathbb{C}^\times\bigr)^n$ be a compact set.
Then there exist positive constants $C_1$ \textup{(}independent of 
$K$\textup{)} and $C_2$ \textup{(}dependent of $K$\textup{)}
such that 
\[ |E(s_{\lambda};x)|\leq C_1^{\,l(u_\lambda)}C_2^{N(\lambda)},\qquad
\forall x\in K,\,\,\, \forall \lambda\in\Lambda.
\]
\end{lem}
Since $l(u_{\lambda})\leq 2n-1$ for all $\lambda\in\Lambda$ and
$E(s_0;\cdot)=1$, 
Proposition \ref{MKbounds} is a direct consequence of Lemma
\ref{tussenstapp}.
We start with a preliminary lemma.
\begin{lem}\label{recu}
For all $\mu\in\Lambda$,
\begin{equation*}
u_nc_0^{\ddagger\sigma}(s_\mu)E(s_{r_0\cdot\mu};x)
=u_n\bigl(c_0^{\ddagger\sigma}(s_{\mu})-1\bigr)E(s_{\mu};x)+
x_1^{-1}\bigl(U_{\ddagger\sigma}\cdot E(\cdot\,;x)\bigr)(s_\mu).
\end{equation*} 
\end{lem}
\begin{proof}
By \eqref{Tdot} and the definition of $T_0$,
\begin{equation}\label{F2}
 \bigl(T_0^{\ddagger\sigma}E(\cdot\,;x)\bigr)(s_\mu)=
\bigl(u_n^{-1}-u_nc_0^{\ddagger\sigma}(s_{\mu})\bigr)E(s_{\mu};x)+
u_nc_0^{\ddagger\sigma}(s_{\mu})E(s_{r_0\cdot\mu};x).
\end{equation}
On the other hand, 
\begin{equation*}
\begin{split}
\bigl(T_0^{\ddagger\sigma}\cdot E(\cdot\,;x)\bigr)(s_\mu)&=
\bigl(T_0^{\ddagger\sigma}\cdot
I\mathfrak{E}_{\mathcal{A}}(\cdot,x)\bigr)(s_\mu)\\
&=\bigl(T_0^{\sigma}{}^{-1}\cdot
\mathfrak{E}_{\mathcal{A}}(\cdot,x)\bigr)(s_\mu^{-1})\\
&=\bigl(\psi_\sigma(T_0^{\sigma}{}^{-1})E(s_\mu;\cdot)\bigr)(x)
\end{split}
\end{equation*}
by Proposition \ref{allrelations}, where $I:
\mathcal{F}(\mathcal{S}_\ddagger)\rightarrow 
\mathcal{F}(\mathcal{S})$
is the ``inversion map'' $(Ig)(s)=g(s^{-1})$.
Now $\psi_\sigma(T_0^{\sigma}{}^{-1})=Ux_1=u_n^{-1}-u_n+x_1^{-1}U^{-1}$,
see \cite{Sa} or \cite[Appendix]{St2} for the second equality,
and $(\dagger_\sigma\circ\psi)(U^{-1})=U_{\ddagger\sigma}$, hence 
\begin{equation}\label{F1}
\bigl(T_0^{\ddagger\sigma}\cdot E(\cdot\,;x)\bigr)(s_\mu)=
(u_n^{-1}-u_n)E(s_{\mu};x)+
x_1^{-1}\bigl(U_{\ddagger\sigma}\cdot E(\cdot\,;x)\bigr)(s_\mu).
\end{equation}
Subtracting \eqref{F2} from \eqref{F1} gives the result.
\end{proof}
In a similar fashion as in the proof of Lemma \ref{recu}
we obtain from Proposition \ref{allrelations} the relations
\begin{equation}\label{Titransform}
 \bigl(T_wE(s;\cdot)\bigr)(x)=\bigl(T_w^{\ddagger\sigma}{}^{-1}
\cdot E(\cdot;x)\bigr)(s)
\end{equation}
for $w\in W_0$ and $s\in\mathcal{S}$.

Let $\lambda\in\Lambda\setminus\{0\}$ and set
$\mu=(r_0u_{\lambda}^{-1})\cdot\lambda$. Observe
that $N(\mu)=N(\lambda)-1$, as follows easily from the definition
of the action $\mathcal{W}$ on $\Lambda$ and from the definition
of $u_\lambda\in W_0$. Lemma \ref{recu} applied to $\mu$ now gives
the recurrence relation
\begin{equation}\label{startingpoint}
\begin{split}
u_nc_0^{\ddagger\sigma}(s_{(r_0u_{\lambda}^{-1})\cdot\lambda})
E(s_{u_{\lambda}^{-1}\cdot\lambda};x)&=
x_1^{-1}\bigl(U_{\ddagger\sigma}\cdot
E(\cdot\,;x)\bigr)
(s_{(r_0u_{\lambda}^{-1})\cdot\lambda})\\
+u_n&\bigl(c_0^{\ddagger\sigma}(s_{(r_0u_{\lambda}^{-1})\cdot\lambda})-
1\bigr)E(s_{(r_0u_{\lambda}^{-1})\cdot\lambda};x).
\end{split}
\end{equation}
Let $\preceq$ be the Bruhat order on $W_0$.
We apply $T_{u_{\lambda}}$, acting upon $\mathcal{A}$
on the $x$-variable,
to both sides of \eqref{startingpoint}. For the left hand
side, we obtain in view of \eqref{Titransform} an expression
of the form
\[
K(\lambda)E(s_{\lambda};x)+\sum_{\stackrel{\scriptstyle{w\in W_0}}
{w\prec u_\lambda^{-1}}}e_w(\lambda)
E(s_{(u_{\lambda}w)^{-1}\cdot\lambda};x).
\]
The coefficient $K(\lambda)$ can be computed explicitly as follows
(compare with the proof of Lemma \ref{Fcyclic}). Let $u_\lambda=
r_{i_1}r_{i_2}\cdots r_{i_l}$ be a reduced expression for
$u_\lambda\in W_0$, so that $l=l(u_\lambda)$, and define
\[\lambda_{j+1}=(r_{i_j}\cdots r_{i_2}r_{i_1})\cdot\lambda,
\qquad j=0,\ldots,l,
\]
with the convention that $\lambda_1=\lambda$. Observe that
$\lambda_{l+1}=u_\lambda^{-1}\cdot\lambda$, and that 
$\lambda_j=r_{i_j}\cdot\lambda_{j+1}\not=\lambda_{j+1}$
for $j=1,\ldots,n$, cf. the proof of Lemma \ref{Fcyclic}.
Then
\begin{equation}\label{Klambda}
K(\lambda)=
u_nt_{u_{\lambda}}c_0^{\ddagger\sigma}
(s_{(r_0u_{\lambda}^{-1})\cdot\lambda})
\prod_{j=1}^{l(u_\lambda)}
c_{i_j}^{\ddagger\sigma}(s_{\lambda_{j+1}})
\end{equation}
and $K(\lambda)\not=0$ due to Lemma \ref{actioncompatibleW}{\bf b)}.
Furthermore,  $|K(\lambda)|^{-1}$ is uniformly
bounded for $\lambda\in\Lambda\setminus \{0\}$ in view of Lemma
\ref{boundedcoeff}. Similarly, $|e_w(\lambda)|$ is uniformly bounded
for $\lambda\in\Lambda\setminus \{0\}$ and $w\in W_0$
with $w\prec u_{\lambda}^{-1}$.
Finally, note that $l(u_\mu)<l(u_\lambda)$
for $\mu=(u_{\lambda}w)^{-1}\cdot\lambda$
with $\lambda\in\Lambda\setminus\{0\}$,
$w\in W_0$ and $w\prec u_{\lambda}^{-1}$.
Hence we conclude that acting by $T_{u_\lambda}$ on the
left hand side of \eqref{startingpoint} yields an expression
of the form
\begin{equation}\label{expression1}
K(\lambda)E(s_\lambda;x)+
\sum_{\stackrel{\scriptstyle{\mu\in W_0\cdot\lambda}}
{l(u_\mu)<l(u_\lambda)}}
c(\lambda;\mu)E(s_\mu;x)
\end{equation}
with $K(\lambda)$ given by \eqref{Klambda} and 
\begin{equation}\label{bound1}
|K(\lambda)|^{-1}\leq L, \qquad |c(\lambda,\mu)|\leq L
\end{equation}
for all $\lambda\in\Lambda\setminus\{0\}$ and 
$\mu\in W_0\cdot\lambda$ such that $l(u_\mu)<l(u_\lambda)$,
for some constant $L>0$.

To deal with the action of 
$T_{u_{\lambda}}$ on the right hand side 
of \eqref{startingpoint}, we need the commutation relations
between the $T_i$ and $x_j$ in $\mathcal{H}$ for $i,j=1,\ldots,n$.
They are given by
\begin{equation*}
\begin{split}
&T_ix_j=x_jT_i,\qquad\qquad\qquad\quad\qquad |i-j|>1,\\
&T_ix_{i-1}=x_{i-1}T_i,\qquad\quad\qquad\qquad i=2,\ldots,n,\\
&T_ix_iT_i=x_{i+1},\qquad\qquad\qquad\qquad i=1,\ldots,n-1,\\
&(x_n^{-1}T_n^{-1}-u_n)(x_n^{-1}T_n^{-1}+u_n^{-1})=0,
\end{split}
\end{equation*}
see e.g. \cite{Sa} or \cite[Prop. 6.5]{St2}. 
In view of \eqref{Titransform}, applying  
$T_{u_{\lambda}}$ to the first term in the right hand side 
of \eqref{startingpoint} then yields an expression of the form
\begin{equation}\label{expression2}
\bigl(X\cdot E(\cdot;x)\bigr)(s_{(r_0u_\lambda^{-1})\cdot\lambda})+
\sum_{\stackrel{\scriptstyle{1\leq j\leq n}}
{\xi=\pm 1}}x_j^{\xi}\bigl(X_{j,\xi}\cdot
E(\cdot\,;x)\bigr)(s_{(r_0u_\lambda^{-1})\cdot\lambda})
\end{equation}
with $X,X_{j,\xi}\in H_0^{\ddagger\sigma}$ independent
of $\lambda$, where (recall) $H_0^{\ddagger\sigma}
\subset \mathcal{H}_{\ddagger\sigma}$
is the subalgebra generated by $T_i^{\ddagger\sigma}$
($i=1,\ldots,n$).
Finally, applying $T_{u_{\lambda}}$ to the second 
term in the right hand side 
of \eqref{startingpoint},  yields 
\begin{equation}\label{expression3}
u_n\bigl(c_0^{\ddagger\sigma}(s_{(r_0u_{\lambda}^{-1})\cdot\lambda})-1
\bigr)
\bigl(T_{u_{\lambda}}^{\ddagger\sigma}{}^{-1}\cdot
E(\cdot\,;x)\bigr)(s_{(r_0u_{\lambda}^{-1})\cdot\lambda}).
\end{equation}
Combining \eqref{expression1}, \eqref{expression2} and 
\eqref{expression3}, applying $T_{u_{\lambda}}$
to \eqref{startingpoint} gives the identity
\begin{equation*}
\begin{split}
&K(\lambda)E(s_{\lambda};x)=-\sum_{
\stackrel{\scriptstyle{\mu\in W_0\cdot\lambda}}{l(u_\mu)<l(u_\lambda)}}
c(\lambda;\mu)E(s_\mu;x)+
\sum_{\stackrel{\scriptstyle{1\leq j\leq n}}{\xi=\pm 1}}x_j^\xi\bigl(
X_{j,\xi}\cdot E(\cdot\,;x)\bigr)(s_{(r_0u_{\lambda}^{-1})\cdot\lambda})\\
&\quad +\bigl(X\cdot E(\cdot;x)\bigr)(s_{(r_0u_\lambda^{-1})\cdot\lambda})
+u_n\bigl(c_0^{\ddagger\sigma}(s_{(r_0u_{\lambda}^{-1})\cdot\lambda})-1
\bigr)
\bigl(T_{u_{\lambda}}^{\ddagger\sigma}{}^{-1}\cdot
E(\cdot\,;x)\bigr)(s_{(r_0u_{\lambda}^{-1})\cdot\lambda}).
\end{split}
\end{equation*}
This formula will now be used inductively. 
Observe that 
the terms 
\[ \bigl(Z\cdot E(\cdot\,;x)\bigr)
(s_{(r_0u_{\lambda}^{-1})\cdot\lambda}),\qquad\qquad
\bigl(T_{u_{\lambda}}^{\ddagger\sigma}{}^{-1}\cdot
E(\cdot\,;x)\bigr)(s_{(r_0u_{\lambda}^{-1})\cdot\lambda})
\]
with $Z=X,X_{j,\xi}$
are linear combinations of $E(s_{\mu};x)$ with 
$\mu\in W_0\cdot((r_0u_\lambda^{-1})\cdot\lambda)$
and with coefficients which are uniformly bounded for all
$\lambda\in\Lambda\setminus\{0\}$ and
$\mu\in W_0\cdot((r_0u_\lambda^{-1})\cdot\lambda)$, in view of
Lemma \ref{boundedcoeff}. Combined with the uniform bounds
\eqref{bound1}, we conclude that there exists a constant $M\geq 1$
such that 
\begin{equation}\label{inductionstartingpoint}
\begin{split}
|E(s_{\lambda};x)|\leq 
M&\sum_{\stackrel{\scriptstyle{\mu\in W_0\cdot\lambda}}
{l(u_{\mu})<l(u_{\lambda})}}
|E(s_{\mu};x)|\\
&+M\Bigl(1+2n\,\underset{j,\xi}{\hbox{max}}(|x_j^{\xi}|)
\Bigr)\sum_{\nu\in 
W_0\cdot ((r_0u_{\lambda}^{-1})\cdot\lambda)}|E(s_{\nu};x)|
\end{split}
\end{equation}
for all
$\lambda\in\Lambda\setminus \{0\}$, where the maximum is taken over
$j\in\{1,\ldots,n\}$ and $\xi\in\{\pm 1\}$. Now fix a compactum
$K\subset \bigl(\mathbb{C}^\times\bigr)^n$, and choose a constant
$C_1\geq 2|W_0|M$ and a constant $C_2\geq 1$ (depending on $K$), such
that
\begin{equation}\label{condC}
 M\bigl(1+2n\,\underset{x,j,\xi}{\hbox{max}}(|x_j^{\xi}|)\bigr)|W_0|
C_1^{2n-1}\leq \frac{C_2}{2},
\end{equation}
where the maximum is taken over $x\in K$, $j\in\{1,\ldots,n\}$
and $\xi\in\{\pm 1\}$.
Fix $\lambda\in\Lambda\setminus \{0\}$, and assume that
\[ |E(s_{\mu};x)|\leq C_1^{l(u_{\mu})}C_2^{N(\mu)}, \qquad \forall\,
x\in K
\]
when $\mu\in\Lambda$ satisfies either $N(\mu)<N(\lambda)$,
or $N(\mu)=N(\lambda)$ and $l(u_{\mu})<l(u_{\lambda})$. 
We use now the fact that 
$N(\nu)=N(w\cdot\nu)$ for all $w\in W_0$ and $\nu\in\Lambda$,
and that $N((r_0u_\lambda^{-1})\cdot\lambda)=N(\lambda)-1$
for $\lambda\in\Lambda\setminus\{0\}$, to derive from 
\eqref{inductionstartingpoint} that
\begin{equation*}
\begin{split}
\underset{x\in K}{\hbox{max}}|E(s_{\lambda};x)|\leq &MC_2^{N(\lambda)}
\sum_{\stackrel{\scriptstyle{\mu\in W_0\cdot\lambda}}
{l(u_{\mu})<l(u_{\lambda})}}
C_1^{l(u_{\mu})}\\
+&M\bigl(1+2n\,\underset{x,j,\xi}{\hbox{max}}(|x_j^{\xi}|)\bigr)
C_2^{N(\lambda)-1}
\sum_{\mu\in W_0\cdot((r_0u_\lambda^{-1})\cdot\lambda)}C_1^{l(u_\mu)}\\
\leq
&M|W_0|\Bigl(C_2^{N(\lambda)}C_1^{l(u_\lambda)-1}
+\bigl(1+2n\,\underset{x,j,\xi}{\hbox{max}}(|x_j^\xi|)\bigr)
C_2^{N(\mu)-1}C_1^{2n-1}\Bigr)\\
\leq &\frac{1}{2}C_2^{N(\lambda)}C_1^{l(u_\lambda)}+
\frac{1}{2}C_2^{N(\mu)}\\
\leq &C_2^{N(\lambda)}C_1^{l(u_\lambda)},
\end{split}
\end{equation*}
where we used
$l(u_\mu)\leq 2n-1$ for $\mu\in\Lambda$ in the second equality,
the estimates \eqref{condC} and 
$C_1\geq 2|W_0|M$ in the third 
equality, and the fact that $C_1\geq 1$. This proves
Lemma \ref{tussenstapp}, and hence Proposition \ref{MKbounds}.


\bibliographystyle{amsplain}

\end{document}